\documentclass[11pt]{amsart}

\usepackage{amssymb,amsfonts}
\usepackage{amsmath}
\usepackage{graphicx}
\usepackage[all,arc]{xy}
\usepackage{enumerate}
\usepackage{hyperref}
\usepackage{mathrsfs}
\usepackage{tikz-cd}
\usetikzlibrary{calc}

\usetikzlibrary{decorations.pathreplacing}

\usepackage{mathtools}
\usepackage{ytableau}
\usepackage[left=1in,top=1in,right=1in,bottom=1in]{geometry}
\DeclarePairedDelimiter{\ceil}{\lceil}{\rceil}

\newtheorem{thm}{Theorem}[section]
\newtheorem*{thm*}{Theorem}

\newtheorem{cor}[thm]{Corollary}
\newtheorem{prop}[thm]{Proposition}
\newtheorem{lem}[thm]{Lemma}

\newtheorem{conj}[thm]{Conjecture}

\theoremstyle{definition}
\newtheorem{defn}[thm]{Definition}

\newtheorem{con}[thm]{Construction}
\newtheorem{exmp}[thm]{Example}

\newtheorem{rem}[thm]{Remark}

\theoremstyle{remark}

\newtheorem{ansatz}[thm]{Ansatz}

\newcommand{\Z}{\mathbb{Z}}

\newcommand{\C}{\mathbb{C}}

\newcommand{\A}{\mathcal{A}}

 \newcommand{\X}{\mathcal{X}}

\newcommand{\colim}{\mathrm{colim}}
\newcommand{\Coh}{\mathrm{Coh}}
\newcommand{\Fuk}{\mathrm{Fuk}}

\renewcommand{\emptyset}{\varnothing}

\makeatletter
\let\c@equation\c@thm
\makeatother
\numberwithin{equation}{section}

\title{Functorial Homological Mirror Symmetry for Acyclic Cluster Varieties}
\author{Josephine Hlavinka}
\address{Josephine Hlavinka, Department of Mathematics, University of California, Berkeley,
         970 Evans Hall, Berkeley, CA 94720-3840, USA}
\email{josephhlavinka@berkeley.edu}
\begin{document}
\begin{abstract}
We introduce the formalism of \textit{cluster models}, a class of log Calabi–Yau varieties obtained from smooth quasi-affine toric varieties by iterated affine blowups along non-toric codimension 2 centers determined by a really full rank quiver. Cluster models include the truncated cluster varieties of Gross, Hacking, and Keel with respect to really full rank quivers, all really full rank acyclic cluster varieties, and all open toric subvarieties of $\mathbb C^n\times(\mathbb C^*)^m$. We then construct mirrors to cluster models via iterated Weinstein handle attachments to toric mirrors, so that the corresponding skeleta acquire compact pieces, and we prove homological mirror symmetry for all cluster models, compatibly with freezing, deleting, and contracting vertices of the quiver. This extends toric mirror symmetry and its naturalities beyond the cotangent bundle. In particular, we obtain homological mirror symmetry for all really full rank acyclic cluster varieties. We develop this package with a view towards conjectures of Marsh, Rietsch, and Williams about mirror symmetry for Grassmannians and closed positroid varieties. Motivated by these, we also conjecture that open positroid varieties are self-mirror. Assuming a natural ansatz, we prove this for a class of positroid varieties that includes nontrivial cells in every Grassmannian and all cells in finite-type Grassmannians. Towards this end we introduce a new combinatorial input, the notion of weak $k$-booleanness for permutations and braids.
    
    %%Many spaces of interest in homological mirror symmetry arise as compactifications of really full rank locally acyclic cluster varieties. Building on work of Zhou et al. for $A_n$ quivers, we exhibit a homological mirror symmetry (HMS) equivalence for all really full rank \textit{acyclic} cluster varieties. We also posit a non-compact version of Marsh, Rietsch, and Williams' homological mirror symmetry conjectures for Grassmannians and closed positroid varieties, and (in tandem with an ansatz) we apply our main theorem to establish many instances of this ``open" conjecture. In particular we deduce a ``self-mirrorness" theorem for every cell in the positroid stratification of a finite cluster type Grassmannian, i.e. all $Gr(2,n)$'s, $Gr(3,6)$, $Gr(3,7)$, and $Gr(3,8)$.
    
    %%The key technical tool is a generalization of Gross, Hacking, and Keel's construction of truncated cluster varieties as non-toric blowups of toric varieties and a description of the mirrors to varieties which arise from it, which we call \textit{partially truncated cluster varieties}. We expect these varieties (examples of which include all open toric subvarieties of $\mathbb{C}^n$, truncated cluster varieties, acyclic cluster varieties, and more) to be of independent interest, and establish HMS for all of them.
\end{abstract}
\maketitle
\tableofcontents
\bigskip

\section{Introduction}
Homological mirror symmetry for toric varieties identifies the wrapped Fukaya category of the cotangent bundle of a torus $T^*T^n$, stopped at the boundary at infinity $\partial^\infty (-\mathbb L_\Sigma)$ of the FLTZ Lagrangian $-\mathbb{L}_\Sigma$, with coherent sheaves on the toric variety $X_\Sigma$, compatibly with the combinatorics of the fan $\Sigma$. In this paper we show that this package survives leaving the cotangent bundle. That is, we study a symplectic formalism which allows the skeleton of a toric mirror to acquire \textit{compact} pieces (via attaching handles whose cores are disks times skeleta of smaller mirrors), and we show that on the mirror side this corresponds to affine blowups of $X_\Sigma$ along non-toric codimension $2$ centers. When the centers of said blowups are defined by the combinatorics of a really full rank quiver $Q$, we call the resulting log Calabi–Yau varieties $\mathcal A(Q)_{\Sigma,A}$ \textit{cluster models}, and for each such variety we construct an explicit stopped Weinstein manifold $\mathcal X(Q)_{\Sigma,A}$ satisfying $$\Fuk(\mathcal X)\simeq \Coh(\mathcal A)$$ \noindent compatibly with freezing, deleting, and contracting vertices of $Q$ (Theorem \ref{natural CM HMS}). Cluster models include all really full rank acyclic cluster varieties, and so these acquire explicit mirrors (Theorem \ref{Acyclic HMS}). In fact we expect $\mathcal X(Q)$ to be Weinstein deformation equivalent to $\mathcal A(Q)$ itself, so that these cluster varieties are ``self-mirror," for reasons which we motivate in the remainder of the introduction via the theory of positroid cells. Applying the theory to positroid cells, we obtain a homological ``self-mirror symmetry" statement for a class of positroid cells including nontrivial cells in every $\mathrm{Gr}(k,n)$, and all cells in finite-type Grassmannians.

\subsection{Scientific Context: Cluster HMS as Toric HMS and Open HMS for Positroid Varieties}

This work was incited by investigations into the geometry of cells in the positroid stratification of the Grassmannian $Gr(k,n)$, also known as (open) \textit{positroid varieties} (as defined in \cite{KLS}). Though their origins are ostensibly in combinatorial and Schubert geometry, positroid varieties are remarkable in that they also arise as moduli spaces of both $A$-branes, namely constructible sheaves with an interesting singular support condition on $\mathbb{R}^2$ (see \cite{STWZ}), and of closed strings, namely stable map spaces of rational curves into $Gr(k,n)$ (see Section 8 of \cite{KLS}). They are smooth, affine, and Calabi-Yau: indeed, the pairs $(\overline{\Pi}_{u,w}, \overline{\Pi}_{u,w}\setminus\Pi_{u,w})$ consisting of their closures in $Gr(k,n)$ and the boundary are log Calabi–Yau compactifications.\\

This furnishes the idea that $Gr(k,n)$, and more generally any closure of a positroid variety, may participate in a homological mirror pairing with a Landau--Ginzburg model with underlying space the original positroid variety. Formalizing this, beautiful work of Marsh, Rietsch, and Williams (\cite{MR} and \cite{RW}) suggests the following

\begin{conj}[Marsh--Rietsch, Rietsch--Williams; the MRW conjecture]\label{MRW}
Given any positroid variety $\Pi_{u,w}$, there exists a superpotential $W: \Pi_{u,w} \rightarrow \mathbb{C}$ and a small desingularization $\widehat{\overline{\Pi}}_{u,w}$ of the closed variety $\overline{\Pi}_{u,w}$ such that
$$\Fuk(\widehat{\overline{\Pi}}_{u,w}) \simeq MF(\Pi_{u,w}, W)\:\:\text{and}\:\:\Coh(\widehat{\overline{\Pi}}_{u,w}) \simeq FS(\Pi_{u,w}, W).$$
\end{conj}

\noindent Here the pair $(\Pi_{u,w}, W)$ forms a \textit{Landau--Ginzburg (LG) model}, and $FS$ and $MF$ refer to the Fukaya--Seidel and matrix factorization categories, respectively, of $(\Pi_{u,w}, W)$.\\

There has been little progress on Conjecture \ref{MRW}, and indeed a precise statement of it for all $\Pi_{u,w}$ hasn't yet appeared in print. The only solutions we are aware of for this problem (besides when $\overline{\Pi}_{u,w}$ is a toric variety) are in the case $\overline{\Pi}_{u,w} = Gr(k,n)$ with $n$ prime, in which fantastic work by Castronovo establishes the equation on the left above \cite{castronovo2020fukaya}.\\

We begin by recalling the full toric HMS package proven by Bondal, Feng-Liu-Treumann-Zaslow (FLTZ), Kuwagaki, Gammage-Shende, Zhou, and others (see \cite{FLTZ}, \cite{Kuwagaki}, \cite{gammage2022mirror}, and \cite{zhou2019twisted}). Given a manifold $X$, we write $Sh_\Lambda(X)$ to denote the category of constructible sheaves on $X$ with singular support along a conic Lagrangian $\Lambda \subset T^*X$, and write $(-)^c$ to denote the compact objects in a large category. Then, also writing $X_\Sigma$ to denote the toric variety corresponding to an $n$-dimensional fan $\Sigma$, we have the following

\begin{thm}[Toric HMS and its Naturality]\label{Toric HMS} There exists a singular support condition $-\mathbb{L}_\Sigma \subset T^*(S^1)^n$ defined explicitly from $\Sigma$ such that
$$(Sh_{-\mathbb{L}_\Sigma}((S^1)^n))^c \simeq \Coh(X_\Sigma).$$
\end{thm}

Further, this equivalence is natural under combinatorial operations on $\Sigma$, in that stop removal is mirror to restriction onto open toric subvarieties (which corresponds combinatorially to passage to a subfan of $\Sigma$), and microlocalization is mirror to restriction onto closed toric subvarieties (which corresponds combinatorially to taking stars of cones in $\Sigma$). We refer the reader to Theorem \ref{toric HMS natural} for a more complete enumeration of the naturality properties that toric HMS enjoys. It is our view that any HMS statement which takes combinatorial objects as input should supply such naturalities, and a large portion of this work is dedicated to doing so for cluster models.\\

Notably for the present argument, FLTZ's proposal crucially relied upon the simpler, singular-support-less fact that a class of open toric varieties was self-mirror, namely that $$Sh_{0_{(S^1)^n}}((S^1)^n) \simeq \Fuk((\mathbb{C}^*)^n) \simeq \Coh((\mathbb{C}^*)^n).$$

Reasoning naively, one can imagine that a closed positroid variety is an object which is quite like a toric variety: after all, a toric variety is a combinatorially defined partial log CY compactification of $(\mathbb{C}^*)^n$, and a closed positroid variety $\overline{\Pi}_{u,w}$ is a combinatorially defined log CY compactification of $\Pi_{u,w}$. To approach Conjecture \ref{MRW}, then, one might imagine it should behave something like HMS for toric varieties, and then decide to work up to it by first proving HMS for the corresponding ``tori": that is, for the $\Pi_{u,w}$.\\

In fact the form of Conjecture \ref{MRW} suggests what the mirror should be even in this ``open" setting. Indeed, in all known cases, the components of $\overline\Pi_{u,w}\setminus\Pi_{u,w}$ biject with the summands of $W$. De-compactifying on both sides, this means that it is natural to expect $\Pi_{u,w}$ to be its \textit{own mirror} when equipped with the Stein structure arising from it being a smooth affine variety. Therefore the discussion above leads us to propose the following

\begin{conj}[Open MRW's Conjecture]\label{Open MRW}
    Let $\Pi_{u,w}$ be any positroid variety. Then $\Pi_{u,w}$ is self-mirror: i.e., $$\Fuk(\Pi_{u,w}) \simeq \Coh(\Pi_{u,w})$$ \noindent where $\Pi_{u,w}$ carries the Stein structure arising from it being a smooth affine variety.
\end{conj}

Two questions arise: whether positroid varieties are torus-like enough for this to be reasonable, and how one might prove it. Recent work by \cite{GaLa} answers both, by showing that positroid varieties are all \textit{cluster varieties}. On the one hand, cluster varieties contain a rich collection of dense cluster tori (so that the $\Pi_{u,w}$ are quite torus-like indeed). On the other, by work of Gross, Hacking, and Keel \cite{GHK}, cluster varieties arise, up to codimension 2, from toric varieties by affine blowups along non-toric codimension 2 centers. This is to say that cluster varieties differ from toric ones only in some simple local modification. The fundamental idea of the present paper is that, taking the analogy with toric geometry above seriously, Conjecture \ref{Open MRW} should therefore be theoretically downstream of a version of toric HMS for such spaces. We develop such a theory here, and as its first application we prove HMS for all really full rank acyclic cluster varieties.\\

This picture is borne out in two known cases. First, work of Gammage and Le \cite{GL} establishes a duality between handle attachments (on the $A$-side) and affine blowups (on the $B$-side) which yields HMS for all of \cite{GHK}'s \textit{truncated cluster varieties}. Gammage and Le prove and leverage the interesting fact that the cluster variety $\mathcal{A}(A_1)$ of the $A_1$ quiver, which is isomorphic to the dense open positroid variety in $Gr(2,4)$ up to torus factors, is indeed self-mirror: $\Fuk(\mathcal{A}(A_1)) \simeq \Coh(\mathcal{A}(A_1))$ for a suitable Weinstein structure.\\

Our key cluster theoretic forerunner is recent work of \cite{SSWZ}, which proposes a symplectic formalism matching cluster varieties with skeleta built by compact cell attachments, and carries it out for all $A_n$ cluster varieties (equivalently, all $\Pi_{2,n}$’s) in the form expected from Conjecture \ref{Open MRW}: they show that $\Fuk(\Pi_{2,n}) \simeq \Coh(\Pi_{2,n})$ (subsequent to appropriate choices) for all $n$. Our construction of skeleta in Section 3 is heavily inspired by theirs, but \cite{SSWZ} develops neither the $B$-side systematically nor the naturalities of Theorem 1.2.\\

While the preceding two works are both cases of HMS for positroid varieties, we see no natural reason to restrict the scope of Conjecture \ref{Open MRW} to this setting: indeed, all positroid varieties are in particular \textit{really full rank, leaf-recurrent} cluster varieties, and we expect the corresponding statement to be true in that generality.

\begin{conj}[Self-Mirrorness for Leaf-Recurrent Cluster Varieties]\label{leaf rec HMS}
    Let $\mathcal{A}(Q)$ be the cluster variety of a really full rank, leaf-recurrent quiver $Q$. Then $\mathcal{A}(Q)$ is self-mirror.
\end{conj}

We expect an even broader extension of this conjecture, namely to the setting of sink-recurrent cluster varieties, though we also expect that in non-simply-laced types that self-mirrorness will fail (given the concomitant failure of $\mathcal{A}$ vs. $\mathcal{X}$ cluster variety auto-duality in this setting). For this reason a statement in this generality is outside of the scope of the present work.\\

From our perspective, the works above show how one can push toric HMS into a non-cotangent setting: indeed, the mirror of $\mathcal A(A_1)$ has a compact disk in its skeleton (see Figure \ref{A1 skeleton} below), which no conic Lagrangian in $T^*T^2$ can produce. What is missing is a theory which does this systematically on both sides, while retaining the naturalities of Theorem \ref{Toric HMS}. The purpose of the next subsection is to supply one.

\subsection{Results: Toric HMS Beyond the Cotangent Bundle}

On the $B$-side, our basic observation is that a modification of Gross, Hacking, and Keel's construction of truncated cluster varieties as affine blowups suffices to produce acyclic RFR cluster varieties \textit{verbatim}. In Subsection 2.2 we show that by varying their construction further, one can in fact produce a class of partial log Calabi–Yau compactifications of (truncated) cluster varieties of RFR quivers which we name \textit{cluster models}. Cluster models, which depend on a triple $(Q, \Sigma, A)$ consisting of an RFR quiver, a fan, and a set of rays of $\Sigma$ called a \textit{cluster datum}, arise as iterated affine blowups of smooth quasi-affine toric varieties $X_\Sigma$ along (strict transforms of) smooth, non-toric codimension $2$ subvarieties $H_i$ embedded in toric divisors $D_i$. They are all quasi-affine and have trivial Chow groups (see Proposition \ref{chow vanishing}), and one can recover $\mathcal{A}(Q)$ as a certain model $\mathcal{A}(Q)_{\Sigma_Q,\Sigma_Q(1)} = \mathrm{Bl}^{aff}_{\cup H_j} X_{\Sigma_Q}$ when $Q$ is acyclic.\\

Locally, cluster models are simple: for each cone $\sigma_I\in\Sigma$, writing $J=I\cap A$ and $K=I\setminus A$, a cluster model has a chart of the form $\mathcal A(A_1)^J\times\mathbb C^K\times(\mathbb C^*)^r$ (Proposition 2.28), where a smooth toric variety would have $\mathbb C^{J\sqcup K}\times(\mathbb C^*)^r$. In this sense the only local difference between smooth toric geometry and the geometry of cluster models is the appearance of factors of $\mathcal A(A_1)\cong\mathbb C^2\setminus\{xy=1\}$.\\

To mirror this construction on the $A$-side we follow Gammage and Le, for whom affine blowups are mirror to handle attachments. We begin from the toric mirror $X^\vee_\Sigma=(T^*T^n,\partial^\infty(-\mathbb L_\Sigma))$ and, for each ray $\rho_i\in A$, attach a handle $h_i=H_i^\vee\times T^*D^2$  along a ribbon of the piece $F_i$ of the stop corresponding to $\rho_i$, producing a stopped Weinstein manifold $\mathcal X=\mathcal X(Q)_{\Sigma,A}$. Here $H_i^\vee$ is the mirror of the center $H_i$, itself a cluster model by Theorem \ref{natural CM HMS}$(3)$. That one can attach along one ribbon without destroying the others is not clear \textit{a priori}, and this is the role of the formalism of $\Delta$-collared sectors in Section 3.\\

\begin{figure}[ht]
    \centering
\includegraphics[width=0.5\linewidth]{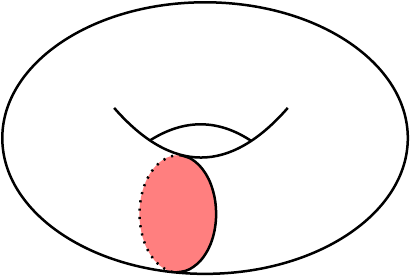}
    \caption{A depiction of the skeleton of the mirror to the cluster variety $\mathcal{A}(A_1)$: a $2$-torus with a disk capping a meridian. The disk piece does not come from any conic Lagrangian in the cotangent bundle $T^* T^2$.}
    \label{A1 skeleton}
\end{figure}

This mirror construction is local in a way akin to cluster models themselves. That is, the mirror $\mathcal X(A_1)$ of $\mathcal{A}(A_1)$ has skeleton a $2$-torus with a disk capping a meridian (Figure \ref{A1 skeleton}), and the charts above have mirrors $\mathcal X(A_1)^J\times(T^*S^1,\infty)^K\times T^*T^r$ (Corollary \ref{cor:charts}). Globally, the skeleton $\mathcal S$ of $\mathcal X$ is obtained from the FLTZ skeleton $-\mathbb L_\Sigma$ by capping the cone over each ray in $A$ with a disk.\\

Our main result globalizes these local pictures into a natural form of HMS, agreeing with Theorem \ref{Toric HMS} where the two overlap.

\begin{thm}\label{natural CM HMS}
    Let $\A := \mathcal{A}(Q)_{\Sigma,A}$ be a cluster model with corresponding stopped Weinstein manifold $\X$. Then $\X$ is stably polarized, with relative skeleton $\mathcal{S} := \mathcal{S}(Q)_{\Sigma,A}$, and there is an equivalence
$$\Phi: \mu sh_{-\mathcal{S}}(-\mathcal{S})^c \simeq \Fuk(\X) \simeq \Coh(\A).$$
\end{thm}

Moreover, $\Phi$ is natural in the following senses.

\begin{enumerate}
    \item Let $S \subset Q^{mut}$ be a set of mutable vertices, and write $Q^{\backslash S}$ for the freezing of $Q$ at $S$. Then there is a stop/handle removal functor $\mathfrak{sr}_S: \Fuk(\X) \to \Fuk(\mathcal{X}(Q^{\setminus S})_{\Sigma_{\setminus S},A_{\setminus S}})$ which is mirror to restriction along a cluster localization $\mathcal{A} \setminus V(\prod_{s\in S} x_s) \cong \mathcal{A}(Q^{\setminus S})_{\Sigma_{\setminus S},A_{\setminus S}}$.
    \item Let $K \subset Q^{mut}$ be an independent set of mutable vertices such that $K \cap A = \emptyset$, and write $Q^{-K}$ for the deletion of these vertices from $Q$. Then the set of such $K$ parametrizes both the boundary strata $D_K$ of $\mathcal A$ and the ends $N_K$ of $\mathcal X$, which have cores $\mathcal S_K\cong\mathcal S(Q^{-K})_{\Sigma_{-K},A_{-K}}\times\sigma_K$, and there is a microlocalization functor $\mu sh_{-\mathcal{S}}(-\mathcal{S}) \rightarrow \mu sh_{-\mathcal{S}_{K}}(-\mathcal{S}_{K})$ which is mirror to restriction along a closed stratum inclusion $\mathcal{A}(Q^{-K})_{\Sigma_{-K},A_{-K}} \hookrightarrow \A$.
    \item The centers $H_i$ of the blowups defining $\A$ are also cluster models: $H_i \cong \mathcal{A}(Q^{/i})_{\Sigma_{/i},A_{/i}}$, where $Q^{/i}$ denotes the contraction of $Q$ at a vertex $i$. Then, passing to the base $H_i^\vee$ of the handle $h_i$ attached to mirror blowing up along $H_i$ is, itself, mirror to restriction from the divisor $D_i$ containing $H_i$ to $H_i$.
    \item $\mathcal{X}(Q)_{\Sigma,A} = \mathcal{X}(Q)_{\Sigma,A\setminus\{i\}} \cup h_i$, and under $\Phi$ the gluing pushout on Fukaya categories induced by this handle attachment is identified with the pushout on categories of coherent sheaves (defined by \cite{GL}) induced by the affine blowup $\mathcal{A}(Q)_{\Sigma,A} = \mathrm{Bl}^{aff}_{H_i}\mathcal{A}(Q)_{\Sigma,A\setminus\{i\}}$.
\end{enumerate}

$\Phi$ is also compatible with products of cluster data, and all of the above is natural in the relevant posets (of freezings of $Q$, of deletions from $Q$, etc.).\\

We note immediately that while $(1)$ and $(2)$ of the above have clear toric analogues (namely, naturality under passing to subfans and stars, respectively), $(3)$ and $(4)$ are new features of cluster HMS. After all, the $H_i$ are pointedly non-toric centers, and handle attachments bring one out of the cotangent bundle -- ergo, out of toric HMS. From this perspective, a major point of the present work is to expand the story of toric HMS beyond the cotangent bundle $T^* T^n$ of the torus while still retaining essentially all of its naturalities \textit{and even introducing new ones}. This is the systematic theory called for at the end of Section 1.1. The happy payoff is that in doing so we can also prove mirror symmetry for RFR acyclic cluster varieties, as we explain thoroughly in Section 1.3.\\

Two further results make the comparison with toric HMS precise. First, $\Fuk(\mathcal{X})$ is the limit, along stop removals, of the Fukaya categories of the mirrors of the toric charts, mirroring Zariski descent (Corollary \ref{descent}). Second, $\partial^\infty\mathcal{S}$ is covered by copies of the skeleta of the mirrors of the boundary strata (Corollary \ref{stop cover}). Using the latter, we then show the following

\begin{thm}
    Assume given a cluster model $\A$ with log CY boundary $\partial \A$, as well as its mirror $\X$ with stop at infinity $\partial^\infty \mathcal{S}$. Then there exists an equivalence $$\mu sh_{-\partial^\infty\mathcal{S}}(\mathcal{-\partial^\infty S})^c \simeq \Fuk(F_D) \simeq \Coh(\partial\mathcal{A}),$$ \noindent where $F_D$ is the (stably polarizable) ribbon of $\partial^\infty \mathcal{S}$ in $\X$.
\end{thm}

Under these equivalences, the cup functor $\Fuk(F_D) \to \Fuk(\X)$ is mirror to pushforward along the boundary inclusion, and stop removal of all of $\partial^\infty \mathcal{S}$ is mirror to restriction to $\A \setminus \partial \A$.\\

A cluster model of the form $\mathcal{A}(Q)_{\Sigma, \emptyset}$ is the same thing as a toric variety $X_\Sigma$, and in this case $\partial\mathcal{A}$ is its toric boundary. This recovers the equivalence $\mu sh \simeq \Coh$ of Theorem 7.4.1 of \cite{gammage2022mirror}. We therefore view $F_D$ as a cluster analogue of a fiber of a mirror potential, and expect that it literally is with respect to some linear combination of theta functions.\\

Moreover, under $\Phi$ a generic cotangent fiber of $\mathcal X$ corresponds to $\mathcal O_{\mathcal A}$. Since cluster models are quasi-affine, this implies that $\Fuk(\mathcal{X})$ is generated by a single cotangent fiber with wrapped endomorphisms $H^*(\mathcal{A},\mathcal{O}_{\mathcal{A}})$ (Corollary \ref{single object gen}). In particular, if $\mathcal A$ is a cluster variety, then the wrapped Floer cohomology of a cotangent fiber is the cluster algebra $\Gamma(\mathcal A)$ itself.\\

Nothing in this argument is special to quivers: we expect the same methods to establish HMS for log Calabi–Yau varieties obtained by sufficiently nice blowups of (not necessarily quasi-affine!) toric varieties along non-toric codimension 2 centers, in the spirit of Gross–Hacking–Keel toric models and Hacking–Keating’s work in dimension 2 \cite{HK}. Cluster models are the case in which the relevant combinatorics is governed by a quiver, which is sufficient to recover the interesting case of acyclic cluster varieties.

\subsection{Applications to Cluster and Positroid Varieties}

Write $\Sigma_Q$ for the fan whose cones $\sigma_I$ are indexed by independent sets $I$ of mutable vertices of $Q$. Then, setting $\Sigma=\Sigma_Q$ and $A=\Sigma_Q(1)$ in Theorem \ref{natural CM HMS} produces a cluster model $\mathcal A(Q)_{\Sigma_Q,\Sigma_Q(1)}$ which is isomorphic to $\mathcal{A}(Q)$ whenever $Q$ is acyclic. Further, in this case the mirror $\mathcal X(Q):=\mathcal X(Q)_{\Sigma_Q,\Sigma_Q(1)}$ carries no stop. This yields the following, which simultaneously generalizes the work of \cite{GL} and \cite{SSWZ}.

\begin{thm}\label{Acyclic HMS}
    Let $Q$ be a really full rank acyclic quiver. Then $\mathcal X(Q)$ is a Weinstein manifold with skeleton combinatorially modeled after $\mathcal A(Q)$, and $$\Fuk(\X(Q)) \simeq \Coh(\mathcal{A}(Q)).$$
\end{thm}

Contemporaneous work of Hu and Shende [21] announces that $\mathcal X(Q)$ can be taken to be $\mathcal A(Q)$ itself. We record this as an ansatz, which we assume throughout without further comment.

\begin{ansatz}[*]\label{*}
    For $Q$ really full rank and acyclic, $\mathcal{X}(Q)$ is Weinstein deformation equivalent to $\mathcal{A}(Q)$, equipped with the Weinstein structure induced by any affine embedding.
\end{ansatz}

The ansatz is used only in statements asserting $\Fuk(\mathcal{A}) \simeq \Coh(\mathcal{A})$ for a cluster variety or positroid variety $\mathcal{A}$ itself. All statements about $\Fuk(\mathcal{X}(Q))$ are unconditional. Therefore, for example, we have the following

\begin{cor}\label{RFR Self HMS}
    Assume (*). Then every RFR acyclic cluster variety $\mathcal{A}(Q)$ is self-mirror: $$\Fuk(\mathcal{A}(Q)) \simeq \Coh(\mathcal{A}(Q)).$$
\end{cor}

Therefore RFR acyclic cluster varieties constitute a considerable class of examples in which Conjecture \ref{leaf rec HMS} holds.\\

We note that if $Q$ and $Q'$ are mutation equivalent RFR acyclic quivers then $\mathcal{A}(Q) \cong \mathcal{A}(Q')$, and so $\Fuk(\mathcal{X}(Q)) \simeq \Fuk(\mathcal{X}(Q'))$ (Corollary \ref{mut invar}). This holds even though the two mirrors are built from different fans. Nevertheless, under $(*)$ they are Weinstein deformation equivalent. It would be interesting to see this directly.\\

We now return to positroid varieties, where we end our discussion of results by initiating the passage from the $\Pi_{2,n}$ case of Conjecture \ref{Open MRW} to the conjecture in general. Since Theorem \ref{Acyclic HMS} and (*) prove Conjecture \ref{Open MRW} for acyclic positroid varieties, the question becomes purely combinatorial: when does $\Pi_{u,w}$ admit an acyclic seed? Using well-known techniques for producing seeds for positroid varieties, we answer this in Section 5 in a range including nontrivial cells in every $\mathrm{Gr}(k,n)$.

\begin{cor}\label{Squareless HMS} Let $u \leq w$ in $S_n$ with $w$ $k$-Grassmannian, and assume that no antidiagonal of the Le-diagram of $u \leq w$ contains 3 or more dots. Then $\Pi_{u,w}$ is acyclic, and so $$\Fuk(\Pi_{u,w}) \simeq \Coh(\Pi_{u,w}).$$
\end{cor}

To our knowledge this is the first homological mirror symmetry statement which applies to non-trivial\footnote{I.e., not arising in the positroid stratification of $Gr(k',n')$ for $k' < k$ or $n' < n$.} positroid cells in \textit{every} $Gr(k,n)$. In particular, in Section \ref{apps} we deduce, still assuming (*), that every positroid cell in a finite-type Grassmannian and every skew-shaped cell $\Pi_{\lambda/\mu}$ whose skew shape contains no $3 \times 3$ square is acyclic and so self-mirror. Beyond positroid varieties, the same methods show that double Bott–Samelson varieties $X(\Delta\beta)$ with $\beta$ weak 2-boolean are self-mirror (Corollary 5.13).
\\

From another perspective, since positroid varieties are really full rank and locally acyclic, each is covered by acyclic cluster localizations, and Corollary \ref{RFR Self HMS} shows that each of these is self-mirror: a local version of Conjecture \ref{Open MRW}.\\

The combinatorial input behind these results is a new notion in Coxeter combinatorics. Namely, we call a permutation or positive braid weak (resp. strong, as introduced by Gao and H\"{a}nni) \textit{$k$-boolean} if some (resp. every) reduced word for it uses each simple transposition at most $k$ times (Definition 5.2). When $k = 1$ these notions agree and recover the notion of boolean permutations -- and in fact, a Schubert variety $X_w$ is toric precisely when $w$ is boolean.\\

We expect that the minimal (weak or strong, since the two notions agree for Grassmannian permutations]) $k$-booleanness of the permutation $w$ in a pair $u \leq w$ with $w$ Grassmannian is a rough measure of the ``cluster complexity" of $\Pi_{u,w}$. Indeed, we introduce weak $k$-booleanness here because it bounds how many vertices of each color appear in the quivers $Q_{u,w}$ of Subsection 2.4 that underlie the cluster structures on positroid varieties, and this number in turn directly impacts how many cycles $Q_{u,w}^{mut}$ has. This is what lets us certify acyclicity in Section \ref{apps}: namely, if $w$ is weak $2$-boolean, then $Q_{u,w}$ is acyclic.\\

Gao and H\"{a}nni show that strong $k$-booleanness is not a pattern-avoidance condition. By contrast, we expect that weak $k$-booleanness \textit{is} a pattern-avoidance condition, albeit one determined by infinitely many subwords. We prove this for $k \le 2$, and give evidence for this claim in general, in Lemma 5.4. A more concerted combinatorial study of weak booleanness is beyond the scope of the present work.

\subsection{Outline and Related Work}

In Section 2, after recalling background on cluster theory and the theory of positroid varieties, we introduce cluster models and establish their basic geometry: charts, strata, centers, freezing, and log Calabi–Yau structure. In Section 3, after recalling symplectic background and discussing $\Delta$-collaredness, we construct cluster model mirrors and their skeleta and discuss concomitant properties. In Section 4 we prove Theorem \ref{natural CM HMS} and its naturalities, deduce the boundary equivalence, and record the consequences above. Finally, in Section 5 we study acyclicity of positroid and double Bott–Samelson varieties, introduce weak $k$-booleanness, and deduce the applications above.\\

Our approach to cluster mirror symmetry through toric models is directly inspired by work of Gammage-Le and, especially, by work of [Peng] (which essentially performs our mirror construction in the case of the $\Pi_{2,n}$). It is also related to Hacking and Keating’s formalism of homological mirror symmetry for log Calabi–Yau surfaces via toric models \cite{HK}: indeed, cluster models were defined to imitate a higher-dimensional analogue of Gross–Hacking–Keel toric models. Again, our methods should apply to establish HMS for all log CY varieties arising from nice enough blowups of (not necessarily quasi-affine!) toric varieties along non-toric codimension $2$ centers -- cluster models are just the case in which the relevant combinatorics is governed by a quiver. Finally, our treatment of gluing along ribbons relies on the work of Ganatra, Pardon, and Shende \cite{GPSI} and of Asplund \cite{Asplund}, to which the $\Delta$-collared formalism of Section 3 reduces.

\subsection{Statement on Contemporaneous Work}
The reader will find multiple publications (\cite{SSWZ}, \cite{HS} and the present work), all posted to the arXiv on the same day, which prove (aspects of) homological mirror symmetry for acyclic cluster varieties and their generalizations. Here we provide some context for this.\\

The author began working on this project in the summer of 2025, and by early 2026 she had proven the theorem on HMS for RFR acyclic cluster varieties. In spring of 2026, the author spoke with V. Shende about her work on homological mirror symmetry for acyclic cluster varieties and communicated the fact that she could prove this theorem. She recalls also sharing her conjectures about the self-mirrorness of positroid cells during said conversation. On June 23rd of the same year, the author gave a recorded lecture on a portion of the present paper's content at the Fields Institute, see https://www.youtube.com/watch?v=KGHenftzyv4.\\

On October 19th, the author received an email from Shende explaining that he and a collaborator had, over the course of the preceding weekend, prompted ChatGPT Astra to solve the same problem ``basically by itself". Shende expressed a plan to publish a human readable account in less than $2$ weeks, offering to involve the author if she had made progress on the problem in a sense he specified. Attached was a writeup of the result by ChatGPT Astra, which cites, among other sources, the recording of the author’s Fields Institute talk and public material related to \cite{SSWZ}.\\

At this time the author had plans to complete the present project by the end of the year (by the end of summer 2026 every major result in the present work was known to her). She requested that she be given the chance to complete the work without the involvement of another party. This request was denied, and the parties instead agreed to post simultaneously.\\

Since the time of that exchange the other party has announced, also using AI tools, homological mirror symmetry for all sink-recurrent cluster varieties, and in particular for all positroid cells (see \cite{HS}). These statements were conjectured by the author in late 2025 and shared with her advisor at that time; a suggestion by her advisor that the author be invited to take part in this work was declined. Not knowing the form of their results, therefore, and in order to record their provenance accurately, we leave our discussion of such statements as conjectures.\\

For the reader interested in a better understanding of the interactions between the author and others participating in the same arXiv drop, the relevant correspondence, with personal content removed, is available at https://sites.google.com/view/josiehlavinka/statement-on-functorial-hms-for-acyclic-cluster-varieties.

\subsection{Acknowledgements} We are grateful to Gabriel Beiner, Sonia Choy, William Fisher, Christian Gaetz, Benjamin Gammage, Yibo Gao, Pavel Galashin, Eugene Gorsky, Soyeon Kim, Melissa Sherman-Bennett, Jos\'{e} Simental, and Brian Sun for helpful discussions. We are also grateful to Beiner, Fisher, Zakkai Goriely, Charlotte Keisling, and Sophie McCormick for their close reading of an earlier draft of this work. Special thanks are due to David Nadler, without whose support this paper could not have been completed, and particularly to Peng Zhou. At an earlier stage of this work the author felt that Zhou's contributions were such that he should be a coauthor on this project, an offer which Zhou politely declined. Since then the project has taken on new aspects, but we would be remiss not to note that most of the results in Subsections 3.2 and 3.3 (and assuredly much more) were learned or initially developed in conversations with him. We are very grateful. The author's research is partially supported by the National Science Foundation Graduate Research Fellowship Program under Grant No. DGE-2146752 and by a Simons Dissertation Fellowship in Mathematics.

\subsection{AI Declaration}
Through October 19, 2026, the author used an LLM on this project only for literature search. In the final weeks of writing this use expanded appreciably, and such tools have since been used as a text editor, proof checker, and prover of technical lemmata. The most consequential result of this is the formalism of $\Delta$-collaredness, which greatly simplified the presentation of Section 3, and we note in particular that Definition 3.3 through Theorem 3.7 were developed in conversations with AI. AI also suggested general organizational streamlinings for Subsection 3.2. The main results, constructions, and conjectures of this paper are the author’s.

\section{Generalized Gross-Hacking-Keel and B-Side Preliminaries}

\subsection{Background on Cluster Varieties}

We work throughout with cluster varieties of skew-symmetric seeds. Let $Q = (V, E, F \subset V)$ denote a quiver with frozen vertices $F$ and skew-symmetric (directed) adjacency matrix $E$, with $|F| = m$ and $|V| = m + n$. Our quivers do not contain $1$- or $2$-cycles. A non-frozen vertex in $Q$ is called \textit{mutable}, and we do not allow edges between frozen vertices. All frozen variables will be inverted.\\

Throughout we (often tacitly) fix a vertex labeling $x$ on $V$ with variables $x_1, ..., x_{m + n}$ generating the field $\mathbb C(x_1,\dots,x_{m+n})$, and we assume that $x(F) = \{x_{n + 1}, ..., x_{m + n}\}$. A pair $(Q, x)$ of a quiver and such a vertex labeling is called a \textit{seed}. Corresponding to such a seed we fix a lattice $M \cong \mathbb{Z}^{m+n}$ with basis $\theta_1, ..., \theta_{m+n}$ and its dual lattice $N \cong M^\vee$ with basis $\theta_1^\vee, ..., \theta_{m+n}^\vee$. We think of $\theta_i$ and $\theta_i^\vee$ as corresponding to the vertex $v_i$ in $Q$ labeled by $x_i$, and given $I \subset V$ we write $M_I :=\bigoplus_{j \notin I} \mathbb{Z} \theta_j$ and $N_I:=\bigoplus_{j\notin I} \mathbb{Z} \theta^\vee_j$.\\

We will write $Q^{mut}$ to denote the induced subgraph of $Q$ on $V \setminus F$, called the \textit{mutable portion} of $Q$. By abuse of notation we will often say that a seed $(Q, x)$ has a property $\mathcal{P}$ if $Q$ does; but, an important exception to this rule is that a seed is called \textit{acyclic} if its mutable portion $Q^{mut}$ is acyclic.\\

We write $\tilde{E}$ to denote the \textit{extended exchange matrix} of $Q$, the leftmost $(m+n) \times n$ submatrix of $E$. In other words, $\tilde{E}$ is the matrix whose top $n \times n$ submatrix is the (directed) adjacency matrix of $Q^{mut}$ and whose last $m$ rows encode the directed adjacency matrix between vertices in $Q^{mut}$ and those in $F$ in $Q$. We will often write \begin{equation}
    c_i := \sum_j \tilde{E}_{j,i} \theta_j \in M
\end{equation} to denote the lattice vector corresponding to the $i$-th column of $\tilde{E}_{j,i}$. Therefore two mutable vertices $v_i, v_j \in Q^{mut}$ are adjacent if and only if $(c_i, \theta_j^\vee) \neq 0$\\

Set $[a]_+ = \text{max}(a, 0)$ and $[a]_- = \text{min}(a, 0)$. Given $v_k$ mutable, we define the \textit{mutation $(Q', x')_{v_k}$ of $(Q, x)$ at $v_k$} to be the seed $((V, E', F), x')$, where
\begin{equation}
    \tilde{E}_{i,j}' =
    \begin{cases}
      -\tilde{E}_{i,j} & \text{if}\ i = k\:\text{or}\ j = k \\
      \tilde{E}_{i,j} + [\tilde{E}_{i,k}]_+[\tilde{E}_{k,j}]_+ - [\tilde{E}_{i,k}]_-[\tilde{E}_{k,j}]_- & \text{otherwise},
    \end{cases}
  \end{equation}
  \noindent and 
  \begin{equation}
  \tilde{x}(v_i) =
      \begin{cases}
      x_i & \text{if}\ i \neq k \\
      \frac{\prod_i x_i^{[\tilde{E}_{i,k}]_+} +\:\prod_i x_i^{-[\tilde{E}_{i,k}]_-}
      }{x_k}& \text{if}\ i = k.
    \end{cases}
  \end{equation}

We note in particular that $\tilde{x}(F) = x(F)$. One also sees that, given $i$, $[\tilde{E}_{i,k}]_+$ (resp. $[\tilde{E}_{i,k}]_-$) is non-zero precisely if the $i$-th coordinate of $c_k$ is positive (resp. negative).\\

Let $\mathcal{Q}$ denote the set of seeds which arise by a sequence of mutations from $(Q, x)$, and let $M(Q)$ denote the set of all elements of $\mathbb{C}(x_1, ..., x_{m+n})$ which are in the image of $x'$ for some $(Q', x') \in \mathcal{Q}$. An $x'$ arising in this way is called a \textit{cluster}, and we call elements of $M(Q)$ \textit{cluster variables}. Subsequent to a choice of labeling $x$, the corresponding seed $(Q, x)$ is called the \textit{initial seed} of $\mathcal{Q}$.

\begin{defn}
    The \textit{cluster algebra} $\Gamma(Q)$ associated to $(Q, x)$ is the $\mathbb{C}$-subalgebra of $\mathbb{C}(x_1, ..., x_{m+n})$ generated by $\mathcal{M}(Q)$, i.e. by the cluster variables, and by $x_{n+1}^{-1}, ..., x_{m+n}^{-1}$. The \textit{cluster variety} $\mathcal{A}(Q)$ associated to $(Q, x)$ is the scheme Spec$(\Gamma(\A(Q))$.
\end{defn}

The cluster $x$ in the initial seed $(Q, x)$ is suppressed from the notation for $\mathcal{A}(Q)$ and $\mathcal{M}(Q)$ but understood from context. In fact, we will often abuse notation and refer to the cluster variety $\mathcal{A}(Q)$ when we have only specified a quiver $Q$. Note that as an abstract variety, $\mathcal{A}(Q)$ only depends on the mutation equivalence class of $Q$.\\

Before giving examples, we fix the convention (to be standing throughout) that if a vertex $v \in Q$ has no neighbors at all then it is assumed to be frozen.

\begin{exmp}\label{a1 quiver}
    Let $A_0$ denote the graph with a single vertex $v$. If we freeze $v$ then $\mathcal{A}(A_0) \cong \mathbb{C}^*$. If we don't freeze $v$, set the initial seed to be $(A_0, x(v) = x)$. Mutating at $v$ gives the seed $(A_0, x(v) = x^{-1})$. Mutating at $v$ again gives the seed $(A_0, x(v) = x)$. Thus we see in either scenario that $\mathcal{A}(A_0) \cong \mathbb{C}^*$, and so we once and for all fix the notation that $A_0$ denotes the quiver containing a single \textit{frozen} vertex. In general, $\mathcal{A}(Q) \cong (\mathbb{C}^*)^{|V(Q)|}$ if and only if $Q$ has no edges.

    \begin{exmp}
        Now let $A_1$ denote the quiver with two vertices, one mutable one $m$ and one frozen one $f$, and an arrow $ar$ from the former to the latter.
    % https://q.uiver.app/#q=WzAsMixbMCwwLCJcXGJ1bGxldCJdLFsxLDAsIlxcYnVsbGV0IixbMTgwLDYwLDYwLDFdXSxbMCwxXV0=
\[\begin{tikzcd}
    A_1 := \underset{m}{\bullet} & \textcolor{rgb,255:red,92;green,214;blue,214}{\underset{f}{\bullet}}
	\arrow[from=1-1, to=1-2, "\textcolor{gray}{ar}"']
\end{tikzcd}\]

Fix an initial seed $(A_1, x(m) = x_1, x(f) = x_2)$. Mutating at $m$ gives the seed $(A_1, x(m) = \frac{1 + x_2}{x_1}, x(f) = x_2)$. Mutating at $m$ again gives the seed $(A_1, x(m) = x_1, x(f) = x_2)$. Setting $x_1' := \frac{1 + x_2}{x_1}$, we therefore see that $\mathcal{A}(A_1) \cong$ Spec $\mathbb{C}[x_1, x_1', x_2^\pm]/(x_1x_1' - 1 = x_2)$, which is isomorphic to Spec $\mathbb{C}[x, y][\frac{1}{xy - 1}]$. This scheme is isomorphic to the complement of a smooth hyperbola in $\mathbb{C}^2$. 
    \end{exmp}

In general, $A_n$ will denote the quiver whose underlying graph consists of a directed path with $n + 1$ vertices, such that the first $n$ vertices are unfrozen and the last is frozen.
\[\begin{tikzcd}
	A_n := \bullet & \bullet & {...} & \bullet & \textcolor{rgb,255:red,92;green,214;blue,214}{\bullet}
	\arrow[from=1-1, to=1-2]
	\arrow[from=1-2, to=1-3]
	\arrow[from=1-3, to=1-4]
	\arrow[from=1-4, to=1-5]
\end{tikzcd}\]

One can show that $\mathcal{A}(A_n) \cong$ Spec $\mathbb{C}[x_1, x_1', \dots, x_n, x_n', q^{\pm}]/(x_ix_i' = x_{i-1} + x_{i+1}\,,\forall i \in [n])$, where we take the convention that $x_0 = 1$ and $x_{n+1} = q$.
\end{exmp}

It is often useful to have a nice store of open subsets of a cluster variety which are themselves cluster varieties. One way to produce them is by freezing extra vertices in a seed for $\Gamma(Q)$. That is, fix any quiver $Q$, and assume that $(Q', x')$ is mutation equivalent to $(Q, x)$. Fix a set $S=\{v_{i_1}, ..., v_{i_k}\}$ of mutable vertices in $Q'$. Then, we write $Q'^{\backslash S}$ to denote the quiver which agrees with $Q'$ save with all the $v_{i_j}$ frozen. Fix an initial seed $(Q'^{\backslash S}, x')$. Then there is a natural containment $$\Gamma(Q'^{\backslash S}) \subset \Gamma(Q)[x_{v_{i_1}}^{-1}, ..., x_{v_{i_k}}^{-1}].$$ 

\begin{defn}
    $\Gamma(Q'^{\backslash S})$ is called a \textit{cluster localization} of $\Gamma(Q)$ if the containment above is an isomorphism. In this case the induced open embedding $\mathcal{A}(Q'^{\backslash S}) \hookrightarrow \mathcal{A}(Q)$ is called a \textit{cluster chart}.
\end{defn}

In general it can be somewhat technical to check that $\Gamma(Q^{\backslash S})$ is a cluster localization of $\Gamma(Q)$. The condition of local acyclicity is useful to guarantee that this is the case.

\begin{defn}
    A cluster variety $\mathcal{A}(Q)$ is called \textit{acyclic} if there exists some acyclic seed $(Q', x') \in \mathcal{Q}$. A cluster variety is called \textit{locally acyclic} if it admits a Zariski cover by a finite number of acyclic cluster charts. 
\end{defn}

The upshot is that if $\mathcal{A}(Q'^{\backslash S})$ is locally acyclic, then $\Gamma(Q'^{\backslash S})$ is a cluster localization of $\Gamma(Q)$ (see \cite{Mul}). It is not always possible to cover a locally acyclic cluster variety $\mathcal{A}(Q)$ with cluster charts simply consisting of cluster varieties for freezings of $Q$, however, i.e. without actually mutating to $Q'$ and then freezing. This is one of the reasons we restrict our attention to acyclic cluster varieties.\\

We will also utilize the following extreme examples of acyclic cluster varieties.

\begin{defn}
    A cluster variety is called \textit{isolated} if some mutation of $Q$ has no edges between mutable vertices (equivalently, if the top $n \times n$ submatrix of $\tilde{E}$ with respect to some mutation of $Q$ is zero).
\end{defn}

For example, $\mathcal{A}(A_1)$ is an isolated cluster variety, but all others $\mathcal{A}(A_n)$'s are only acyclic. Locally acyclic cluster varieties are in fact covered by cluster charts of \textit{isolated} cluster varieties (see \cite{Mul}).\\

Cluster varieties are, indeed, varieties (integral affine schemes). They can fail to be smooth, though. This occurs for $\mathcal{A}(Q_{Mar})$, where $Q_{Mar}$ denotes following quiver:

% https://q.uiver.app/#q=WzAsNCxbMiwwLCJcXGJ1bGxldCJdLFszLDEsIlxcYnVsbGV0Il0sWzEsMSwiXFxidWxsZXQiXSxbMCwxLCJcXG1hdGhjYWx7TX0gOj0iXSxbMCwxXSxbMCwxLCIiLDIseyJvZmZzZXQiOi0zfV0sWzIsMF0sWzIsMCwiIiwyLHsib2Zmc2V0IjotM31dLFsxLDJdLFsxLDIsIiIsMSx7Im9mZnNldCI6LTN9XV0=
\[\begin{tikzcd}
	&& \bullet \\
	{\mathcal{M} :=} & \bullet && \bullet
	\arrow[from=1-3, to=2-4]
	\arrow[shift left=3, from=1-3, to=2-4]
	\arrow[from=2-2, to=1-3]
	\arrow[shift left=3, from=2-2, to=1-3]
	\arrow[from=2-4, to=2-2]
	\arrow[shift left=3, from=2-4, to=2-2]
\end{tikzcd}\]

\noindent which is not even Noetherian. If $\mathcal{A}(Q)$ is locally acyclic, however, then certain conditions on the extended exchange matrix $\tilde{E}$ of a(n initial) seed $(Q, x)$ can ameliorate this situation. %%In the following, we write $\tilde{E}$ to denote the $n \times n$ submatrix of $E$ consisting only of columns which correspond to mutable vertices of $(Q, x)$.

\begin{defn}
    A seed $(Q, x)$ is \textit{full rank} if $\tilde{E}$ has full rank as a matrix over $\mathbb{Q}$. If moreover the rows of $\tilde{E}$ span $\mathbb{Z}^n$ as a $\mathbb{Z}$-module, then $(Q, x)$ is \textit{really full rank}. A cluster variety $\mathcal{A}(Q)$ is said to be (really) full rank (or (R)FR for short) if there exists some $(Q', x') \in \mathcal{Q}$ which is (really) full rank, respectively.
\end{defn}

$Q$ is RFR if and only if its columns $c_i$ are part of a $\mathbb Z$-basis of $M$. Equivalently, and to be used in the sequel, this is the same as saying that one can choose a basis of elements $c_i^\vee\in N$ with $(c_k,c_i^\vee)=\delta_{ik}$. Further, (really) full rank-ness turns out to be a mutation invariant of $\tilde{E}$, and so will hold for $(Q, x)$ itself if it holds for some $(Q', x') \in \mathcal{Q}$. A particular upshot of the condition of full rankness is the following

\begin{thm}[\cite{Mul} Corollary 7.9.] \label{loc. acyclic smoothness}
    A full rank locally acyclic cluster variety is smooth.
\end{thm}

The class of \textit{really} full rank quivers in particular is closed under some further combinatorial operations besides freezing, to which we will soon give geometric interpretations.

\begin{defn}\label{contraction}
    Assume given an RFR acyclic quiver $Q$ and an independent set, i.e. set containing no edges, of mutable vertices $I \subset Q^{mut}$.

    \begin{enumerate}
        \item The \textit{deletion} of $I$ from $Q$, denoted $Q^{-I}$, is the quiver which arises from $Q$ by deleting every vertex in $I$ and freezing all of their neighbors. 
        \item Assume $I = \{s\}$ contains a single element. The \textit{contraction} of $s$ from $Q$, denoted $Q^{/s}$, is $Q^{-s}$ with its frozen lattice quotiented by $c_s$. That is, identifying the data of $Q^{-s}$ with that of $\tilde{E}(Q^{-s})$, $Q^{/s}$ is the quiver determined by the extended exchange matrix given by composing $\tilde{E}(Q^{-s})$ with the quotient by $c_s$.
    \end{enumerate}
\end{defn}

All of $Q^{\backslash S}$, $Q^{-I}$, and $Q^{/s}$ are acyclic and RFR if $Q$ was (see Section 4.2 of \cite{hlavinka2026chow}). If $M$ was the lattice corresponding to $Q$, it is not hard to see that $M_I$ is the one corresponding to $Q^{-I}$ and that $M_s/\mathbb{Z}c_s$ is the one corresponding to $Q^{/s}$. Also, it is not hard to see that $(Q^{/i})^{-K}$ is $(Q^{-K})^{/i}$.\\

We will utilize and generalize a stratification which exists for any cluster variety but which takes an extremely nice form for really full rank acyclic ones. All of the following facts are due to Lam and Speyer, see \cite{LSII}.

\begin{defn}\label{v/d defn}
    Given an initial seed $(Q, x)$ defining a cluster variety $\mathcal{A}(Q)$, the \textit{$V/D$ decomposition} of $\mathcal{A}(Q)$ is the decomposition of $\mathcal{A}(Q)$ into subsets on which a given set of initial cluster variables $x_i$ is identically $0$ and the remaining variables are invertible. That is, the $V/D$ decomposition comes from writing $$\mathcal{A}(Q) = \bigsqcup_{I = \{i_1, ..., i_k\} \subset [n]} V(x_{i_1})\cap ...\cap V(x_{i_k}) \cap D(\prod_{j \notin I} x_j),$$ \noindent with the understanding that many of these pieces are empty. The non-empty pieces of this expression are the \textit{V/D pieces} of $\mathcal{A}(Q)$.
\end{defn}

In general the $V/D$ decomposition can be somewhat of a mess. The first sign of its usefulness for us is that the indexing set of its nonempty strata simplifies when $(Q, x)$ is acyclic. In this case the $V/D$ pieces are parametrized by independent subsets of (the underlying graph of) $Q$. 

\begin{lem}[\cite{LSII} Lemma 3.1.]\label{independent support}
    Let $(Q, x)$ be an acyclic seed. Then for any point $x = (x_1, ..., x_{m+n}) \in \mathcal{A}(Q)$, the set of indices $i$ for which $x_i = 0$ labels an independent set of vertices in $Q^{mut}$.
\end{lem}

Following \cite{LSII}, given an acyclic seed $(Q, x)$ we write $\mathcal{O}_\mathcal{I}$ for the subset of $\mathcal{A}(Q)$ on which the variables $x_i$ with $i \in \mathcal I$ are identically $0$ and the rest are invertible. Then, in the case that $(Q, x)$ is really full rank, Lam and Speyer show that the $\mathcal{O}_I$ themselves are also extremely nice.

\begin{thm}[\cite{LSII} Proposition 3.5.]\label{V/D pieces} Assume that $(Q, x)$ is really full rank and acyclic. Then for any independent set $\mathcal{I} \subset Q$ with $|\mathcal{I}| = k$, $\mathcal{O}_\mathcal{I}$ is isomorphic to $\mathbb{A}^k \times \mathbb{G}_m^{m + n - 2k}$.
\end{thm}

\begin{exmp}
    Thinking of $\mathcal{A}(A_1)$ as $\mathbb{C}^2$ minus the hyperbola $xy = 1$, the $V/D$ pieces are the locus on which $x \neq 0$ and the locus on which $x = 0$, respectively. The former is $\mathcal{A}(A_1)$ minus the $y$-axis, which is isomorphic to a copy of $(\mathbb{C}^*)^2$; the latter is the $y$-axis.
\end{exmp}

Another perk of the acyclic setting is that one can write down the cluster variety corresponding to an acyclic quiver explicitly in coordinates.

\begin{exmp}\label{complete int}
    Let $(Q, x)$ denote an acyclic seed with $|V| = m + n$, and write $\tilde{E}$ for the extended exchange matrix of $Q$. Assume that the non-frozen vertices of $Q$ are labeled by $x_1, ..., x_n$. Then by \cite{LSII} Theorem 2.4., the cluster variety $\mathcal{A}(Q)$ is isomorphic to $$\text{Spec}\:\mathbb{C}[x_1, ..., x_n, x_1', ..., x_n', x_{n+1}^\pm, ..., x_{m+n}^\pm]/(x_j x_j' = \prod_{i=1}^{m+n} x_i^{[\tilde{E}_{i,j}]_+} + \prod_{i=1}^{m+n} x_i^{-[\tilde{E}_{i,j}]_-} , \forall j),
    $$
\end{exmp}

This exhibits any acyclic cluster variety as a complete intersection in $\mathbb{C}^{2n} \times (\mathbb{C}^*)^m$.\\

The defining equations of $\mathcal{A}(Q)$ look as if they are coming from affine blowups along centers of the form $V(x_j) \cap V(\prod_{i=1}^{m+n} x_i^{[\tilde{E}_{i,j}]_+} + \prod_{i=1}^{m+n} x_i^{-[\tilde{E}_{i,j}]_-})$, and we close this subsection by defining them as such. To do so, we fix notation and consider a sequence of closed immersions $H \subset D \subset X$ where $H$ is a smooth divisor contained in a single component of  $D$ and $D$ is a simple normal crossings divisor in $X$.

\begin{defn}
    The \textit{affine blowup} Bl$_H^{\text{aff}}X$ of $X$ along $H \subset D$ is defined to be $\text{Bl}_H X \setminus \tilde{D}$, where $\tilde{D}$ is the proper transform of the component of $D$ containing $H$.
\end{defn}

In all cases we will consider, $D$ will be (the proper transform of) a toric boundary divisor in some toric variety, and $H$ will explicitly be defined as the intersection of such a $D$ with some vanishing locus. As such $D$ is suppressed from the notation but will be clear from context.\\

We will actually have to perform \textit{iterated} affine blowups along (proper transforms of) certain centers to get $\mathcal{A}(Q)$. It will also turn out that in all cases we will consider, the order of these blowups doesn't matter. Therefore, when considering multiple codimension $2$ subsets $H_i \subset D_i$ which we perform iterated blowups along to define a variety $\text{Bl}_{\tilde{H}_{i_{n}}}^{\text{aff}}\dots \text{Bl}_{H_{i_1}}^{\text{aff}} X$, we will set $H = \bigcup_i H_i$ and $D = \bigcup_i D_i$ and instead write $\mathrm{Bl}_H^{\text{aff}} X$ for the same space, i.e. in such situations $$\mathrm{Bl}_H^{\text{aff}} X := \text{Bl}_{\tilde{H}_{i_{n}}}^{\text{aff}}\dots \text{Bl}_{H_{i_1}}^{\text{aff}} X.$$

That $\mathrm{Bl}_H^{\text{aff}} X$ does not actually denote the affine blowup along the snc subvariety $H \subset X$ will not cause confusion if the reader keeps in mind that we \textit{always} blow up along smooth centers in this work. We hope its employment will be clear from context.\\

We may now express any acyclic cluster variety as an iterated affine blowup of a toric variety. Given an acyclic quiver $Q$, define a fan $\Sigma_Q \subset \mathbb{R}^{m+n}$ whose cones are all the $\sigma_I = \mathbb{R}_{\geq 0} \theta_{i_1}^\vee + ... + \mathbb{R}_{\geq 0} \theta_{i_k}^\vee \subset N_\mathbb{R} \cong \mathbb{R}^{n+m}$ such that $I = \{i_1, ..., i_k\}$ labels an independent set of $Q^{mut}$. We then write $X_{\Sigma_Q}$ for the toric variety of $\Sigma_Q$. Equivalently, $X_{\Sigma_Q}$ is the union of $(\mathbb{C}^*)^{n+m}$ with all loci of the form $\bigcap_{i \in I} V(x_i)\:\cap \bigcap_{j \notin I} D(x_j)$ inside of $\mathbb{C}^{n+m}$.\\

Set $D_i = V(x_i)$ and $H_j = V(\prod_{i=1}^{m+n} x_i^{[\tilde{E}_{i,j}]_+} + \prod_{i=1}^{m+n} x_i^{-[\tilde{E}_{i,j}]_-}) \cap D_j$. Then, we have the following

\begin{thm}
\label{blaff acyclic}
    Let $Q$ be an acyclic RFR quiver. Then $$\mathcal{A}(Q) = \text{Bl}_{\tilde{H}_{i_{n}}}^{\text{aff}}\dots \text{Bl}_{H_{i_1}}^{\text{aff}} X_{\Sigma_Q} = Bl_H^{\text{aff}} X_{\Sigma_Q},$$

    where the blowup order may be chosen arbitrarily.
\end{thm}

\begin{proof}
   We show first that the order of blowups doesn't matter. To do so, set $P_j := \prod_{i=1}^{m+n} x_i^{[\tilde{E}_{i,j}]_+} + \prod_{i=1}^{m+n} x_i^{-[\tilde{E}_{i,j}]_-}$. By construction $X_{\Sigma_Q}$ is covered by affine opens of the form $\mathbb{C}^{|I|} \times (\mathbb{C}^*)^{n + m - |I|}$, where $I = \{i_1, ..., i_k\}$ is an independent set of $Q$ and the coordinates on the affine space factor are $x_{i_1}, ..., x_{i_k}$, and so we may check commutativity of the blowups on these charts. Assume $i, j \in I$, so that $H_i \cap \mathbb{C}^{|I|} \times (\mathbb{C}^*)^{n + m - |I|}$ and $H_j \cap \mathbb{C}^{|I|} \times (\mathbb{C}^*)^{n + m - |I|}$ are non-empty. Then by definition, $\mathrm{Bl}_{H_j}^{\text{aff}} (\mathbb{C}^{|I|} \times (\mathbb{C}^*)^{n + m - |I|})$ is isomorphic to $$\text{Spec}\:\mathbb{C}[x_{i_1}, ..., x_{i_k}, x_j', x_{n+1}^\pm, ..., x_{n+m}^\pm][x_i^{-1}, \forall i \notin I]/(x_j x_j' = P_j, \forall j),$$ and the preimages of $D_i$ and $H_i$ are still $V(x_i)$ and $V(P_i) \cap D_i$. Since $P_i$ does not contain $x_j$ by the independence assumption, this means that $\mathrm{Bl}_{\tilde{H}_i}^{\text{aff}}Bl_{H_j}^{\text{aff}} (\mathbb{C}^{|I|} \times (\mathbb{C}^*)^{n + m - |I|})$ is isomorphic to
   $\text{Spec}\:\mathbb{C}[x_{i_1}, ..., x_{i_k}, x_i', x_j', x_1^\pm, ..., x_{n+m}^\pm]/(x_j x_j' = P_j, x_ix_i' = P_i)$.\\

   One sees that the same argument, save now employing the non-presence of $x_i$ in $P_j$, tells us that blowing up in the opposite order gives the same variety. Therefore the claim is true when only performing two blowups, and the argument for the general case proceeds similarly.\\

   For two independent subsets $I$ and $J$ of $Q$, we have that $\mathbb{C}^{|I|} \times (\mathbb{C}^*)^{n + m - |I|}$ and $\mathbb{C}^{|J|} \times (\mathbb{C}^*)^{n + m - |J|}$ intersect precisely in $\mathbb{C}^{|I \cap J|} \times (\mathbb{C}^*)^{n + m - |I \cap J|}$. Therefore the argument above tells us that $\mathrm{Bl}_H^{\text{aff}} X_{\Sigma_Q}$ is covered by affine opens of the form $\mathrm{Bl}_H^{\text{aff}} (\mathbb{C}^{|I|} \times (\mathbb{C}^*)^{n + m - |I|})$ spanning over all independent sets $I \subset V$, each of which is isomorphic to $A_I:=\text{Spec}\:\mathbb{C}[x_{i_1}, ..., x_{i_k}, x_{i_1}', ..., x_{i_k}', x_{1}^\pm, ..., x_{m+n}^\pm]/(x_j x_j' = P_j, \forall j \in I)$, and that the overlap of two such affine opens $A_I$ and $A_J$ in $\mathrm{Bl}_H^{aff} X_{\Sigma_Q}$ satisfies $A_I \cap A_J = A_{I \cap J}$.\\

   By Lemma \ref{independent support}, we know that $\mathcal{A}(Q)$ is covered by opens of the form $$\mathcal{A}_I :=\text{Spec}\: \mathbb{C}[x_1, ..., x_n, x_1', ..., x_n', x_{n+1}^\pm, ..., x_{m+n}^\pm][x_{i}^{-1},\forall i \notin I]/(x_j x_j' = P_j, \forall j).$$ This implies that $\mathcal{A}_I \cap \mathcal{A}_J = \mathcal{A}_{I \cap J}$. But since $x_i$ is invertible for any $i \notin I$, one sees that on the locus $\mathcal{A}_I$ we have $x_i' = \frac{P_i}{x_i}$; i.e., that the variables $x_i'$ for $i \notin I$ are already in the corresponding ring. Therefore $\mathcal{A}_I \cong A_I$ compatibly with overlaps, and this completes the proof.
\end{proof}

\begin{rem}
    It is not hard to see that
    the $V/D$ decomposition of $\mathcal{A}(Q)$ is the preimage of the stratification of $X_{\Sigma_Q}$ into torus orbits under the iterated blowup maps.
\end{rem}

\begin{exmp}\label{isolated b-side local model}
    We previously described $\mathcal{A}(A_1)$ as $\mathbb{C}^2\setminus\{xy = 1\}$. Theorem \ref{blaff acyclic} implies that also $$\mathcal{A}(A_1) \cong Bl^{\text{aff}}_{(0, -1)} (\mathbb{C} \times \mathbb{C}^*).$$ 
    
    We leave it to the reader to check that these two spaces are the same.
\end{exmp}

\subsection{Cluster Models}

Given any quiver $Q$, Gross, Hacking, and Keel gave a construction of a \textit{truncated cluster variety} $\overline{\mathcal{A}(Q)}$ associated to $Q$ by performing an affine blowup of a toric variety.

\begin{defn}[GHK's Construction of Truncated  Varieties, \cite{GHK}]
    Let $(Q, x)$ be an RFR seed with mutable vertices labeled $\{x_j\}$, and write $\overline{\Sigma_Q} \subset N_\mathbb{R} \cong \mathbb{R}^{m+n}$ for the fan consisting only of cones $\rho_j = \mathbb{R}_{\geq 0} \theta_j^\vee$. We write $V(x_j)$ to denote the torus orbit closure corresponding to $\rho_j$, and we consider the \textit{exchange polynomial} $P_j := \prod_{i=1}^{m+n} x_i^{[\tilde{E}_{i,j}]_+} + \prod_{i=1}^{m+n} x_i^{-[\tilde{E}_{i,j}]_-}$%%; note that the vanishing locus of $P_j$ is the closure of the vanishing locus of the character $\chi_j = \frac{\prod_{i=1}^{m+n} x_i^{[\tilde{E}_{i,j}]_+}}{ \prod_{i=1}^{m+n} x_i^{-[\tilde{E}_{i,j}]_-}}$ on the torus $(\mathbb{C}^*)^{m + n}$
    . Write $H_j$ to denote the locus $V(P_j) \cap V(x_j) \subset X_{\overline{\Sigma_Q}}$. Then, the \textit{truncated cluster variety} $\overline{\mathcal{A}(Q)}$ associated to $Q$ is defined as $$\overline{\mathcal{A}(Q)} := Bl_{\bigcup_j H_j}^{aff} X_{\overline{\Sigma_Q}},$$

    \noindent where we implicitly use that the order of the appropriate blowups doesn't matter, as above. It is not hard to check that $H_j$ is a smooth hypersurface in $V(x_j)$. Further, because $Q$ is RFR the column $c_j$ corresponding to $P_j$ is primitive, i.e. $P_j$ is irreducible, and so $H_j$ is connected.
\end{defn}

Gross, Hacking, and Keel showed that $\overline{\mathcal{A}(Q)}$ is an open subset of the cluster variety $\mathcal{A}(Q)$ whose complement has codimension at least $2$, and that the affinization of $\overline{\mathcal{A}(Q)}$ is $\mathcal{A}(Q)$ when $Q$ is locally acyclic. Unfortunately, $\overline{\mathcal{A}(Q)}$ and $\mathcal{A}(Q)$ differ in almost all examples. They are different varieties even for $Q = A_1 \sqcup A_1$: $\mathcal{A}(Q)$ is $\mathcal{A}(A_1) \times \mathcal{A}(A_1)$, but $\overline{\mathcal{A}(Q)} = \mathcal{A}(Q) \setminus (V(x_1) \cap V(x_2))$ (per \cite{GL}).\\

Our observation is that Theorem \ref{blaff acyclic} implies that if one considers fans with higher dimensional cones, a version of GHK's construction suffices to construct acyclic cluster varieties \textit{verbatim}. To interpret this we fix the following notation.

\begin{defn}\label{sigmaq}
    Given a seed $(Q, x)$ with extended exchange matrix $\tilde{E}$, we construct a fan $\Sigma_Q$ as follows. As in GHK's construction, the elements of $\Sigma_Q(1)$ are rays $\rho_j = \mathbb{R}_{\geq 0} \theta_j^\vee$ for every $x_j$ a label of a mutable vertex of $Q$.\\
    
    The departure from before is that we also include higher dimensional cones in $\Sigma_Q$, which we define to be the fan consisting of all cones $\sigma_I = \mathbb{R}_{\geq 0} \theta_{i_1}^\vee + ... + \mathbb{R}_{\geq 0} \theta_{i_k}^\vee$ such that $I = \{i_1, ..., i_k\}$ labels an independent set of $Q^{mut}$, i.e. such that $I$ contains no two adjacent vertices.
\end{defn}

Note in particular that $\overline{\Sigma_Q} = \Sigma_Q(1) \subset \Sigma_Q \subset N_\mathbb{R}$. In the sequel, we will think of the fan corresponding to a quiver $Q$ as living in the vector space $N_\mathbb{R}$. Now we generalize GHK's construction. From now on, unless otherwise stated, all quivers are RFR and acyclic.

\begin{con}[Generalized GHK's Construction/Cluster Models]\label{GHK+}
    Fix a quiver $Q$, a subfan $\Sigma \subset \Sigma_Q$, and a subset $A \subset \Sigma(1)$. We will call this last piece of data a set of \textit{active rays}, and the triple $(Q, \Sigma, A)$ a \textit{cluster datum}.\\
    
    The \textit{cluster model} $\mathcal{A}(Q)_{\Sigma, A}$ of $Q$ with respect to $\Sigma$ and $A$ is defined to be $$\mathcal{A}(Q)_{\Sigma, A} :=\text{Bl}_{\bigcup_{\rho_j \in A} H_j}^{\text{aff}} X_\Sigma.$$
\end{con}

Again the order of the blowups does not matter, essentially by the proof of Theorem \ref{blaff acyclic}. The motivating examples are the following.

\begin{exmp}
    When $A = \emptyset$, Construction \ref{GHK+} produces the toric variety $X_\Sigma$: $\mathcal{A}(Q)_{\Sigma, \emptyset} \cong X_{\Sigma}$.
\end{exmp}

\begin{exmp}
    When $\Sigma = \Sigma_Q$, $A = \Sigma_Q(1)$ and $Q$ is acyclic, Theorem \ref{blaff acyclic} tells us that Construction \ref{GHK+} produces $\mathcal{A}(Q)_{\Sigma_Q,A} \cong \mathcal{A}(Q)$.
\end{exmp}

\begin{exmp}
    When $\Sigma = \overline{\Sigma_Q}$ and $A = \Sigma_Q(1)$, the main theorem of \cite{GHK} tells us that for any quiver $Q$ the output is $\mathcal{A}(Q)_{\overline{\Sigma_Q}, A} \cong \overline{\mathcal{A}(Q)}$. That is, $\mathcal{A}(Q)_{\overline{\Sigma_Q}, \overline{\Sigma_Q}(1)} \cong \overline{\mathcal{A}(Q)}$.
\end{exmp}

Given a property $\mathcal{P}$ previously defined for quivers or cluster varieties (acyclicity, isolatedness, RFR-ness, ...), a cluster model is said to have property $\mathcal{P}$ if $Q$ or $\mathcal{A}(Q)$ does. Also, in the sequel we will identify rays $\rho_i$ of $\Sigma$ with the corresponding vertices $i \in Q$ with abandon, and so will write things like $I \cap A$ (where $I \subset Q^{mut}$ and $A$ is an active set of rays) to denote $\{\rho_i\:|\: i \in I\} \cap A$.\\

To begin understanding the geometry of cluster models in general we record a proposition. In it and in the sequel we will find it convenient to think of the loci $V(P_i)$, when restricted to certain open charts in $\mathcal{A}(Q)_{\Sigma, A}$, instead as the locus on which a monomial is set to be $-1$. The monomial in question will be $y_i = x^{c_i}$, the product of the coordinate functions with exponents coming from the column $c_i$ of $\tilde{E}$. This is analogous to the definition of contracted quivers in Definition \ref{contraction}: it is in some sense the geometric result of trivializing the frozen character corresponding to $c_i$.

\begin{prop}\label{chart lemma}
    Assume given a cluster datum $(Q, \Sigma, A)$ and $\sigma_I \in \Sigma$. Set $J := I \cap A$ and $K := I \setminus A$, and write $\mathcal{A}(Q)_{\Sigma, A}|_{U_{\sigma_I}}$ for the preimage of the open toric subvariety $U_{\sigma_I} \subset X_\Sigma$ under the blowup map. For $i \in I$ write $y_i = x^{c_i}$. Then there is an isomorphism $$\mathcal \A(Q)_{\Sigma,A}\big|_{U_{\sigma_I}} \cong \prod_{j \in J} \mathcal A(A_1) \times \prod_{k \in K} (\mathbb C_{x_k} \times \mathbb{C}_{y_k}^*)\times (\mathbb C^*)^{m+n-2|I|}$$ \noindent where the $j$-th copy of $\mathcal{A}(A_1)$ is the affine blowup of $\mathbb{C}_{x_j} \times \mathbb{C}^*_{y_j}$ at the point $(0, -1)$, and where $x_k$ and $y_i$ are coordinates on the indicated factors.
\end{prop}

\begin{proof}
For $j \in J$, write $P_j = u_j(1+y_j)$ with $u_j:=x^{[-c_j]_+}$. Since the factor $u_j$ \textit{a priori} involves only frozen variables and variables
$x_k$ with $v_k$ adjacent to $v_j$, the independence of $I$ tells us that $u_j$ is a unit on
$U_{\sigma_I}$. Therefore $H_j\cap U_{\sigma_I}=\{x_j=0,\ y_j=-1\}$. Further, $x_a$ is a unit on
$U_{\sigma_I}$ for $a \in A\setminus I$, and so $H_a \cap U_{\sigma_I} = \emptyset$. We see also that $(c_i,\theta_k^\vee) =\tilde{E}_{k,i}=0$ for all $i,k\in I$ since $I$ is an independent set, and
so $c_i\in \bigoplus_{v_j \notin I} \mathbb {Z} \theta_j = M_I$ for all $i \in I$.\\

Since $Q$ is RFR the $c_i$, $i\in I$, are part of a $\Z$-basis of $M$, and so of the direct summand
$M_I$. Choose a complement $B$ for the span of the $c_i$ inside of $M_I$, i.e. with $M_I=\bigoplus_{i\in I}\Z c_i \oplus B$, so that the rank of $B$  is $m+n-2|I|$.
Then $U_{\sigma_I}\cong\prod_{i\in I}(\C_{x_i}\times\C^*_{y_i})\times (\mathbb{C}^*)^{m + n - 2|I|}$, and each
$H_j$ is $\{(0,-1)\}$ in the $j$-th factor times everything else. The result then follows from the fact that affine blowups commute with products, and from Example \ref{isolated b-side local model}.
\end{proof}

By the preceding proposition, if $i \notin A$ then the preimage of $V(x_i)$ under the blowup map, i.e. the vanishing locus of $x_i$ in $\mathcal{A}(Q)_{\Sigma, A}$, is precisely its strict transform, which we denote by $D_i$. We call these the $D_i$ the \textit{inactive divisors} of $\mathcal{A}(Q)_{\Sigma, A}$. On the other hand, if $i \in A$, then the preimage of $V(x_i)$ is the exceptional divisor $E_i \cong H_i \times \mathbb{A}^1$ of the affine blowup at $H_i$. Whenever $j \notin A$ we will write $H_j$ for the preimage of the center $H_j \subset X_\Sigma$ under the blowup map $\mathcal{A}(Q)_{\Sigma, A} \rightarrow X_\Sigma$, and continue to refer to it as a \textit{center} in $\mathcal{A}(Q)_{\Sigma, A}$. This is also the strict transform of $H_j$, and one sees that $\mathcal \A(Q)_{\Sigma,A\cup\{j\}}=\mathrm{Bl}^{\mathrm{aff}}_{H_j}\mathcal \A(Q)_{\Sigma,A}$ by definition.\\

From the above we learn the following 

\begin{cor}\label{transversality}
    Assume given a cluster model $\mathcal{A}(Q)_{\Sigma, A}$ with inactive divisors $D_i$. Then, every intersection $D_I = \bigcap_{i \in I \subset \Sigma(1) \setminus A} D_i$ is transverse and each such $D_I$ is smooth. Further, given $\rho_i \notin A$, we have $H_j \subset D_j$, while $H_j$ meets $D_I$ $(j \notin I)$ if and only if $\sigma_{j \cup I}\in\Sigma$ (in which case it also does so transversely).
\end{cor}

Further, since the proposition implies that any cluster model $\mathcal{A}(Q)_{\Sigma, A}$ has an affine cover by distinguished opens on which certain global functions (namely, the $D(\prod_{\rho_k\in\Sigma(1),\ k\notin I}x_k$) are invertible, we learn that any cluster model is quasi-affine.\\ 

To parlay this chartwise understanding of cluster models into a global combinatorial and geometric one we need one final auxiliary lemma, whose proof is left to the reader.

%%The last proposition begins to suggest a loose slogan: that cluster models are to (smooth, quasi-affine) toric varieties as (RFR acyclic) cluster varieties are to tori. Therefore we expect that cluster models admit interesting recursive structures with respect to combinatorial operations on $Q$. To see that this is the case we first record a lemma, whose proof we leave to the reader.

\begin{lem}\label{fan operations}
    Assume given an acyclic RFR cluster quiver $Q$, $S \subset Q^{mut}$ any set of vertices, and $I \subset Q^{mut}$ an independent set of vertices. Then $$\Sigma_{Q^{\backslash S}} = \Sigma_Q \setminus \{\sigma_J\:|\: J \cap S \neq \emptyset\} \subset N_\mathbb{R}\:\:\text{and}\:\: \Sigma_{Q^{-I}} =\{\sigma_J\:|\: J \cap I = \emptyset\:\text{and}\:\sigma_{I \cup J} \in \Sigma_Q\} \subset (N_I)_\mathbb{R}.$$
\end{lem}

In other words, $\Sigma_{Q^{-I}}$ is the star of $\sigma_I$ in $\Sigma_Q$ modulo the span of $\sigma_I = \theta_{i_1}^\vee, ..., \theta_{i_{|I|}}^\vee$; this is the sense in which it naturally lives in $(N_I)_\mathbb{R}$. Further, one sees that $X_{\Sigma_{Q^{-I}}}$ can be identified with the closed toric subvariety $\bigcap_{i \in I} V(x_i)$ in $X_{\Sigma_Q}$.\\

In general, given $\Sigma \subset \Sigma_Q$ we write $\Sigma_{\backslash S} \subset N_\mathbb{R}$ and $\Sigma_{-I} \subset (N_I)_\mathbb{R}$ to denote fans defined similar to those in the first and second equations of Lemma 2.30 (resp.), save with every instance of $\Sigma_Q$ replaced by $\Sigma$ itself. Given $A \subset \Sigma(1)$, we write $A_{\backslash S} := A \cap \Sigma_{\backslash S}(1)$ and $A_{-I} := A \cap \Sigma_{-I}(1)$.\\

In fact, we note that every ray of $\Sigma_{Q^{-i}}$ pairs to zero with $c_i$ -- this is true because the mutable vertices of $Q^{-i}$ are exactly those of $Q$ which aren't adjacent to $i$. Therefore $\Sigma_{Q^{-i}}$ actually lies in the sublattice $c_i^\perp=(M_i/\mathbb Zc_i)^\vee\subset N_i$. One can check that, when considered as a fan in this \textit{sublattice}, we have $\Sigma_{Q^{-i}} \cong \Sigma_{Q^{/i}}$. We will always think of this fan as living in the larger lattice when deleting and in the smaller one when contracting. In any event, this means that when $I = \{i\}$ we will also use the notation $\Sigma_{/i} := \Sigma_{-i}$ and $A_{/i} := A_{-i}$.\\

We define an isomorphism of cluster data $(Q,\Sigma,A)\to(Q',\Sigma',A')$ to be a lattice isomorphism $\varphi\colon M\to M'$, together with a bijection $\sigma$ of mutable vertices, such that $\varphi^\vee(\theta'^\vee_{\sigma(i)})=\theta_i^\vee$ and $\varphi(c_i)=c'_{\sigma(i)}$ for each mutable $v_i$, and such that $\varphi^\vee$ carries $\Sigma'$ and $A'$ onto $\Sigma$ and $A$. Such a $\varphi$ induces isomorphisms of toric varieties preserving the centers $H_i$, and hence of cluster models. We are therefore allowed to specify cluster data, up to isomorphism, by instead specifying a lattice with marked columns, a fan and an active set. A particular example of this that we will use in the sequel is the cluster datum $(\mathbb Zc,\{0\},\emptyset)$, the datum of the quiver consisting of a single frozen vertex frozen variable $x^c$.\\

Given two cluster data $(Q, \Sigma, A)$ and $(Q', \Sigma', A')$, we refer to the datum $(Q \sqcup Q', \Sigma \times \Sigma', A \sqcup A')$ as the \textit{product} of the two data. It will be useful in the sequel that the datum corresponding to a deleted quiver is, in fact, isomorphic to the product of a trivial datum with that of a contracted quiver.

\begin{lem}\label{product data decomposition}
Assume given $v_i \in Q^{mut}$. Then there is $c_i^\vee\in N_i$ with $(c_i,c_i^\vee)=1$ and $(c_k,c_i^\vee)=0$ for every mutable $v_k$ of $Q^{-i}$. For any such $c_i^\vee$, put $B_i:=M_i\cap\ker c_i^\vee$. Then $M_i=\Z c_i\oplus B_i$, and $(Q^{-i},\Sigma_{-i},A_{-i})$ is isomorphic to the product of the trivial datum $(\Z c_i,\{0\},\emptyset)$ with $(Q^{/i},\Sigma_{/i},A_{/i})$ on $B_i\cong M_i/\Z c_i$.
\end{lem}
 
\begin{proof}
The $c_k$ lie in $M_i$, as $v_k$ is not adjacent to $v_i$, and so does $c_i$ by definition. Since $Q$ is RFR they span a saturated sublattice of $M$, hence of $M_i$, and so one may simply take $c_i^\vee$ to be dual to $c_i$. Then $c_k\in B_i$ for all $k$, and these are the columns of $Q^{/i}$ by definition. But each ray $\theta_k^\vee$ of $\Sigma_{-i}$ has $v_k$ not adjacent to $v_i$ also by definition, and so $\Sigma_{-i}\subset c_i^\perp=B_i^\vee$. This completes the proof.
\end{proof}

Now we may give the following recursive geometric understanding of freezing, deletion, and contraction of quivers in terms of cluster models.

\begin{prop}\label{quiver geometry}
    Assume given a cluster datum $(Q, \Sigma, A)$.
    
    \begin{enumerate}
        \item For $S \subset Q^{mut}$ a set of vertices we have $$\mathcal{A}(Q)_{\Sigma, A}\setminus V(\prod_{i \in S}x_i) \cong \mathcal{A}(Q^{\backslash S})_{\Sigma_{\backslash S}, A_{\backslash S}}.$$
        \item For $\sigma_I \in \Sigma$ such that $I \cap A = \emptyset$, the $D_i$ with $i \in I$ intersect transversally and we have $$\bigcap_{i \in I} D_i = \mathcal{A}(Q^{-I})_{\Sigma_{-I}, A_{-I}}.$$
        \item For $i \in \Sigma(1) \setminus A$ we have $$D_i \cong \mathbb{C}^*_{y_i} \times \mathcal{A}(Q^{/i})_{\Sigma_{/i}, A_{/i}},$$ where $y_i = x^{c_i}$ is as in the previous proposition. Therefore $$\mathcal{A}(Q^{/i})_{\Sigma_{/i}, A_{/i}} \cong H_i \subset \mathcal{A}(Q)_{\Sigma, A}.$$
    \end{enumerate}
\end{prop}

\begin{proof}
    For $(1)$ we note that $X_\Sigma \setminus V(\prod_{i\in S} x_i) = X_{\Sigma_{\backslash S}}$ by definition. The claim then follows from the fact that affine blowups commute with open immersions, and that the centers of the $H_i$ with $i \in S$ lie in $V(\prod_{i \in S} x_i)$.\\

    For $(2)$ we induct on $|A|$. When $A = \emptyset$ this is the identity $X_{\Sigma_{-I}} = \bigcap_{i\in I}V(x_i)\subset X_\Sigma$. Inducting, assume given $j \in A$, and write $\A^-:=\mathcal \A(Q)_{\Sigma,A\setminus \{j\}}$. Further, let $D^-_I:=\bigcap_{i \in I} D^-_i\subset \A^-$. By the hypothesis we know that $D_I^-\cong\mathcal A(Q^{-I})_{\Sigma_{-I},(A\setminus \{j\})_{-I}}$. Further, one has $\mathcal \A(Q)_{\Sigma,A}=\mathrm{Bl}^{\mathrm{aff}}_{H_j}\A^-$ by definition, and we see that $\bigcap_{i \in I} D_i$ is the preimage of $D_I^-$ since $D_i$ is the zero locus of $x_i$ and by the paragraph after Proposition \ref{chart lemma}. Therefore the claim would follow if one could commute the affine blowup along $H_j$ past the inclusion of $D_I^-$ into $\A^-$, which follows from Corollary \ref{transversality}.\\

    For $(3)$, we know from $(2)$ that $D_i\cong\mathcal A(Q^{-i})_{\Sigma_{-i},A_{-i}}$. By Lemma \ref{product data decomposition}, before any blowups $X_{\Sigma_{-i}}=\mathbb C^*_y\times X_{\Sigma_{/i}}$, where $y:=x^{c_i}$. For $k\in A_{-i}$ we have $c_k\in B_i$, so $V(x_k)$ and $\{x^{c_k}=-1\}$ are $\mathbb C^*_y$ times the corresponding loci in $X_{\Sigma_{/i}}$. Therefore Proposition 2.28 implies that each center in $X_{\Sigma_{-i}}$ is $\mathbb C^*_y\times H'_k$, where $H'_k$ is the corresponding center in $X_{\Sigma_{/i}}$. Since affine blowups commute with products this completes the proof.
\end{proof}

\begin{rem}\label{GV duality remark}
    For $K\subset K'\subset \Sigma(1)\setminus A$ with $\sigma_{K'}\in \Sigma$, write $k_{K\subset K'}\colon D_{K'}\hookrightarrow D_K$, where $D_\emptyset:=\mathcal \A(Q)_{\Sigma,A}$. Then Proposition \ref{chart lemma} implies that $D_{K'}$ is cut out of $D_K$ by a regular sequence of regular functions, and so its normal bundle is trivial. Grothendieck-Verdier duality then gives
$$(k^*_{K\subset K'})^L\simeq k_{K\subset K'*}[-|K'\setminus K|].$$

Similarly, since $E_i=V(x_i)$ for $i\in A$ we have that $(\iota_{E_i}^*)^L\simeq\iota_{E_i*}[-1]$.
\end{rem}

The preceding propositions suggest a loose slogan: that cluster models are to (smooth, quasi-affine) toric varieties as (RFR acyclic) cluster varieties are to tori. In particular, $(1)$ and $(2)$ of the above show that the usual recursive structure of toric varieties, wherein subfans produce open subvarieties/complements of certain boundary divisors and stars of rays produce divisors, persists. $(3)$ is new to cluster models: indeed, the center $H_i$ is not a toric subvariety of $X_\Sigma$. Nevertheless, the analogous subvariety of a cluster model $\mathcal{A}(Q)_{\Sigma, A}$ \textit{is} again a cluster model.\\

Cementing the slogan, we show now that cluster models are partial log Calabi-Yau compactifications of acyclic cluster varieties.

\begin{prop}\label{CY prop}
    Let $\A = \mathcal{A}(Q)_{\Sigma, A}$ be a cluster model, and set $D := \sum_{\rho_j \in \Sigma(1) \setminus A} D_j$. Then $D$ is a simple normal crossings divisor, the form $\omega := \bigwedge_{v \in V} dlog x_v$ on the standard torus $(\mathbb{C}^*)^V$ induced from $X_\Sigma$ extends to a meromorphic volume form on $\A$ with simple poles along $D$, and $(\A, D)$ is a log Calabi-Yau pair.
\end{prop}

By combining this with $(1)$ of Proposition \ref{quiver geometry} we learn that any cluster model $\mathcal{A}(Q)_{\Sigma, A}$ can be identified with a log CY compactification of the all-active-rayed cluster model $\mathcal{A}(Q^{\backslash (\Sigma(1) \setminus A)})_{\Sigma_{\backslash (\Sigma(1) \setminus A)}, A}$. In the sequel we refer to the latter as the \textit{reduced datum} corresponding to $(Q, \Sigma, A)$: i.e., for $S = \Sigma(1) \setminus A$ we write $$(Q^{red}, \Sigma^{red}, A^{red}) := (Q^{\setminus S},\Sigma_{\setminus S},A),$$

\noindent and we write $\A^{red}_{\Sigma, A}$ to denote the corresponding cluster model. Therefore, for example, we see that any $\A^{red}_{\Sigma_Q, A}$ agrees with the cluster variety $\mathcal{A}(Q^{\backslash S})$.

\begin{proof}
The SNC-ness claim follows from Proposition \ref{chart lemma}, since it tells us that locally on a collection of charts of $\mathcal{A}(Q)_{\Sigma, A}$ the divisors in question identify with (products with) coordinate axes. To show the rest we proceed by induction on $|A|$, with base case $A = \emptyset$ the known case of toric varieties rel their toric boundaries.\\

Inducting, assume given $i\in A$. Consider the cluster model $\A^-:=\mathcal{A}(Q)_{\Sigma,A\setminus\{i\}}$ with center $H_i\subset D_i$ (here $H_i$ and $D_i$ denote the strict transforms of $H_i$ and $D_i$ in $\A^-$, as above). By induction we know that $div(\omega)=-\sum_{\rho_j \in \Sigma(1) \setminus (A \setminus\{i\})} D_j$. Let $\pi$ be the (non-affine) blowup along $H_i$, and let $F_i$ denote its exceptional divisor. Since $H_i$ is smooth of codimension $2$, we have that $div (\pi^*\omega)=\pi^*div(\omega)+F_i$. Further, from Proposition \ref{chart lemma} we see that $H_i\not\subset D_j$ for $j\neq i$. Therefore $\pi^*D_i=\tilde D_i+F_i$ and $\pi^*D_j=\tilde D_j$ (where tildes denote strict transforms, contrary to our usual abuse of notation for $D_i$), and so
$$div(\pi^*\omega)=-\sum_{\rho_j \in \Sigma(1) \
\setminus (A \setminus \{i\})} \tilde D_j-F_i+F_i=-\sum_{\rho_j \in \Sigma(1) \
\setminus (A \setminus \{i\})} \tilde D_j.$$
Removing $\tilde D_i$ gives $\mathcal{A}(Q)_{\Sigma, A}$, on which we may identify the remaining portion of $\tilde{D}_j$ with the divisor $D_j$. Therefore we see that $div(\pi^* \omega|_{\mathcal{A}(Q)_{\Sigma, A}}) = -\sum_{\rho_j \in \Sigma(1) \setminus A} D_j$. Since the blowdown from $\mathcal{A}(Q)_{\Sigma, A}$ onto $X_\Sigma$ is clearly an isomorphism onto the latter's central torus the claim follows.
\end{proof}

Finally, we note that boundary divisors on a cluster model are in fact log CY \textit{with volume forms induced from boundary maps.} That is, assume given $\sigma_K \in \Sigma$ such that the vertices in $K$ are inactive. Then, from $(2)$ of Proposition \ref{quiver geometry} we learn that $D_K$ is itself a cluster model with respect to the cluster datum $(Q^{-K},\Sigma_{-K},A_{-K})$. In particular, its own inactive divisors $D_j$ are the intersections $D_j\cap D_K$ for inactive $j \notin K$ satisfying $\sigma_{j \cup K} \in \Sigma$. Further, one can compute that its volume form is in fact the residue of $\omega$:
$$Res_{D_K}\omega=\pm\bigwedge_{v\notin K}d\log x_v.$$

In short, therefore, cluster models and their inactive boundary divisors behave like toric varieties and their toric boundary strata: every boundary stratum is again one of them, and residues of their natural volume forms induce the natural forms on smaller such varieties.\\

We end this section with a computation of the Chow groups and, as a corollary, the $0$-th K-theory of cluster models. We find that they are trivial, a result which we encourage the reader to compare to (and which generalizes the acyclic case of) the main theorem of \cite{hlavinka2026chow}.

\begin{prop}\label{chow vanishing}
    Assume given a cluster datum $(Q, \Sigma, A)$. Then $CH_{m + n}(\mathcal{A}(Q)_{\Sigma, A}) = \mathbb{Z}$ and $CH_i(\mathcal{A}(Q)_{\Sigma, A}) = 0$ for all $i < m +n$, where $CH_i$ denotes the $i$-th integral (homological) Chow group. This implies that $K_0(\mathcal{A}(Q)_{\Sigma, A}) = 
    \mathbb{Z}$. 
\end{prop}

\begin{proof}
    By excision, inducting on both $m + n$ and $|A|$. For $A = \emptyset$ this follows from the fact that $X_\Sigma$ is an open subset of affine space.\\

    Inducting, we note that $\mathcal{A}(Q)_{\Sigma, A \setminus \{j\}} \setminus D_j$ and $\mathcal{A}(Q)_{\Sigma, A} \setminus E_j$ are isomorphic and that $E_j = V(x_j)$ in the full cluster model. Therefore the excision sequence for the inclusion $E_j \hookrightarrow \mathcal{A}(Q)_{\Sigma, A}$ takes the form $$CH_i(E_j) \rightarrow CH_i(\mathcal{A}(Q)_{\Sigma, A}) \rightarrow CH_i(\mathcal{A}(Q)_{\Sigma, A \setminus \{j\}} \setminus D_j) \rightarrow 0.$$

    $E_j \cong \mathbb{A}^1 \times H_j$ by construction. Therefore the induction hypothesis implies that the only possible contributions to $CH_*(\mathcal{A}(Q)_{\Sigma, A})$ come from the fundamental classes of $\mathbb{A}^1 \times H_j$ and of $\mathcal{A}(Q)_{\Sigma, A \setminus \{j\}} \setminus D_j$. But $E_j$ is principal, and so the former dies, implying the Chow portion.\\

    The claim about $K_0$ follows from the fact that the coniveau filtration on $G_0$ has associated graded given by subquotients of $CH_i$, since cluster models are smooth and so $G_0 \cong K_0$.
\end{proof}

\subsection{Braid and double Bott--Samelson Varieties.}

Though our focus is not braid varieties, it will help (especially in Section 5) for us to have algorithms which pertain to them at hand for computation. We therefore briefly recall their theory, which also gives us a chance to fix otherwise useful notation. This subsection and the following one will not be used explicitly until Section 5.\\

Fix standard and opposite Borels $B$ and $B_-$, respectively, in $G := GL_n$. Writing $B = \langle e_1, ..., e_n \rangle$, we also fix ``cyclic shifts" of $B$ by setting $B_i = \langle e_i, ..., e_n, e_1, ..., e_{i - 1}\rangle$.\\

Let $\lambda$ denote a partition whose Young diagram is contained in a $k \times (n-k)$ rectangle. We write $X_{\lambda, i}^\circ$ to denote the Schubert cell in $Gr(k,n)$ corresponding to the Young diagram $\lambda$ \textit{with respect to $B_i$}, and $X_{\lambda, i}$ to denote the corresponding Schubert variety. Given a pair of Bruhat related permutations $u \leq w$ in $S_n$, we write $R_{u,w}^\circ$ and $R_{u,w}$ to denote the corresponding open and closed (resp.) Richardson varieties in the flag variety $Fl(n) = G/B$.\\

We will use $\beta \in$ $Br_n^+$ to denote a positive braid or braid word on $n$-strands with Demazure product $\delta(\beta) = w_0$. We will never consider non-positive braid words, and we will often (likely tacitly) fix particular braid words $\beta = \sigma_{i_1} \cdots \sigma_{i_r}$ for a given braid $\beta$ throughout. For example, $\Delta$ will always denote the half-twist in the form $\Delta = \sigma_1(\sigma_2 \sigma_1)\cdots(\sigma_{n-1}\dots \sigma_1)$.\\

Let $i \in [n-1]$. Following \cite{CGGLSS}, we define the \textit{braid matrix} $B_i(z)$  associated to the twist $\sigma_i$ by
\begin{equation} B_i(z) = I_{i-1} \oplus B(z) \oplus I_{n - i - 1},
\end{equation}

\noindent where $B(z) = \begin{bmatrix}
    z      & -1  \\
    1       & 0
\end{bmatrix}.$
Therefore $B_i(z)$ is the matrix which is the identity away from a $2 \times 2$ subblock equal to $B(z)$, and with upper left corner at $(i,i)$. Given $\beta = \sigma_{i_1}  \in Br^+_n$ and $(z_1, ..., z_r) \in \mathbb{C}^r$, we write 
\begin{equation}
    B_\beta(z) := B_{i_1}(z_1) \cdots B_{i_r}(z_r).
\end{equation}

\begin{defn}
    The \textit{braid variety} $X(\beta)$ associated to a positive braid word $\beta = \sigma_{i_1} \cdots \sigma_{i_r}$ is the subvariety of $\mathbb{C}^r$ on which $B_\beta(z)$ is upper anti-triangular, i.e.
    \begin{equation}
        X(\beta) := \{(z_1, ..., z_r) \in \mathbb{C}^r\:|\:w_0B_\beta(z) \in B\}.
    \end{equation}
\end{defn}

$X(\beta)$ is an affine variety cut out by the vanishing of certain generalized minors of $B_\beta(z)$. By \cite{GLSBSI} Corollary 5.8. and Theorem 5.9., it is also naturally a really full rank locally acyclic cluster variety.

\begin{defn}
    A \textit{double Bott--Samelson variety} is a braid variety of the form $X(\Delta \beta)$.
\end{defn}

We recall an algorithm, introduced in \cite{CGGLSS} but interpreted in the following from in \cite{castronovo2024cluster}, which constructs a quiver $Q_\beta$ for the cluster structure on any braid double Bott--Samelson variety $X(\Delta \beta)$.

\begin{con}\label{con:DES Construction}
    Fix a positive braid word $\beta = \sigma_{i_1}\cdots\sigma_{i_j}$ on $n$ strands, and set $\beta' = \sigma_{i_1}\cdots\sigma_{i_{j-1}}$. We associate a vertex colored quiver $Q_\beta$ to $\beta$ as follows. First, a vertex will be colored $i$ when it corresponds to an appearance of $\sigma_i$ in $\beta$. Further, an $i$-colored vertex is frozen if and only if it corresponds to the rightmost appearance of $\sigma_i$ in $\beta$.\\
    
    The construction then proceeds iteratively. That is, assume that the vertex-colored quiver $Q_{\beta'}$ has been constructed. Then to construct $Q_\beta$:

    \begin{enumerate}
        \item Thaw the frozen vertex $\bullet$ of $Q_{\beta'}$ corresponding to the rightmost appearance of $\sigma_{i_j}$ in $\beta'$ (if any exist).

        \item Add a new frozen vertex $\square$ and an edge $\bullet \rightarrow \square$. $\square$ will correspond to the rightmost appearance of $\sigma_{i_j}$ in $\beta$, and so will be colored $i_j$.

        \item For $r \in \{i_j - 1, i_j + 1\}$, let $\triangle$ be the $r$-colored vertex in $Q_{\beta'}$ corresponding to the rightmost appearance of $\sigma_r$. If $\triangle$ already has an edge to a vertex of color $i_j$, do nothing; otherwise, add an edge $\triangle \rightarrow \bullet$ to the newly thawed vertex $\bullet$ above.
    \end{enumerate}
\end{con}

\begin{thm}[\cite{CGGLSS} Corollary 4.45.]
    We have an isomorphism $$X(\Delta \beta) \simeq \mathcal{A}(Q_\beta).$$
\end{thm}

Most positroid varieties are not double Bott--Samelson varieties. We recall a different algorithm which constructs a quiver for the cluster structure on any positroid variety in the next subsection.

\begin{exmp}
These examples follow the ones given in Example 4.39 of \cite{castronovo2024cluster} almost verbatim, since they are ones of essential importance to us.

\begin{enumerate}
    \item For the word $\beta = \sigma_1^r$, the quiver $Q_{\beta}$ is $A_{r-1}$.

    \item For the word $\beta = \sigma_1^{2}\sigma_2\sigma_1^{2}\sigma_2$, the quiver $Q_{\beta}$ is the following.

\begin{center}
\begin{tikzpicture}
\node at (1,0) {$\circ$};
\draw [->] (1.25, 0) -- (1.75,0);
\node at (2,0) {$\circ$};
\draw [->] (2.25, 0) -- (2.75,0);
\node at (3,0) {$\circ$};
\draw [->] (3.25, 0) -- (3.75,0);
\node at (4,0) {$\square$};

\node at (3,-1) {$\circ$};
\draw [->] (3.25, -1) -- (3.75,-1);
\node at (4,-1) {$\square$};

\draw [->] (4, -0.25) -- (3.15, -0.85);
\draw[->] (2.85, -0.85) -- (2.15, -0.25);
\end{tikzpicture}
\end{center}
This quiver is clearly acyclic. In fact, the mutable portion of $Q_\beta$ is of type $D_4$. More generally, for $\beta = \sigma_1^{n-2}\sigma_2\sigma_1^{2}\sigma_2$, $Q_{\beta}$ is a quiver of type $D_n$ with two frozen vertices.

\item For the words $\beta = (\sigma_1\sigma_2)^{4}$ and $(\sigma_1\sigma_2)^{5}$, the $Q_{\beta}$'s are the following two quivers (resp.).
\begin{center}
\begin{tikzpicture}
\node at (0,0) {$\circ$};
\draw[->] (0.25,0)--(0.75,0);
\node at (1,0) {$\circ$};
\draw [->] (1.25, 0) -- (1.75,0);
\node at (2,0) {$\circ$};
\draw [->] (2.25, 0) -- (2.75,0);
\node at (3,0) {$\square$};

\node at (0,-1) {$\circ$};
\draw[->] (0.25,-1)--(0.75,-1);
\node at (1,-1) {$\circ$};
\draw [->] (1.25, -1) -- (1.75,-1);
\node at (2,-1) {$\circ$};
\draw [->] (2.25, -1) -- (2.75,-1);
\node at (3,-1) {$\square$};

\draw[->] (0, -0.85) -- (0, -0.15);
\draw[->] (1, -0.85) -- (1, -0.15);
\draw[->] (2, -0.85) -- (2, -0.15);

\draw[->] (0.75, -0.25)--(0.25, -0.75);
\draw[->] (1.75, -0.25)--(1.25, -0.75);
\draw[->] (2.75, -0.25)--(2.25, -0.75);
\end{tikzpicture}

\begin{tikzpicture}
\node at (0,0) {$\circ$};
\draw[->] (0.25,0)--(0.75,0);
\node at (1,0) {$\circ$};
\draw [->] (1.25, 0) -- (1.75,0);
\node at (2,0) {$\circ$};
\draw [->] (2.25, 0) -- (2.75,0);
\node at (3,0) {$\circ$};
\draw[->] (3.25, 0) -- (3.75, 0);
\node at (4,0) {$\square$};

\node at (0,-1) {$\circ$};
\draw[->] (0.25,-1)--(0.75,-1);
\node at (1,-1) {$\circ$};
\draw [->] (1.25, -1) -- (1.75,-1);
\node at (2,-1) {$\circ$};
\draw [->] (2.25, -1) -- (2.75,-1);
\node at (3,-1) {$\circ$};
\draw[->] (3.2, -1) -- (3.8, -1);
\node at (4,-1) {$\square$};

\draw[->] (0, -0.85) -- (0, -0.15);
\draw[->] (1, -0.85) -- (1, -0.15);
\draw[->] (2, -0.85) -- (2, -0.15);
\draw[->] (3, -0.85) -- (3, -0.15);

\draw[->] (0.75, -0.25)--(0.25, -0.75);
\draw[->] (1.75, -0.25)--(1.25, -0.75);
\draw[->] (2.75, -0.25)--(2.25, -0.75);
\draw[->] (3.75, -0.25)--(3.25, -0.75);
\end{tikzpicture}
\end{center}
\end{enumerate}

It is not \textit{a priori} obvious that the latter two examples are mutation equivalent to quivers of type $E_6$ and $E_8$ resp., but Proposition 7.5. of \cite{SSBW} explains how they are. In fact, $\mathcal{A}(A_{n-1}) \times (\mathbb{C}^*)^n \cong \Pi_{2,n + 2}$, $\mathcal{A}(D_4) \times (\mathbb{C}^*)^3 \cong \Pi_{3,6}$, and $\mathcal{A}(E_8) \times (\mathbb{C}^*)^5 \cong \Pi_{3,8}$.
\end{exmp}

\subsection{Positroid Varieties.}

The positroid stratification of the Grassmannian $Gr(k,n)$ is a reasonable coarsening of its matroid stratification whose closed strata still enjoy the nice smoothness properties of closed Schubert and Richardson varieties (e.g., they all have rational singularties).

\begin{defn}
    Let $Gr(k,n) = \bigsqcup_{(\lambda_1, ..., \lambda_n)} \bigcap_i X_{\lambda_i, i}^\circ$ be the decomposition of $Gr(k,n)$ into cyclic intersections of Schubert cells. A non-empty piece of this decomposition is called a \textit{positroid cell} or a (in some sources, \textit{open}\footnote{Indeed, we warn the reader that the convention of referring to positroid cells as positroid varieties without the modifier ''open" is the less common one in the literature. We do so throughout nonetheless because we will consider closed positroid varieties only vanishingly often.}) \textit{positroid variety}.
\end{defn}

It is non-obvious when a positroid cell corresponding to a tuple of partitions $(\lambda_1, ..., \lambda_n)$ is non-empty. Positroid varieties are bijectively parametrized by several other pieces of combinatorial information, however. For example, they biject with pairs $u \leq w$ of Bruhat related permutations such that $w$ is $k$-Grassmannian. For this reason we will call such pairs $(u, w)$ \textit{positroid pairs.}

\begin{thm}[\cite{KLS} Theorem 5.9.]\label{projected}
    Let $\pi: Fl(n) \rightarrow Gr(k,n)$ denote the projection map. Then, for every positroid cell $X \subset Gr(k,n)$, there exists a unique pair of permutations $u \leq w$ in $S_n$ such that $w$ is $k$-Grassmannian and $\pi|_{R_{u,w}^\circ}$ maps isomorphically onto $X$.
\end{thm}

The positroid variety associated to a positroid pair $u \leq w$ is denoted $\Pi_{u,w}$. The only exception to this rule is that we write $\Pi_{k,n} = \bigcap_i X_{\emptyset, i}^\circ$ to denote the open positroid cell in $Gr(k,n)$.\\

Given a reduced expression for a permutation $w$, let $\beta(w)$ denote a positive braid lift of said word. Then, the remarkable works \cite{CGGLSS} and \cite{GLSBSI} show that every positroid variety is a braid variety, with the former exhibiting an isomorphism of the form $$\Pi_{u,w} \cong X(\beta(u^{-1}w_0)\beta(w_0ww_0)).$$

In particular, all positroid varieties are really full rank locally acyclic cluster varieties.\\

We now begin to recall \cite{GaLa}'s construction of a cluster structure on any positroid variety $\Pi_{u,w}$. To do so we must introduce a collection of decorated Young diagrams which also bijectively parametrize positroid pairs $u \leq w$.

\begin{defn}
    Let $\lambda$ be a Young diagram inside a $k \times (n-k)$ rectangle. A \textit{Le-diagram} $D$ for $\lambda$ is a filling of the boxes of $\lambda$ with dots satisfying the following property: if a box $b$ is both directly below and directly to the right of a dot, then $b$ also has a dot.
\end{defn}

An example of such a $D$ is given in the middle of Figure \ref{Le depiction}. Note that the Le condition implies that the $D$ depicted there \textit{must} have a dot in the bottom right box.\\

\def\borderInnerSep{2pt}
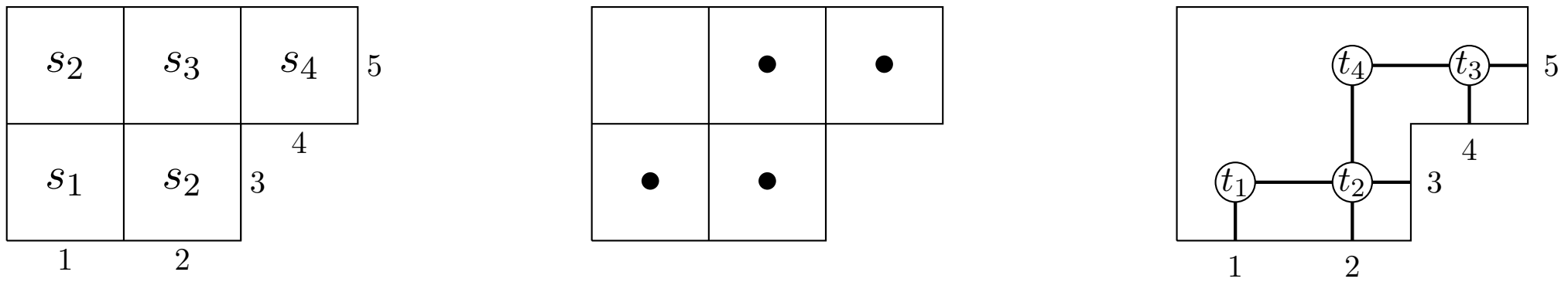
\begin{figure}
\begin{tikzpicture}
\draw (0,0)-- (0,2)--(3,2)--(3,1)--(2,1)--(2,0)--(0,0);
\draw (1,0)--(1,2);
\draw (2,0)--(2,2);
\draw (0,1)--(2,1);

\def\nodescl{0.7}
\node[anchor=north, inner sep=\borderInnerSep] (A) at (0.5,0) {\scalebox{\nodescl}{$1$}};
\node[anchor=north, inner sep=\borderInnerSep] (A) at (1.5,0) {\scalebox{\nodescl}{$2$}};
\node[anchor=north, inner sep=\borderInnerSep] (A) at (2.5,1) {\scalebox{\nodescl}{$4$}};
\node[anchor=west, inner sep=\borderInnerSep] (A) at (2,0.5) {\scalebox{\nodescl}{$3$}};
\node[anchor=west, inner sep=\borderInnerSep] (A) at (3,1.5) {\scalebox{\nodescl}{$5$}};
\node (A) at (0.5,0.5) {$s_1$};
\node (A) at (1.5,0.5) {$s_2$};
\node (A) at (0.5,1.5) {$s_2$};
\node (A) at (1.5,1.5) {$s_3$};
\node (A) at (2.5,1.5) {$s_4$};
\begin{scope}[shift={(5,0)}]
\draw (0,0)-- (0,2)--(3,2)--(3,1)--(2,1)--(2,0)--(0,0);
\draw (1,0)--(1,2);
\draw (2,0)--(2,2);
\draw (0,1)--(2,1);
\def\nodescl{0.7}
\node (A) at (0.5,0.5) {$\bullet$};
\node (A) at (1.5,0.5) {$\bullet$};
\node (A) at (1.5,1.5) {$\bullet$};
\node (A) at (2.5,1.5) {$\bullet$};
\end{scope}

\begin{scope}[shift={(10,0)}]

\def\nodescl{0.7}
\draw (0,0)-- (0,2)--(3,2)--(3,1)--(2,1)--(2,0)--(0,0);
\node[anchor=north] (A) at (0.5,0) {\scalebox{\nodescl}{$1$}};
\node[anchor=north] (A) at (1.5,0) {\scalebox{\nodescl}{$2$}};
\node[anchor=north] (A) at (2.5,1) {\scalebox{\nodescl}{$4$}};
\node[anchor=west] (A) at (2,0.5) {\scalebox{\nodescl}{$3$}};
\node[anchor=west] (A) at (3,1.5) {\scalebox{\nodescl}{$5$}};
\draw[line width=0.8pt] (0.5,0.5)--(2,0.5);
\draw[line width=0.8pt] (1.5,1.5)--(3,1.5);
\draw[line width=0.8pt] (1.5,1.5)--(1.5,0);
\draw[line width=0.8pt] (0.5,0.5)--(0.5,0);
\draw[line width=0.8pt] (2.5,1.5)--(2.5,1);
\node[draw,circle,scale=0.8,inner sep=0pt,fill=white] (A) at (0.5,0.5) {$t_1$};
\node[draw,circle,scale=0.8,inner sep=0pt,fill=white] (A) at (1.5,0.5) {$t_2$};
\node[draw,circle,scale=0.8,inner sep=0pt,fill=white] (A) at (2.5,1.5) {$t_3$};
\node[draw,circle,scale=0.8,inner sep=0pt,fill=white] (A) at (1.5,1.5) {$t_4$};
\end{scope}

\end{tikzpicture}
  \caption{\label{Le depiction} The Young diagram $\lambda$, Le-diagram $D$, and graph $G(D)$ corresponding to $(u,w) = (s_2,s_2s_1s_4s_3s_2)$.}
\end{figure}

We quickly remind the non-combinatorially oriented reader of a procedure which will be useful to explain our interest in Le-diagrams.

\begin{con}\label{**}
    There is a bijection from Young diagrams $\lambda$ inside a $k \times (n-k)$ rectangle to $k$-Grassmannian permutations $w \in S_n$ described as follows. Given such a Young diagram $\lambda$, write the transposition $s_1$ in the bottom left entry of $\lambda$. Iteratively, in every box above or to the right of one marked with $s_i$, write the tranposition $s_{i+1}$. For example, the top right box of $\lambda$ will contain the transposition $s_{a(\lambda) + \ell(\lambda) +1}$, where $a$ and $\ell$ denote the arm and leg of $\lambda$. Then, reading the boxes in $\lambda$ from right to left, starting with the bottom row and going to the top, produces a reduced expression for a $k$-Grassmannian word $w_\lambda$.\\
    
    The result of the ``labeling" portion of this procedure is illustrated in the first item of Figure \ref{Le depiction}; the word which one should ``read off" in this example is $s_2s_1s_4s_3s_2$. That this procedure gives a bijection between Young diagrams and $k$-Grassmannian words is standard.
\end{con}

Continuing to write $w_\lambda$ for the permutation corresponding to $\lambda$, we may now establish the following

\begin{prop}
    There is a bijection between Le-diagrams of shape $\lambda$ and positroid pairs $u \leq w_\lambda$.
\end{prop}

\begin{proof}
    Given a Le-diagram $D$ for $\lambda$, one recovers a reduced expression for some $u \leq w_\lambda$ via the reading portion of Construction \ref{**}, save that one \textit{skips the transpositions in any box with a dot}. For example, the Le-diagram in Figure \ref{Le depiction} corresponds to the positroid pair $u = s_2 \leq w_\lambda = s_2s_1s_4s_3s_2$.\\

    Conversely, to associate a Le-diagram to a positroid pair $u \leq w_\lambda$, one finds the unique rightmost reduced expression for $u$ inside of the reading procedure expression for $w_\lambda$. The Le-diagram $D = D_{u,w_\lambda}$ corresponding to $u \leq w_\lambda$ is then given by taking the Young diagram of the partition $\lambda$ corresponding to $w_\lambda$, doing the labeling portion of Construction \ref{**} in $\lambda$, and then marking all of the letters which are not used by this rightmost reduced expression for $u$ with a dot.\\

    We leave it to the reader to check that the two directions are mutually inverse.
\end{proof} 

We now follow \cite{GaLa} in describing a way to associated a quiver to any positroid cell. To start, to a Le-diagram $D$ for $\lambda$ we associate a planar graph $G(D)$ embedded in a disk (whose boundary we identify with the boundary of $\lambda$) by ``drawing a southeast pointing hook from each vertex." More precisely, from each dot we draw a line segment going east and south until it hits the boundary of $\lambda$. The vertices of $G(D)$ are the dots of $D$ and the places where the line segments hit the boundary; two vertices are connected by an edge if they are connected by a line segment. Finally, we label each of the interior vertices with entries $t_i$. This procedure is illustrated in the final entry of Figure \ref{Le depiction}.\\

We will refer to the face of $G(D)$ whose top-left vertex is $t_r$ as $F_r$. The top-left face of $G(D)$ will be called $F_0$.

\begin{con}\label{Pasha Construction}
    Let $D$ be a Le-diagram with corresponding planar graph $G(D)$. We construct a quiver $Q_D$ whose vertices are the faces of $G(D)$ excluding $F_0$. By abuse of notation, we also refer to said vertices by the names of the corresponding faces in what follows.\\
    
    To do so, for each $F_r$ we draw edges to adjacent faces $F_i$ (save $F_0$) depending on the local structure of $G(D)$ near the vertex labeled $t_r$. If $F_r$ only borders a single other face $F_c$, at $t_r$ then we draw an edge from $F_r$ to $F_c$. If $F_r$ is to the right of a face $F_b$ and they are both below a face $F_c$ at $t_r$, we draw an edge from $F_r$ to $F_b$. If $F_r$ is below a face $F_a$ and they are both to the right of a face $F_c$ at $t_r$, we draw an edge from $F_r$ to $F_a$. If $F_r$ is the bottom right entry in a set of four faces around $t_r$, we draw edges from $F_r$ to the top right and bottom left ones and an edge from the top left one to $F_r$. Finally, we set $F_i$ to be frozen if it borders the bottom right of the disk and to be mutable otherwise. $Q_D$ is defined to be the resulting quiver. An illustration of this procedure is given in Figure \ref{local moves}.
\def\leop{0.4}
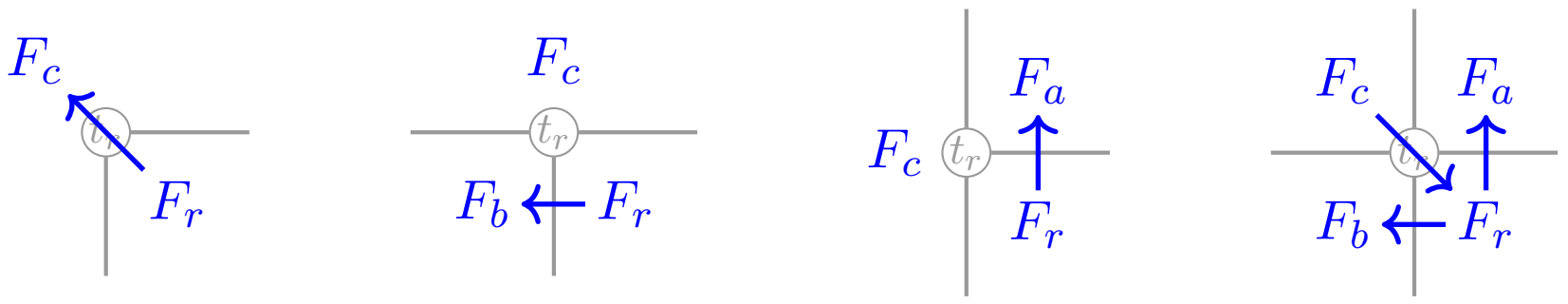
\begin{figure}
\begin{tikzpicture}
\def\qlw{1pt}
\def\qop{1}

\def\Finnersep{2pt}
\def\stp{3}

\node(A) at (0,0){
\begin{tikzpicture}
\draw[line width=0.8pt,opacity=\leop] (0,0)--(1,0);
\draw[line width=0.8pt,opacity=\leop] (0,0)--(0,-1);
\node[draw=white,circle,scale=0.8,inner sep=0pt,fill=white] (A) at (0,0) {\textcolor{white}{$t_r$}};
\node[draw,circle,scale=0.8,inner sep=0pt,fill=white,opacity=\leop] (A) at (0,0) {$t_r$};
\node[blue,inner sep=\Finnersep] (Fr) at (0.5,-0.5) {$F_r$};
\node[blue,inner sep=\Finnersep] (Fc) at (-0.5,0.5) {$F_c$};
\draw[->,blue,opacity=\qop,line width=\qlw] (Fr)--(Fc);
\end{tikzpicture}
};

\node(A) at (\stp,0){
\begin{tikzpicture}
\draw[line width=0.8pt,opacity=\leop] (-1,0)--(1,0);
\draw[line width=0.8pt,opacity=\leop] (0,0)--(0,-1);
\node[draw=white,circle,scale=0.8,inner sep=0pt,fill=white] (A) at (0,0) {\textcolor{white}{$t_r$}};
\node[draw,circle,scale=0.8,inner sep=0pt,fill=white,opacity=\leop] (A) at (0,0) {$t_r$};
\node[blue,inner sep=\Finnersep] (Fr) at (0.5,-0.5) {$F_r$};
\node[blue,inner sep=\Finnersep] (Fc) at (0,0.5) {$F_c$};
\node[blue,inner sep=\Finnersep] (Fb) at (-0.5,-0.5) {$F_b$};
\draw[->,blue,opacity=\qop,line width=\qlw] (Fr)--(Fb);
\end{tikzpicture}
};

\node(A) at (2*\stp,0){
\begin{tikzpicture}
\draw[line width=0.8pt,opacity=\leop] (0,0)--(1,0);
\draw[line width=0.8pt,opacity=\leop] (0,1)--(0,-1);
\node[draw=white,circle,scale=0.8,inner sep=0pt,fill=white] (A) at (0,0) {\textcolor{white}{$t_r$}};
\node[draw,circle,scale=0.8,inner sep=0pt,fill=white,opacity=\leop] (A) at (0,0) {$t_r$};
\node[blue,inner sep=\Finnersep] (Fr) at (0.5,-0.5) {$F_r$};
\node[blue,inner sep=\Finnersep] (Fa) at (0.5,0.5) {$F_a$};
\node[blue,inner sep=\Finnersep] (Fc) at (-0.5,0) {$F_c$};
\draw[->,blue,opacity=\qop,line width=\qlw] (Fr)--(Fa);
\end{tikzpicture}
};

\node(A) at (3*\stp,0){
\begin{tikzpicture}
\draw[line width=0.8pt,opacity=\leop] (-1,0)--(1,0);
\draw[line width=0.8pt,opacity=\leop] (0,1)--(0,-1);
\node[draw=white,circle,scale=0.8,inner sep=0pt,fill=white] (A) at (0,0) {\textcolor{white}{$t_r$}};
\node[draw,circle,scale=0.8,inner sep=0pt,fill=white,opacity=\leop] (A) at (0,0) {$t_r$};
\node[blue,inner sep=\Finnersep] (Fr) at (0.5,-0.5) {$F_r$};
\node[blue,inner sep=\Finnersep] (Fa) at (0.5,0.5) {$F_a$};
\node[blue,inner sep=\Finnersep] (Fc) at (-0.5,0.5) {$F_c$};
\node[blue,inner sep=\Finnersep] (Fb) at (-0.5,-0.5) {$F_b$};
\draw[->,blue,opacity=\qop,line width=\qlw] (Fr)--(Fa);
\draw[->,blue,opacity=\qop,line width=\qlw] (Fr)--(Fb);
\draw[->,blue,opacity=\qop,line width=\qlw] (Fc)--(Fr);
\end{tikzpicture}
};

\end{tikzpicture}
\caption{\label{local moves} Local rules for constructing a quiver from a Le-diagram.}
\end{figure}

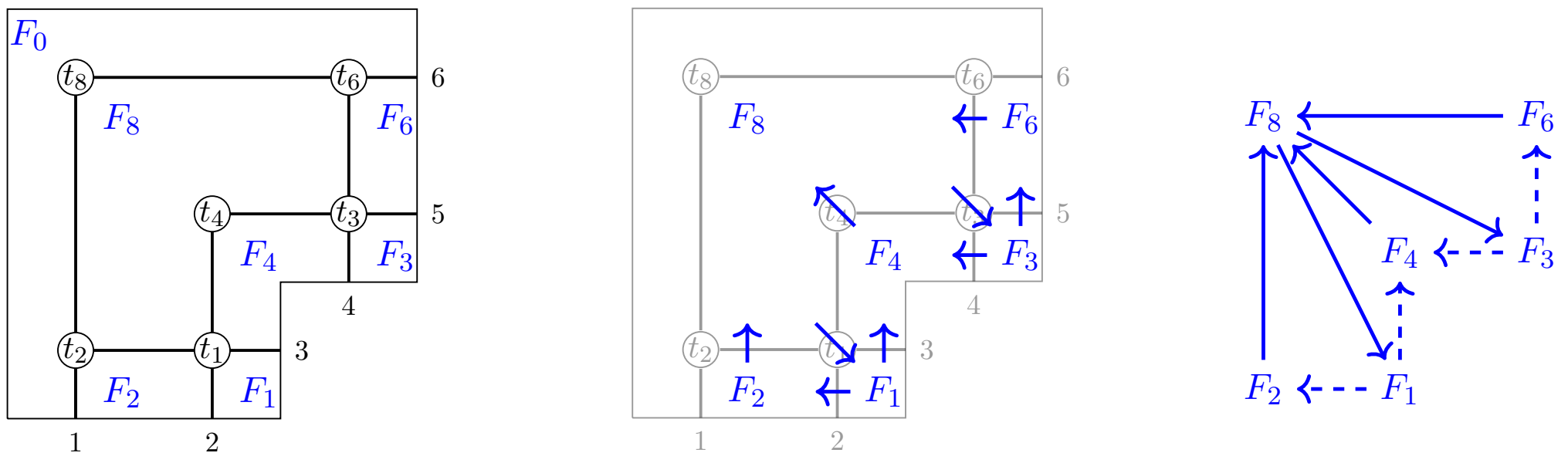
\begin{figure}
\scalebox{0.9}{
\begin{tikzpicture}

\def\scl{1.3}
\def\Fscl{0.9}
\def\Fsep{4pt}
\def\arrlen{0.5}

\node(X) at (0,0){
\begin{tikzpicture}[scale=\scl]
\def\nodescl{0.7}
\draw (0,0)-- (2,0)--(2,1)--(3,1)--(3,3)--(0,3)--(0,0);
\node[anchor=north] (A) at (0.5,0) {\scalebox{\nodescl}{$1$}};
\node[anchor=north] (A) at (1.5,0) {\scalebox{\nodescl}{$2$}};
\node[anchor=north] (A) at (2.5,1) {\scalebox{\nodescl}{$4$}};
\node[anchor=west] (A) at (2,0.5) {\scalebox{\nodescl}{$3$}};
\node[anchor=west] (A) at (3,1.5) {\scalebox{\nodescl}{$5$}};
\node[anchor=west] (A) at (3,2.5) {\scalebox{\nodescl}{$6$}};
\draw[line width=0.8pt] (0.5,0.5)--(2,0.5);
\draw[line width=0.8pt] (1.5,1.5)--(3,1.5);
\draw[line width=0.8pt] (1.5,1.5)--(1.5,0);
\draw[line width=0.8pt] (0.5,0.5)--(0.5,0);
\draw[line width=0.8pt] (2.5,2.5)--(2.5,1);
\draw[line width=0.8pt] (0.5,2.5)--(3,2.5);
\draw[line width=0.8pt] (0.5,2.5)--(0.5,0);

\node[draw,circle,scale=0.8,inner sep=0pt,fill=white] (A2) at (0.5,0.5) {$t_2$};
\node[draw,circle,scale=0.8,inner sep=0pt,fill=white] (A1) at (1.5,0.5) {$t_1$};
\node[draw,circle,scale=0.8,inner sep=0pt,fill=white] (A3) at (2.5,1.5) {$t_3$};
\node[draw,circle,scale=0.8,inner sep=0pt,fill=white] (A4) at (1.5,1.5) {$t_4$};
\node[draw,circle,scale=0.8,inner sep=0pt,fill=white] (A6) at (2.5,2.5) {$t_6$};
\node[draw,circle,scale=0.8,inner sep=0pt,fill=white] (A8) at (0.5,2.5) {$t_8$};

\node[scale=\Fscl,blue,anchor=north west,inner sep=\Fsep] (F1) at (A1.south east) {$F_1$};
\node[scale=\Fscl,blue,anchor=north west,inner sep=\Fsep] (F2) at (A2.south east) {$F_2$};
\node[scale=\Fscl,blue,anchor=north west,inner sep=\Fsep] (F3) at (A3.south east) {$F_3$};
\node[scale=\Fscl,blue,anchor=north west,inner sep=\Fsep] (F4) at (A4.south east) {$F_4$};
\node[scale=\Fscl,blue,anchor=north west,inner sep=\Fsep] (F6) at (A6.south east) {$F_6$};
\node[scale=\Fscl,blue,anchor=north west,inner sep=\Fsep] (F8) at (A8.south east) {$F_8$};
\node[scale=\Fscl,blue,anchor=south east,inner sep=\Fsep] (outer) at (A8.north west) {$F_0$};
\end{tikzpicture}
};

\node(X) at (6,0){
\begin{tikzpicture}[scale=\scl]
\def\nodescl{0.7}
\begin{scope}[opacity=\leop]
\draw (0,0)-- (2,0)--(2,1)--(3,1)--(3,3)--(0,3)--(0,0);
\node[anchor=north] (A) at (0.5,0) {\scalebox{\nodescl}{$1$}};
\node[anchor=north] (A) at (1.5,0) {\scalebox{\nodescl}{$2$}};
\node[anchor=north] (A) at (2.5,1) {\scalebox{\nodescl}{$4$}};
\node[anchor=west] (A) at (2,0.5) {\scalebox{\nodescl}{$3$}};
\node[anchor=west] (A) at (3,1.5) {\scalebox{\nodescl}{$5$}};
\node[anchor=west] (A) at (3,2.5) {\scalebox{\nodescl}{$6$}};
\draw[line width=0.8pt] (0.5,0.5)--(2,0.5);
\draw[line width=0.8pt] (1.5,1.5)--(3,1.5);
\draw[line width=0.8pt] (1.5,1.5)--(1.5,0);
\draw[line width=0.8pt] (2.5,2.5)--(2.5,1);
\draw[line width=0.8pt] (0.5,2.5)--(3,2.5);
\draw[line width=0.8pt] (0.5,2.5)--(0.5,0);
\end{scope}
\node[circle,scale=0.8,inner sep=0pt,fill=white] (A2) at (0.5,0.5) {$\strut$};
\node[circle,scale=0.8,inner sep=0pt,fill=white] (A2) at (0.5,0.5) {$\strut$};
\node[circle,scale=0.8,inner sep=0pt,fill=white] (A1) at (1.5,0.5) {$\strut$};
\node[circle,scale=0.8,inner sep=0pt,fill=white] (A3) at (2.5,1.5) {$\strut$};
\node[circle,scale=0.8,inner sep=0pt,fill=white] (A4) at (1.5,1.5) {$\strut$};
\node[circle,scale=0.8,inner sep=0pt,fill=white] (A6) at (2.5,2.5) {$\strut$};
\node[circle,scale=0.8,inner sep=0pt,fill=white] (A8) at (0.5,2.5) {$\strut$};

\begin{scope}[opacity=\leop]
\node[draw,circle,scale=0.8,inner sep=0pt,fill=white] (A2) at (0.5,0.5) {$t_2$};
\node[draw,circle,scale=0.8,inner sep=0pt,fill=white] (A1) at (1.5,0.5) {$t_1$};
\node[draw,circle,scale=0.8,inner sep=0pt,fill=white] (A3) at (2.5,1.5) {$t_3$};
\node[draw,circle,scale=0.8,inner sep=0pt,fill=white] (A4) at (1.5,1.5) {$t_4$};
\node[draw,circle,scale=0.8,inner sep=0pt,fill=white] (A6) at (2.5,2.5) {$t_6$};
\node[draw,circle,scale=0.8,inner sep=0pt,fill=white] (A8) at (0.5,2.5) {$t_8$};
\end{scope}

\node[scale=\Fscl,blue,anchor=north west,inner sep=\Fsep] (F1) at (A1.south east) {$F_1$};
\node[scale=\Fscl,blue,anchor=north west,inner sep=\Fsep] (F2) at (A2.south east) {$F_2$};
\node[scale=\Fscl,blue,anchor=north west,inner sep=\Fsep] (F3) at (A3.south east) {$F_3$};
\node[scale=\Fscl,blue,anchor=north west,inner sep=\Fsep] (F4) at (A4.south east) {$F_4$};
\node[scale=\Fscl,blue,anchor=north west,inner sep=\Fsep] (F6) at (A6.south east) {$F_6$};
\node[scale=\Fscl,blue,anchor=north west,inner sep=\Fsep] (F8) at (A8.south east) {$F_8$};
\def\qlw{1pt}
\def\qop{1}
\draw[->,blue,opacity=\qop,line width=\qlw] (F2)--($(F2)+(0,\arrlen)$);
\draw[->,blue,opacity=\qop,line width=\qlw] (F1)--($(F1)+(0,\arrlen)$);
\draw[->,blue,opacity=\qop,line width=\qlw] (F1)--($(F1)+(-\arrlen,0)$);
\draw[->,blue,opacity=\qop,line width=\qlw] ($(F1)+(-\arrlen,\arrlen)$)--(F1);
\draw[->,blue,opacity=\qop,line width=\qlw] (F4)--($(F4)+(-\arrlen,\arrlen)$);
\draw[->,blue,opacity=\qop,line width=\qlw] (F3)--($(F3)+(0,\arrlen)$);
\draw[->,blue,opacity=\qop,line width=\qlw] (F3)--($(F3)+(-\arrlen,0)$);
\draw[->,blue,opacity=\qop,line width=\qlw] ($(F3)+(-\arrlen,\arrlen)$)--(F3);
\draw[->,blue,opacity=\qop,line width=\qlw] (F6)--($(F6)+(-\arrlen,0)$);
\end{tikzpicture}
};

\node(X) at (11,0){
\begin{tikzpicture}[scale=\scl]
\def\nodescl{0.7}

\node[circle,scale=0.8,inner sep=0pt,fill=white] (A2) at (0.5,0.5) {$\strut$};
\node[circle,scale=0.8,inner sep=0pt,fill=white] (A2) at (0.5,0.5) {$\strut$};
\node[circle,scale=0.8,inner sep=0pt,fill=white] (A1) at (1.5,0.5) {$\strut$};
\node[circle,scale=0.8,inner sep=0pt,fill=white] (A3) at (2.5,1.5) {$\strut$};
\node[circle,scale=0.8,inner sep=0pt,fill=white] (A4) at (1.5,1.5) {$\strut$};
\node[circle,scale=0.8,inner sep=0pt,fill=white] (A6) at (2.5,2.5) {$\strut$};
\node[circle,scale=0.8,inner sep=0pt,fill=white] (A8) at (0.5,2.5) {$\strut$};

\node[scale=\Fscl,blue,anchor=north west,inner sep=\Fsep] (F1) at (A1.south east) {$F_1$};
\node[scale=\Fscl,blue,anchor=north west,inner sep=\Fsep] (F2) at (A2.south east) {$F_2$};
\node[scale=\Fscl,blue,anchor=north west,inner sep=\Fsep] (F3) at (A3.south east) {$F_3$};
\node[scale=\Fscl,blue,anchor=north west,inner sep=\Fsep] (F4) at (A4.south east) {$F_4$};
\node[scale=\Fscl,blue,anchor=north west,inner sep=\Fsep] (F6) at (A6.south east) {$F_6$};
\node[scale=\Fscl,blue,anchor=north west,inner sep=\Fsep] (F8) at (A8.south east) {$F_8$};
\def\qlw{1pt}
\def\qop{1}
\draw[->,blue,opacity=\qop,line width=\qlw] (F2)--(F8);
\draw[->,blue,opacity=\qop,line width=\qlw,dashed] (F1)--(F4);
\draw[->,blue,opacity=\qop,line width=\qlw,dashed] (F1)--(F2);
\draw[->,blue,opacity=\qop,line width=\qlw] (F8)--(F1);
\draw[->,blue,opacity=\qop,line width=\qlw] (F4)--(F8);
\draw[->,blue,opacity=\qop,line width=\qlw,dashed] (F3)--(F4);
\draw[->,blue,opacity=\qop,line width=\qlw,dashed] (F3)--(F6);
\draw[->,blue,opacity=\qop,line width=\qlw] (F8)--(F3);
\draw[->,blue,opacity=\qop,line width=\qlw] (F6)--(F8);
\end{tikzpicture}
};

\end{tikzpicture}
}
\caption{\label{Le quiver} An implementation of the rules from Figure \ref{local moves}.}
\end{figure}
\end{con}

\begin{exmp}
    Let $u = e$ and $w = (n - k + 1, ..., n, 1, ..., n - k)$, the longest $k$-Grassmannian word in $S_n$. Then the Le-diagram $D_{e,w}$ corresponding to this positroid pair is a $k \times (n-k)$ rectangular Young diagram with a dot in every box. Further, the construction applied to $D_{e,w}$ produces exactly Rietsch and Williams' rectangular plabic graph quiver for $\Pi_{k,n}$, save with all of the edges \textit{reversed}. Note that reversing the orientation of every edge in a quiver $Q$ does not change the corresponding cluster variety structure on $\mathcal{A}(Q)$.
\end{exmp}

\section{A-Side and Mirror Symmetric Preliminaries}

\subsection{Functoriality of Wrapped Fukaya Categories}

We freely utilize the sectorial descent methods of \cite{GPSI} to compute Fukaya categories of Weinstein manifolds, and we will both prove and use mirror symmetry statements for Weinstein sectors which are not Weinstein manifolds. Apart from the notion of $\Delta$-collaring and Lemma \ref{handle ends}, everything in this subsection is due to the seminal work of Ganatra-Pardon-Shende (\cite{GPS1}, \cite{GPSI}, \cite{GPSIII}) and to Asplund (see \cite{Asplund}), and we recall it only in the form we need. In particular, we do not recall the complete theory of Weinstein sectors here, and instead content ourselves to discuss the corresponding theory of \textit{stopped Weinstein manifolds}.\\

A stopped Weinstein manifold is a pair $(W,\mathfrak f)$ with $W$ Weinstein and $\mathfrak f\subset\partial^\infty W$ a closed, mostly Legendrian subset called a \textit{stop}. The \textit{skeleton} or \textit{core} of a stopped Weinstein manifold is the closed conic subset $core(W)\cup\mathbb R_{>0}\cdot\mathfrak f$, the core of $W$ together with the cone over $\mathfrak f$ under the Liouville flow. It is mostly Lagrangian, since $\mathfrak f$ is mostly Legendrian. A \textit{ribbon} for $\mathfrak f$ is a Weinstein hypersurface $F\subset\partial^\infty W$ with $core(F)=\mathfrak f$. Not every stop admits a ribbon, but the existence of one equips $\mathfrak f$ with a \textit{standard neighborhood}, the cone over a contact neighborhood of $F$, of the form $F \times\mathbb C_{\mathrm{Re}\le0}$ with $\mathfrak f$ at $\mathrm{Re}= -\infty$ (up to Weinstein deformation). Deleting such a neighborhood from $W$ gives a \textit{Weinstein sector} $W^\circ:= W \setminus \operatorname{Nbd}F$, whose symplectic boundary, the reduction of $\partial W^\circ$, is $F$. Conversely, $W$ is the convexification of $W^\circ$, and this correspondence is bijective up to contractible choice by Lemma 2.32 of \cite{GPS1}.\\

Associated to a stopped Weinstein manifold (or Weinstein sector, or pair of a Weinstein manifold and a Weinstein hypersurface $F$) is a \textit{partially wrapped Fukaya category}, which in the stopped case is denoted $\Fuk(W,\mathfrak f)$. Since there are natural functors $\Fuk(W^\circ) \to \Fuk(W,F) \to \Fuk(W,\mathfrak f)$ which are all equivalences (Corollary 3.9 of \cite{GPSI}) we elide a definition of the first two terms. Once a ribbon $F$ for $\mathfrak{f}$ is fixed we will pass freely between $W^\circ$, $(W,F)$ and $(W,\mathfrak f)$, however, and will tacitly employ the identification between the three corresponding categories.\\

To compute Fukaya categories we will utilize the microlocal and sectorial theory of \cite{GPSIII} and \cite{GPSI}. As such we will rarely actually consider Lagrangians or branes. The key theorem which allows us to do this is the following.

\begin{thm}[\cite{GPSIII} Theorem 1.1, in the form of \cite{gammage2022mirror} Theorem 2.5.5] Let $M$ be a real analytic manifold and $\Lambda\subset S^*M$ a closed subanalytic isotropic subset. Then $$\Fuk(T^*M,-\Lambda)\simeq\mathrm{Sh}_\Lambda(M)^c.$$ This equivalence carries the cotangent fiber at a point not under $\Lambda$ to a corepresentative of the stalk functor there. It carries the linking disk at a smooth point of $\Lambda$ to a corepresentative of the microstalk there. \end{thm}

There is a plethora of functors relating Fukaya categories of Weinstein sectors. For example, given an inclusion of Weinstein sectors $W' \hookrightarrow W$, \cite{GPSI} establishes that one receives a functor $\Fuk(W') \rightarrow \Fuk(W)$ by pushing forward Lagrangians. It turns out that Fukaya categories satisfy Mayer-Vietoris for such embeddings.

\begin{thm}[\cite{GPSI} Theorem 1.35]\label{sectorial}
   Let $\{U_1,\dots,U_k\}$ be a Weinstein sectorial cover of a Liouville sector $W$, i.e. the symplectic reductions of all strata of $\partial W,\partial U_1,\dots,\partial U_k$ are Weinstein up to deformation (see Definition 1.32 of \cite{GPSI}). If $W$ is stopped, assume the stop is mostly Legendrian and disjoint from every $\partial U_i$ and from a neighborhood of $\partial W$. Then the induced functor $$\colim_{\emptyset \neq I \subset \{1, ..., k\}} \Fuk(\bigcap_{i \in I} U_i) \rightarrow \Fuk(W)$$ \noindent is an equivalence of stable $\infty$-categories.
\end{thm}

A key use of this is that a ribbon $F$ for a piece $\mathfrak{f}$ of a stop $\mathfrak{g}$ determines a sectorial inclusion $F \times T^*[0,1]\hookrightarrow(W,\mathfrak g)$, and hence a functor $$\iota: \Fuk(F)\simeq \Fuk(F\times T^*[0,1])\rightarrow \Fuk(W, \mathfrak{g})$$

\noindent where the first equivalence is the K\"{u}nneth functor of \cite{GPSI}, since $\Fuk(T^*[0,1]) \simeq \mathrm{Perf}(k)$. The subsector $F \times T^*[0,1]$ of $(W, \mathfrak{g})$ is called a \textit{collar} for $\mathfrak{f}$, and the functor above is called a \textit{cup} functor. Conversely, any subsector $N\cong F\times T^*[0,1]$ of $(W,\mathfrak{g})$ such that $\mathfrak{g}\cap(\partial^\infty N)^\circ$ is the stop coming from $F$, and such that each $\partial^\infty(F\times T^*_t[0,1])$ links $\mathfrak{f}$, is a collar for $\mathfrak{f}$.\\

In any event, cupping sends cocores of $F$ to linking disks of $\mathfrak f$. When $F$ is Weinstein and stably polarized (the latter being automatic if $W$ is itself stably polarized) then $\mu sh_{-\mathfrak f}(\mathfrak -f)^c \simeq \Fuk(F)$ by the main theorem of \cite{GPSIII}, in which case $\iota$ is left adjoint to microlocalization along $\mathfrak f$, $\Fuk(W,\mathfrak g)\to \mu sh_{\mathfrak -f}(\mathfrak -f)^c$.\\ 

One can also compare Weinstein sectors with the same underlying manifold by shrinking the stop. That is, if $\mathfrak{f} \subset \mathfrak{g}$ are a pair of closed Legendrians in $\partial^\infty W$ such that $\mathfrak{g} \setminus \mathfrak{f} \subset \partial^\infty W \setminus \mathfrak{f}$ is Legendrian  then there is a \textit{stop-removal} functor $$\mathfrak{sr}: \Fuk(W, \mathfrak{g}) \rightarrow \Fuk(W, \mathfrak{f})$$ \noindent which is a Verdier localization with kernel generated by the linking disks to the critical parts of $\mathfrak{g} \setminus \mathfrak{f}$, see \cite{GPSI} Theorem 1.20. In particular, if $\mathfrak f$ admits a ribbon $F$, then the kernel of $\mathfrak{sr}\colon \Fuk(W,\mathfrak f)\to \Fuk(W)$ is generated by the image of the cup functor.\\

The only non-standard notion we introduce in this subsection is that of $\Delta$-collaring, which we think of as an A-side analogue of a simple normal crossings boundary with strata indexed by $\Delta$. This is essentially a repackaging of Asplund’s simplicial decompositions (see \cite{Asplund}) designed so that handle attachments can be iterated (see Lemma \ref{handle ends}).\\

To give the definition, recall from Section 12 of \cite{GPSI} that a \textit{Liouville sector-with-corners} $W^\circ$ is a Liouville manifold-with-corners whose boundary, viewed as an immersed hypersurface, is sectorial in the sense of Definitions 12.2 and 12.14 of \cite{GPSI}. The \textit{boundary faces} $\partial_v$ of such a $W^\circ$ are defined to be the local branches of this hypersurface. All the facts we use about sectors-with-corners are from Section 12 of \cite{GPSI}.\\

With these assumptions each corner stratum $\partial_K:=\bigcap_{v\in K}\partial_v$ is coisotropic, and the corresponding symplectic reduction $F^\circ_K$ (\textit{a symplectic boundary stratum} or \textit{symplectic face}, if $|K|=1$) is again a Liouville sector-with-corners of dimension $\dim W-2|K|$ by Lemma 12.4 of \cite{GPSI}. Near $\partial_K$ there are coordinates $F_K^\circ\times T^*\mathbb{R}^K_{\le 0}$ given by defining functions $I_v, t_v$, with respect to which the faces are the preimages of the coordinate hyperplanes (Lemma 12.8 of \cite{GPSI}). Here each $I_v$ is chosen with respect to the face $\partial_v$ as in Definition 12.2 of \cite{GPSI}, and $t_v$ is then determined by $I_v$ by Lemma 12.7 of the same. Thus for $K\subset K'$ the coordinates near $\partial_{K'}$ are those near $\partial_K$ together with $(I_v,t_v)_{v\in K'\setminus K}$. After deforming the Liouville form, we may moreover assume that these coordinates split it as $\lambda=\lambda_{F_K^\circ}+\lambda_{T^*\mathbb{R}^K}$, for all $K$ simultaneously (see Lemma 12.12 of \cite{GPSI}). In the sequel we call these coordinates \textit{standard coordinates}.\\

Now assume given a finite set $V$ and a simplicial complex $\Delta$ on $V$, and write $\Delta^+:=\Delta\setminus\{\varnothing\}$. In our use, $\Delta$ will usually be the complex of cones of a simplicial fan $\Sigma$.

\begin{defn}\label{collared}
A Weinstein sector-with-corners $W^\circ$ is \textit{$\Delta$-collared} if its boundary faces are embedded and indexed by $V$, and its corner strata $\partial_K$ are nonempty exactly for $K\in\Delta^+$.
\end{defn}

By using standard coordinates one sees that any $F_K^\circ$ in a $\Delta$-collared $W^\circ$ is itself $\mathrm{link}_\Delta(K)$-collared. Further, smoothing the corners of $W^\circ$ gives a Weinstein sector whose symplectic boundary $F_D$ is Weinstein and is covered by slight enlargements of the symplectic faces $F_v^\circ$ (by Example 12.24 and Lemma 12.26 of \cite{GPSI}). Let $W$ be its convexification. We call $(W,\mathfrak g)$, with $\mathfrak g:=core(F_D)$, \textit{$\Delta$-collared}, and we refer to $W^\circ$ as its \textit{collared sector}. Applying this recursively we see that each convexified stratum $F_K$ is itself a $link_\Delta(K)$-collared stopped manifold $(F_K,\mathfrak g_K)$.\\

Let $(W,\mathfrak g)$ be $\Delta$-collared with collared sector $W^\circ$. For $v\in V$, the \textit{piece} $l_v\subset\mathfrak g$ is the relative core of the stopped manifold $(F_v,\mathfrak g_v)$ (i.e., $F_v$ is a ribbon for $l_v$). For $K\in\Delta^+$, the \textit{end} $N_K \subset W^\circ$ is the part of the coordinate neighbourhood $F^\circ_K\times T^*\mathbb R^K_{\le0}$ of $\partial_K$ lying over $[-1,0]^K$, so that

$$N_K\cong F^\circ_K\times(T^*[0,1])^K \quad \text{and} \quad (T^*[0,1])^K:=\textstyle\prod_{k\in K}T^*[0,1].$$

Since $\Fuk(F^\circ_K)\simeq \Fuk(F_K)$, we often abuse notation and write $N_K\cong F_K\times(T^*[0,1])^K$. With this notation fixed we may give a lemma about the geometry of $\Delta$-collared sectors.

\begin{lem}\label{delta collard features}
    Assume given a $\Delta$-collared $(W,\mathfrak g)$.

    \begin{enumerate}
        \item $\mathfrak g=\bigcup_{v\in V}l_v$, and each end $N_K$ is a subsector of $(W,\mathfrak g)$. In particular $N_{\{v\}}$ is a collar of $F_v$, whose cotangent fibers are linking arcs of $l_v$.
        \item $N_K\cap N_{K'}=N_{K\cup K'}$, which is empty if $K\cup K'\notin\Delta$. For $K\subset K'$, the inclusion $N_{K'}\subset N_K$ is the end of $F_K$ at $K'\setminus K$ times the identity on $(T^*[0,1])^K$.
        \item A product of collared stopped manifolds is collared by the join of their complexes.
    \end{enumerate}

    Conversely, assume given subsectors $N_K\cong F_K\times(T^*[0,1])^K$ of a stopped Weinstein manifold $(W,\mathfrak g)$ that are nonempty exactly for $K\in\Delta^+$, satisfy (2), and whose ribbons $F_v$ glue to a ribbon of $\mathfrak g$. Then they are the ends of a $\Delta$-collared structure on $(W,\mathfrak g)$.
\end{lem}

\begin{proof}
    The ends are subsectors satisfying $(2)$, since standard coordinates are compatible for nested strata and split the Liouville form. In particular $N_{\{v\}}$ satisfies the criterion for collars given above. Since $F_D$ is covered by the $F_v$, $\mathfrak g=\bigcup_vl_v$ by Lemma 2.17 of \cite{Asplund}. Item $(3)$ is Example 12.16 of \cite{GPSI}. For the converse: by $(2)$, slight enlargements of the $F_v$ form a good sectorial cover of the ribbon with nerve $\Delta$ by Definition 3.7 of \cite{Asplund}, and the same work also shows that such a cover is equivalent to a $\Delta$-collared sector (see Lemma 3.9 and Theorem 3.11 of \cite{Asplund}).
\end{proof}

Now we may define the kind of Weinstein handle attachments we utilize.

\begin{defn}\label{collared handles}
Let $(W,\mathfrak g)$ be $\Delta$-collared, fix $u\in V$, and suppose $F^\circ_u\cong F'^\circ\times T^*S^1$ \textit{compatibly with collars}, i.e. $F'^\circ$ is $link_\Delta(u)$-collared and $F^\circ_{K\cup u}\cong F'^\circ_K\times T^*S^1$ compatibly with standard coordinates.\\

Set $h:=F'\times T^*D^2$. Its sector $h^\circ:=F'^\circ\times T^*D^2$ has boundary faces $\partial'_v\times T^*D^2$ and $F'^\circ\times T^*D^2|_{\partial D^2}$, and the latter has symplectic face $F'^\circ\times T^*S^1\cong F_u^\circ$ by $(3)$ of Lemma \ref{delta collard features}. The \textit{collared handle attachment} $(W\cup_uh)^\circ$ glues $W^\circ$ and $h^\circ$ along their faces with symplectic face $F^\circ_u$ (as in Construction 12.18 of \cite{GPSI}), and $W\cup_uh$ is its convexification.
\end{defn}

The assumptions ensure that $W\cup_uh$ looks like $F'$ times the $A_1$ handle attachment of \cite{GL} near $N_u$.\\

The following lemma is the main new input of this subsection; its single-face, single-handle case is Construction 12.18 of \cite{GPSI}. To state it, write $\Delta_{\setminus u}$ for the subcomplex of faces not containing $u$. For $K \in link_\Delta(u)$, the stratum of $F^\circ_K$ at $u$ is $F^\circ_{K\cup u}\cong F'^\circ_K\times T^*S^1$, compatibly with collars, so $F_K\cup_u(F'_K\times T^*D^2)$ is defined.
 
\begin{lem}\label{handle ends}
Assume that $W^\circ$ is $\Delta$-collared and that $(W\cup_uh)^\circ$ is the result of a collared handle attachment. Then the latter is $\Delta_{\backslash u}$-collared. Its faces are the faces $\partial_v$, $v\ne u$, of $W^\circ$, glued for $v\in\mathrm{link}_\Delta(u)$ to the faces $\partial'_v\times T^*D^2$ of $h^\circ$, and its stratum at $K\in(\Delta_{\backslash u})^+$ is $F_K\cup_u(F'_K\times T^*D^2)$ if $K\in\mathrm{link}_\Delta(u)$ and is $F_K$ otherwise.\\

In particular this implies that $W\cup_uh$ is $\Delta_{\backslash u}$-collared, and that its stop is obtained from $\mathfrak g$ by removing $l_u$ and extending each $l_v$, $v\in\mathrm{link}_\Delta(u)$, by $l'_v\times D^2$, i.e. by attaching a parametrized disk.
\end{lem}
 
\begin{proof}
This is the many-face version of Construction 12.18 of \cite{GPSI}. Near the glued faces, $W^\circ$ and $h^\circ$ are $F^\circ_u\times T^*\mathbb R_{\ge0}$ and $F^\circ_u\times T^*\mathbb R_{\le0}$ in standard coordinates. For $v \in link_\Delta(u)$, compatibility with collars means the two sides also share the normal form $F^\circ_{\{u,v\}}\times T^*\mathbb R^2$ near $\partial_{\{u,v\}}$. So the defining functions $I_v,t_v$ glue, and $\partial_v\cup(\partial'_v\times T^*D^2)$ is a sectorial face. Faces away from $\partial_u$ are unaffected.\\

For the strata, take $K\in(\Delta_{\setminus u})^+$. If $K\notin link_\Delta(u)$, then $\partial_K$ misses $h^\circ$. Otherwise, its reduction is $F^\circ_K$ glued to $F'^\circ_K\times T^*D^2$ along $F^\circ_{K\cup u}$. By induction on dimension, this is the collared handle attachment $F_K\cup_u(F'_K\times T^*D^2)$. That attachment is a Weinstein connected sum along a Weinstein face, see Definition 2.7 of \cite{Asplund}, so it is Weinstein.
\end{proof}

As a final note on general symplectic geometry we record the following

\begin{thm}\label{thm:handle-pushout}

Assume given a collared handle attachment $W\cup_uh$ as in Definition \ref{collared handles}. Let $\iota\colon \Fuk(F_u)\to \Fuk(W,\mathfrak g)$ and $r\colon \Fuk(F_u)\to \Fuk(h)$ be the cup functors of the face along which $W^\circ$ and $h^\circ$ are glued. Let $\jmath\colon W^\circ\hookrightarrow(W\cup_uh)^\circ$ and $\kappa\colon h^\circ\hookrightarrow(W\cup_uh)^\circ$ be the inclusions. Then

\[\begin{tikzcd}
	\Fuk(F_u) \arrow[r, "\iota"] \arrow[d, "r"] & \Fuk(W, \mathfrak{g}) \arrow[d, "\jmath_*"] \\
	\Fuk(h) \arrow[r, "\kappa_*"]& \Fuk(W \cup_u h)
\end{tikzcd}\]

is a pushout of stable $\infty$-categories.
\end{thm}

\begin{proof}
    Slight enlargements of $W^\circ$ and $h^\circ$ form a sectorial cover of $(W\cup_uh)^\circ$ whose intersection is the collar $F_u\times T^*[0,1]$ of the glued face (as in Lemma \ref{handle ends}). Its strata have reductions $F_K$ and $F'_K\times T^*D^2$, which are Weinstein by Definitions \ref{collared} and \ref{collared handles}, so the cover is Weinstein. The theorem then follows from Theorem 3.2 and the identification $\Fuk(F_u\times T^*[0,1])\simeq \Fuk(F_u)$ from K\"{u}nneth.
\end{proof}

Here is our main example, and an excuse to fix notation. We fix character and cocharacter lattices $M = \langle \theta_i\rangle$ and $N = \langle \theta_i^\vee\rangle$ and an $(m+n)$-torus $T \cong M \otimes \mathbb{R}/2\pi\mathbb{Z}$ throughout. By abuse of notation we will also denote the angular coordinates on $T$ by $\theta_i$. In particular, with these conventions we may write $T^* T \cong T \times N_\mathbb{R}.$ Also, in the case where $\Sigma(1)$ consists only of coordinate vectors $\theta_i^\vee \in N_\mathbb{R}$, given a cone $\langle \theta_i^\vee \:|\: i \in I\rangle = \sigma_I \in \Sigma$, we write $$T_I = \{\theta \in T\:|\: \langle \theta, \theta_i^\vee\rangle = 0 \:\forall i \in I\} = M_I \otimes \mathbb{R}/2\pi\mathbb{Z}.$$

Now let $X_\Sigma$ be the toric variety  corresponding to a fan $\Sigma \subset N \otimes \mathbb{R} := N_\mathbb{R}$. Associated to $\Sigma$, \cite{FLTZ} defined a conic (singular) \textit{FLTZ Lagrangian} in $T^* T$ by $$\mathbb{L}_\Sigma = \bigcup_{\sigma \in \Sigma} \sigma^\perp \times -\sigma \subset T \times N_\mathbb{R} \cong T^* T.$$

In the case that $\Sigma(1)$ consists only of coordinate rays, we note that this can be rewritten as $$\mathbb{L}_\Sigma = \bigcup_{\sigma_I \in \Sigma} T_I \times (-\sigma_I),$$

\noindent and in this case we will use the decomposition $T^* T \cong T^* T_i \times T^* S^1_{\theta_i}$ induced by the splitting $M \cong M_i \times \mathbb{Z}\theta_i$ often.\\

From the FLTZ Lagrangian we receive a sector $X_\Sigma^\vee := (T^* T, -\mathbb{L}_\Sigma \cap \partial^\infty T^* T)$, where we define $\Lambda_\Sigma := -\mathbb{L}_\Sigma \cap (\partial^\infty T^* T)$ in the sequel. Then, the main theorems of \cite{FLTZ} and \cite{Kuwagaki}, which are \textit{the} key inputs to the remainder of this study, say that $$(Sh_{-\Lambda_\Sigma} ( T))^c \simeq \Fuk(X_\Sigma^\vee) \simeq \Coh(X_\Sigma) .$$

For our applications the reader need only know a single example, though we encourage the reader to consult \cite{GammageThesis} for a beautiful and more sustained discussion of the FLTZ formalism.

\begin{exmp}\label{affine space toric example}
    If $\Sigma$ is a single non-zero ray in $\mathbb{R}$ then $X_\Sigma = \mathbb{A}^1$. In this case, the prescription above unravels to say that $X^\vee_\Sigma=(T^*S^1,\partial^\infty(\{0\}\times\mathbb R_{\ge0}\theta_1^\vee))$, which we abbreviate to $(T^*S^1,\infty)$. One can visualize the corresponding skeleton as a circle with ray glued onto $0$.
\end{exmp}

We note that the $X^\vee_\Sigma$ just considered is collared over the single vertex $\rho_1$: removing a standard neighborhood of the stop point leaves a sector with one face, whose reduction is a point.

\subsection{Generalities on Toric HMS}

Together, \cite{FLTZ}, \cite{Kuwagaki}, and \cite{gammage2022mirror} establish homological mirror symmetry for (toric boundary divisors in) toric varieties. In the process, they also establish several naturality/functoriality statements in toric HMS. Distilled to its essentials, their combined works show the following package of theorems.

\begin{thm}[Toric HMS and its Naturality]\label{toric HMS natural}

Assume given a smooth toric variety $X_\Sigma$ with mirror $X_\Sigma^\vee = (T^*T,\Lambda_\Sigma)$ with $\Lambda_\Sigma:=\partial^\infty(-\mathbb{L}_\Sigma)$ the FLTZ stop, and let $\Phi_\Sigma \colon \Fuk(X^\vee_\Sigma) \simeq \Coh(X_\Sigma)$ be Kuwagaki’s equivalence \cite{Kuwagaki}.

\begin{enumerate}
    \item A subfan inclusion $\Sigma' \subset \Sigma$ induces a commutative diagram \begin{center}
    \begin{tikzcd}
\Fuk(X_\Sigma^\vee) \arrow[r, "\mathfrak{sr}"'] \arrow[d, "\Phi_\Sigma"] & \Fuk((X_{\Sigma'})^\vee) \arrow[d, "\Phi_{\Sigma'}"] \\
\Coh(X_\Sigma) \arrow[r, "j^*"']                                          & \Coh(X_{\Sigma'}),                                 
\end{tikzcd}
\end{center}
intertwining open restriction and stop removal (noting that $\Lambda_{\Sigma'} \subset \Lambda_{\Sigma}$). These diagrams are natural in subfan inclusions.

\item Write $V(\sigma_I) := X_{\Sigma_{-I}}\hookrightarrow X_\Sigma$ for the orbit closure corresponding to a cone $\sigma_I \in \Sigma$, embedded via $k_I$. Then for any such $\sigma_I$ there is an induced sectorial pushforward $\iota_I: \Fuk(X^\vee_{\Sigma_{-I}})\to \Fuk(X^\vee_\Sigma)$, with right adjoint $\iota_I^R$ participating in a commutative diagram 
\begin{center}
\begin{tikzcd}
\Fuk(X_\Sigma^\vee) \arrow[r, "\iota_I^R"'] \arrow[d, "\Phi_{\Sigma}"] & \Fuk(X_{\Sigma_{-I}}^\vee) \arrow[d, "\Phi_{\Sigma_{-I}}"] \\
\Coh(X_\Sigma) \arrow[r, "k_I^*"']                                     & {\Coh(X_{\Sigma_{-I}}),}                                  
\end{tikzcd}
\end{center}

\noindent intertwining closed restriction and the right adjoint. Further, in this case the functor $\iota_I^R$ is identified with the microlocalization functor between the corresponding microsheaf categories. These diagrams are natural in inclusions of cones.
\item For $x \in T \setminus \pi(\Lambda_{\Sigma})$ (where $\pi$ is the projection from $T^* T$ to the base $T$), the cotangent fiber $T^*_x T$ at $x$ is mapped via $\Phi_\Sigma$ to the structure sheaf of $X_\Sigma$.

\item The toric boundary $\partial X_\Sigma$ of $X_\Sigma$ is mirror to the Legendrian at infinity of $X_\Sigma^\vee$: $$\mu sh_{-\Lambda_\Sigma} (-\Lambda_\Sigma)^c \simeq \Coh(\partial X_\Sigma).$$
\end{enumerate}

Finally, the above is compatible with lattice isomorphisms and with products: $\Phi_{\Sigma\times\Sigma'}\simeq\Phi_\Sigma\boxtimes\Phi_{\Sigma'}$ under Künneth.
\end{thm}

\begin{proof}[Proof Sketch]
    $(1)$ is standard, and both $(3)$ and $(4)$ are Section 7 of \cite{gammage2022mirror}. $(2)$ and the product claim follow from unraveling Section 7.2 of \cite{gammage2022mirror}, though it is not stated as such there.
\end{proof}

An important input to the proofs in \cite{gammage2022mirror} is a recursive combinatorial understanding of the structure of the FLTZ stops $\mathbb{L}_\Sigma$. In words, it turns out that if one fixes a coordinate $\theta_i^\vee \in N$ with corresponding ray $\rho_i \in \Sigma$, then the intersection of $-\mathbb{L}_\Sigma \subset (M \otimes \mathbb{R}/2\pi\mathbb{Z}) \times N_\mathbb{R}$ with the region $\theta_i^\vee > 0$ looks like $-\mathbb L_{\Sigma_{-i}}$ times a ray. This means that the part of the stop $\Lambda_\Sigma$ lying over $\{\theta_i^\vee>0\}$ is a copy of $-\mathbb L_{\Sigma_{-i}}$, the core of the Weinstein sector $X^\vee_{\Sigma_{-i}}$ which is the mirror to the torus orbit closure $V(x_i)$. \\

To state this formally, write $\Sigma_{-I}$ for the deleted fan/quotiented star with respect to a cone $\sigma_I \in \Sigma$ (as in Lemma \ref{fan operations}), and continue to fix a basis of $\theta_i$'s and $\theta_i^\vee$'s for $M$ and $N$ with $\rho_i \in \Sigma(1)$. We split $M=M_I\oplus\bigoplus_{i\in I}\mathbb Z\theta_i$, and similarly for $N$, giving $T^*T\cong T^*T_I\times T^*(S^1)^I\cong(T_I\times(N_I)_{\mathbb R})\times((\mathbb R/2\pi\mathbb Z)^I\times\mathbb R^I)$. With the above fixed, we learn the following

\begin{prop}[\cite{gammage2022mirror} Lemma 4.3.1 and Corollary 4.3.2]\label{toric ends}

$$-\mathbb{L}_\Sigma \cap \{\theta_i^\vee > 0\} \cong -\mathbb{L}_{\Sigma_{-i}} \times (\{0\}\times\mathbb R_{>0}),$$ where we have utilized the splitting above and written the product of $-\mathbb{L}_{\Sigma_{-i}}$ in the $T^* S^1_{\theta_i}$ factor and $\{0\}\times\mathbb R_{>0}$ in the $T^* \overline{T_i} \cong T^* S^1$ factor. More generally, given $\sigma_I \in \Sigma$, we have $$-\mathbb{L}_\Sigma \cap \{\theta_i^\vee > 0\:|\:i \in I\} \cong -\mathbb{L}_{\Sigma_{-I}} \times \prod_{i \in I} (\{0\}\times\mathbb R_{>0})$$ \noindent with coordinates internal to the decomposition above.
\end{prop}

Further, one sees that if $\rho_i$ and $\rho_j$ don't span a cone of $\Sigma$, then $(-\mathbb{L}_\Sigma)\cap\{\theta_i^\vee>0\}\subset\{\theta_j^\vee=0\}$.\\

Let $\Sigma$ be a fan in $N_{\mathbb R}$ whose rays are parametrized by a set $J$ and are amongst the basis vectors $\theta_i^\vee$ of $N$. Let $\Sigma_{\max}$ be the fan of all cones $\sigma_I$ with $I\subset J$. Then $X_{\Sigma^{max}} \cong \mathbb C^{J}\times(\mathbb C^*)^{rk(N) - |J|}$ and $X^\vee_{\Sigma_{\max}}=\prod_{j\in J}(T^*S^1_j, \infty)\times T^*T^{rk(N) - |J|}$. For $\sigma_I\in\Sigma^{max}$ we define a subsector $N_I$ corresponding to $\sigma_I$ as the product of the strips $T^*[0,1]\subset(T^*S^1_i,\infty)$ about the stop points $(\{0\}\times\mathbb R_{>0}) \cap \partial^\infty$ over $i \in I$, all further producted with the remaining factors of $X^\vee_{\Sigma^{max}}$. Therefore $N_I=X^\vee_{\Sigma_{-I}^{max}}\times(T^*[0,1])^I$.\\

Given $\Sigma \subset \Sigma^{max}$ and $\sigma_I \in \Sigma$, we claim that $N_I$ naturally restricts to a subsector of $X_\Sigma^\vee$. Indeed, by the last equation of the preceding proposition, the stop of $X^\vee_\Sigma$ meets $N_I$ exactly in the stop of $X^\vee_{\Sigma_{-I}}$ times the strips. Therefore we see that $N_I$ restricted to $X_\Sigma^\vee$, which by abuse of notation we also denote by $N_I$, satisfies $N_I\cong X^\vee_{\Sigma_{-I}}\times(T^*[0,1])^I$. From this one receives a functor $i_I \colon \Fuk(X_{\Sigma_{-I}}^\vee) \to \Fuk(X^\vee_\Sigma)$ via K\"{u}nneth and sectorial pushforward.\\

Finally, for $I=\{i\}$, the chart $\{\theta_i^\vee=1\}$ of $\partial^\infty T^*T$ is $T^*T_i\times S^1_z$, where $z$ is the angle along $S^1_{\theta_i}$, with contact form $dz+\lambda_{T^*T_i}$. By the preceding proposition again, the stop piece $l_i$ near $+\theta_i^\vee$ is $-\mathbb L_{\Sigma_{-i}}$ inside the Liouville hypersurface $\{z=0\}\cong T^*T_i$. This is the relative core of $X^\vee_{\Sigma_{-i}}$, since $-\mathbb L_{\Sigma_{-i}}$ is conic. Therefore $X^\vee_{\Sigma_{-i}}$ is a ribbon of $l_i$ in $X^\vee_{\Sigma}$.\\

The above discussion is already nearly sufficient to exhibit certain collaredness properties of $X^\vee_\Sigma$. Namely, writing $\Delta_\Sigma$ to denote the simplicial complex of cones of a simplicial fan $\Sigma$, we show now that toric mirrors are in fact $\Delta_\Sigma$-collared.

\begin{prop}\label{lem:toric-collared}
Let $\Sigma$ be a fan whose rays are among the basis vectors $\theta_v^\vee$ of $N$. Then $X^\vee_\Sigma$ is $\Delta_\Sigma$-collared with $F_K=X^\vee_{\Sigma_{-K}}$ and ends $N_I$ for all $\sigma_I \in \Sigma$.
\end{prop}
 
\begin{proof}
For the fan $\Sigma_{\max}$ this is immediate. Indeed, in this case $X^\vee_{\Sigma_{\max}}=\prod_k(T^*S^1_{\theta_k},\infty)\times T^*T'$ is a product of sectors, each collared over a single vertex (Example \ref{affine space toric example}), and so it is $\Delta_{\Sigma_{\max}}$-collared by Lemma \ref{delta collard features}(3), with the ends $N_I$ described above.\\

Now assume given $\Sigma\subset\Sigma_{\max}$. We have seen that for $\sigma_I\in\Sigma$ the $N_I$ restrict to subsectors $N_I\cong X^\vee_{\Sigma_{-I}}\times(T^*[0,1])^I$ of $X^\vee_\Sigma$, and by Proposition \ref{toric ends} these are nonempty exactly when $\sigma_I\in\Sigma$. Since they are all cut out by the splitting $T^*T\cong T^*T_I\times T^*(S^1)^I$, their product coordinates agree on overlaps, i.e. they satisfy Lemma \ref{delta collard features}(2). Further, the ribbons $X^\vee_{\Sigma_{-i}}$ of the pieces $l_i$ glue to a ribbon of $\Lambda_\Sigma$. Indeed, the ribbon for $\Sigma_{\max}$ is built by attaching the handles indexed by each $K$ along those indexed by the faces of $K$, by Lemma 2.17 of \cite{Asplund}. Since $\Delta_\Sigma$ is closed under passing to faces, the handles indexed by $\Delta_\Sigma$ form a Weinstein subdomain, whose core is $\Lambda_\Sigma$ by Proposition \ref{toric ends}. Therefore the converse in Lemma \ref{delta collard features} applies, and $X^\vee_\Sigma$ is $\Delta_\Sigma$-collared with $F_K=X^\vee_{\Sigma_{-K}}$.
\end{proof}
 
In particular we learn that, near each cone, $X^\vee_\Sigma$ looks like the local model of a toric stratum: i.e., that the end $N_K$ is the mirror of $V(\sigma_K)$ times one strip for each normal direction, and compatibly so.\\

Our main geometric goal in the remainder of this work is to establish analogues of the naturality statements given above in the context of cluster models. The crucial $B$-sided input for this goal is the following.

\begin{thm}[\cite{GL} Theorem 4.16]\label{affine blowup} Assume given a smooth variety $X$, a divisor $k\colon D\hookrightarrow X$, and a smooth subvariety $i_H\colon H\hookrightarrow D$ of codimension two in $X$. Let $\pi\colon\mathrm{Bl}^{\mathrm{aff}}_HX\to X$ be the affine blowup, with exceptional divisor $E:=\pi^{-1}(D)$, inclusion $\iota_E\colon E\hookrightarrow\mathrm{Bl}^{\mathrm{aff}}_HX$, and $\mathbb A^1$-bundle projection $q\colon E\to H$. Then the square

\[\begin{tikzcd}
	\Coh(D) \arrow[r, "k_*"] \arrow[d, "i^*_H"] & \Coh(X) \arrow[d, "\pi^*"] \\
	\Coh(H) \arrow[r, "\iota_{E*}q^* "]& \Coh(\text{Bl}_H^{\text{aff}}X)
\end{tikzcd}\]

\noindent is a pushout of stable $\infty$-categories.
\end{thm}

\begin{rem}
    Since $(k^*)^L\simeq k_*[-1]$ by Remark \ref{GV duality remark}, the natural transformation arising from applying the shift $[1]$ on $\Coh(D)$ and $\Coh(H)$, and the identity on $\Coh(X)$ induces an equivalence of the square in the previous theorem with the square

\[\begin{tikzcd}
	\Coh(D) \arrow[r, "(k^*)^L"] \arrow[d, "i^*_H"] & \Coh(X) \arrow[d, "\pi^*"] \\
	\Coh(H) \arrow[r, "{\iota_{E*}q^*[-1]}"] & \Coh(\text{Bl}_H^{\text{aff}}X)
\end{tikzcd}\]

which is therefore also a pushout. Here $\iota_{E*}q^*[-1]\simeq(\iota_E^*)^Lq^*$.
\end{rem}

Theorems \ref{thm:handle-pushout} and \ref{affine blowup} have the same shape: in each, the category of the new space is a pushout of the old one and a new piece, along a common face. Therefore one may hope to mirror each affine blowup of Construction 2.24 by a handle attachment. Gammage and Le carry this out for $\mathcal A(A_1)$, as follows. First, they show in Theorem 4.12 of \cite{GL} that the following two spans are equivalent: $$\begin{tikzcd}
\Fuk(T^* D^2) \arrow[d, "\simeq"] &  & \Fuk(T^*S^1) \arrow[ll, "r"'] \arrow[rr, "\iota"] \arrow[d, "\simeq"]               &  & \Fuk((\mathbb{C} \times \mathbb{C}^*)^\vee) \arrow[d, "\simeq"] \\
{\Coh(\{(0, -1)\})}               &  & \Coh(\{0\} \times \mathbb{C}^*) \arrow[rr, "(k^*)^L"] \arrow[ll, "i^*_H"'] &  & \Coh(\mathbb{C} \times \mathbb{C}^*).                          
\end{tikzcd}$$

They then notice that by Theorem \ref{affine blowup}, the pushout of the bottom line is $\Coh(Bl_H^{\text{aff}}(\mathbb{C} \times \mathbb{C}^*))$. But Construction \ref{GHK+} tells us this is $\Coh(\mathcal{A}(A_1)_{\Sigma_{A_1}, \Sigma_{A_1}(1)}) = \Coh(\mathcal{A}(A_1))$! Therefore we learn that $\mathcal{A}(A_1)$ is, indeed, mirror to a Weinstein manifold with skeleton as depicted in Figure \ref{A1 skeleton}.\\

Gammage and Le actually construct a similar such equivalence, with $A$-side presented via handle attachments, for any truncated cluster variety -- in our parlance, for any cluster model $\mathcal{A}(Q)_{\Sigma(1), \Sigma(1)}$. The game, therefore, is to show that something similar works for \textit{any} cluster model $\A$: i.e., that $\A^\vee$ can be constructed by iterated Weinstein handle attachments from the pre-blowup toric model for $\A$. To do this we begin with some preliminaries.
%%We expect that the varieties $X$ we consider are self-mirror, meaning that the Fukaya category of the Stein manifold underlying $X$ and the derived category of coherent sheaves on $X$ are equivalent. This leads us to seek comparisons between $\Fuk(\X)$ and $Coh(X)$ for a good supply of open subsets of $X$. The following two computations are the key inputs of this form; we refer the reader to \cite{Auroux} for their proof.

%%\begin{prop}\label{torus HMS}
    %%$\Fuk((\mathbb{C}^*)^m) \simeq \Coh((\mathbb{C}^*)^m)$ and 
    %%$\Fuk(\mathbb{C}^n) \simeq \Coh(\mathbb{C}^n).$
%%\end{prop}

%%This implies, for example, that every piece of the $V/D$ decomposition of a really full rank acyclic cluster variety is self-mirror.\\

\subsection{Mirrors to Cluster Models}

Assume given a cluster datum $(Q, \Sigma, A)$. This subsection begins with final combinatorial inputs which allow us to define the mirror to $\mathcal{A}(Q)_{\Sigma, A}$ and then describes the actual symplectic geometry of this construction.\\

The work in this subsection was developed in conversations with Peng Zhou, and we think of Construction \ref{Weinstein from Fan} as essentially due to \cite{SSWZ} (if not verbatim so in the case $\Sigma = \Sigma_Q$ and $A = \Sigma_Q(1)$). Nevertheless the reader will find substantial differences in appearance between this construction and the one in \cite{SSWZ}. Indeed, we have found certain modifications and generalizations of their formalism necessary to accommodate both the combinatorics and symplectic topology of the mirrors to general cluster models, and we are unsure of how to prove our main theorems without these changes to their construction.\\

Recall that we often write $i \in A$ to mean $\rho_i \in A$ (or similar such things). Inspired by Gammage and Le's work, we seek to construct the mirror $\mathcal{A}(Q)_{\Sigma, A}^\vee := \mathcal{X}(Q)_{\Sigma, A}$ through a series of Weinstein handle attachments. To start we identify the circles in $T := M\otimes\mathbb R/2\pi\mathbb Z$ which the handles will cap.

\begin{defn}
    Given the columns $c_i$ of $\tilde{E}$, consider the corresponding element $c_i \in M$ and exchange polynomial $P_i = x^{[-c_i]_+}(x^{c_i} + 1) := u_i(y_i + 1)$. Then, the \textit{cap cocharacter} with respect to $i \in A$ is the map $$\chi_i: S^1_i \hookrightarrow T$$ \noindent sending $t \in S^1$ to $tc_i \in M \otimes \mathbb{R}/2\pi\mathbb
    Z$, with image denoted $S^1_{\chi_i}$. This is injective because $c_i$ is primitive. More generally, given $J \subset A$ we will write $\chi_J: \prod_{j \in J} S^1_j \rightarrow T$ for the product morphism and $T_{\chi_J}$ for its image.
\end{defn}

Philosophically, since $c_i$ was used as a character on the $B$-side it plays the role of a cocharacter here.\\

The handle attachments we will do will, in general, require choices of pieces of conormals to the $T_I$ which must satisfy some compatibility in order to be sensible. To achieve this we will fix a \textit{dual system} for the given quiver $Q$: elements $c_j^\vee\in N$, one for each $v_j \in Q^{mut}$, satisfying $(c_i,c_j^\vee)=\delta_{ij}$. Such a thing exists because $Q$ is RFR (for example, by just extending the duals of the $c_i \in M$ to a basis of $N$). The following facts are then clear, and were largely re-derived or already used tacitly in the proof of Propositions \ref{chart lemma} and \ref{quiver geometry}.

\begin{lem}\label{dual system lemma}
    Assume given a dual system $\{c_j^\vee\}$ for an RFR quiver $Q$ and $I \subset Q^{mut}$ independent. For $J \subset I$, fix $C_J:=\bigoplus_{j\in J}\Z c_j$ and $B_{I,J}:=M_I\cap\bigcap_{j\in J}\ker c_j^\vee$. Then $C_J\subset M_I$ and $M_I=C_J\oplus B_{I,J}$. Further, if $I'\subset I$ and $J'\subset J$, then $B_{I,J}\subset B_{I',J'}$, $B_{I,J'}=C_{J\setminus J'}\oplus B_{I,J}$ and $B_{I',J'}\cap M_I=B_{I,J'}$.
\end{lem}

We will continue to fix the notation $C_J$ for spans of corresponding columns and $B_{I, J}$ for shared orthogonals restricted to $M_I := \bigoplus_{j \notin I} \mathbb{Z} \theta_j$ throughout. We will also abuse notation and write $B_{i} := B_{\{i\}, \{i\}}$. A final lemma, which was tacitly derived in Section 2 and whose proof is therefore omitted, pertains to the behavior of dual systems with respect to the splitting of Lemma \ref{product data decomposition}.

\begin{lem}
  Fix a dual system for an RFR quiver $Q$ and $v_i\in Q^{mut}$. Then $c_i^\vee|_{M_i}$ satisfies the hypotheses of Lemma \ref{product data decomposition} and the corresponding complement is $B_i=M_i\cap\ker c_i^\vee$. Further, under the resulting splitting $M_i\cong\mathbb Zc_i\oplus B_i$, the dual system $\{c_j^\vee|_{M_i}\}$ of $Q^{-i}$ vanishes on $\mathbb Zc_i$ and restricts to the dual system $\{c_j^\vee|_{B_i}\}$ of $Q^{/i}$.
\end{lem}

Given a cluster datum $(Q, \Sigma, A)$, $\sigma_I \in \Sigma$, and $J \subset I \cap A$, we write $\overline{T}_I := B_{I, I \cap A} \otimes \mathbb{R}/2\pi\mathbb{Z} \subset B_{I, J} \otimes \mathbb{R}/2\pi\mathbb{Z} := \overline{T}_{I, J}$. Note in this case that $T_{\chi_J} \subset T_I$ by the independence of $I$. Then, the point is that the preceding combinatorial lemmata allow us to compatibly split up $T$ into pieces parametrized by the $I \subset Q^{mut}$.  

\begin{cor}\label{compatible splitting}
    Assume given a dual system on an acyclic RFR quiver $Q$. Then, the map $\chi_J\colon(S^1)^J\to T_{\chi_J}$ is an isomorphism, and $T_I=T_{\chi_J}\times \overline{T}_{I,J}$, for any $J \subset I \cap A$. Further, these splittings are compatible in the sense that if $I'\subset I$ and $J'\subset J$, then $\overline{T}_{I,J}\subset \overline{T}_{I',J'}$.\\
    
    In fact $\overline{T}_{I,J'}=T_{\chi_{J\setminus J'}}\times \overline{T}_{I,J}$ and $\overline{T} _{I',J'}\cap T_I= \overline{T}_{I,J'}$, and so the inclusions $T_I\hookrightarrow T_{I'}$ send the $S^1_{\chi_j}$ factors to each other by the identity.
\end{cor}

\begin{proof}
    Tensor the preceding two lemmas with $\mathbb{R}/2\pi\mathbb{Z}$.
\end{proof}

Now we define a skeleton for a prospective $\mathcal{X}(Q)_{\Sigma, A}$.

\begin{con}[Skeleta of Cluster Mirrors]\label{Weinstein from Fan}

Fix a cluster datum $(Q, \Sigma, A)$. Given $J \subset A$, we write $D^J:=\prod_{j\in J}D^2_j$ and identify $\partial D^2_j$ with $S^1_{\chi_j}$. For $\sigma_I\in\Sigma$ with $J=I\cap A$ and $K=I\setminus A$, the \emph{cell} of $\sigma_I$ relative to $A$ is $$
\mathcal S_I:=D^J\times \overline{T}_I\times\sigma_K.$$
Let $I'\subset I$. The face of $\mathcal S_I$ where $z_j\in\partial D^2_j$ for $j\in J\setminus J'$, and the cone coordinate lies in $\sigma_{K'}$, is
\[
D^{J'}\times T_{\chi_{J\setminus J'}}\times \overline{T}_I\times\sigma_{K'}=D^{J'}\times \overline{T}_{I,J'}\times\sigma_{K'}\subset D^{J'}\times \overline{T}_{I'}\times\sigma_{K'}=\mathcal S_{I'},
\]
by Corollary \ref{compatible splitting}. Here $J'=I'\cap A$ and $K'=I'\setminus A$. The \emph{skeleton} $\mathcal S(Q)_{\Sigma,A}$ is the union of the $\mathcal S_I$, glued along these inclusions.
\end{con}

 In the sequel we will write $\mathcal S_Q:=\mathcal S(Q)_{\Sigma_Q,\Sigma_Q(1)}$, keeping in mind that $\mathcal{S}(Q)_{\Sigma, A}$ is the core of a purported mirror to $\mathcal{A}(Q)_{\Sigma, A}$ and that $\mathcal{A}(Q)_{\Sigma_Q, \Sigma_Q(1)} \cong \mathcal{A}(Q)$ for any RFR acyclic $Q$. Also, when all $J$ are trivial, i.e. when $A = \emptyset$, the skeleton just constructed consists only of cones of the fan $\Sigma$: that is, that $\mathcal{S}(Q)_{\Sigma, \emptyset}$ agrees with the FLTZ Lagrangian $-\mathbb{L}_\Sigma$. It is also not hard to see that $dim (\mathcal{S}_I) = m + n$ and that, for a subfan $\Sigma_0 \subset \Sigma$, we have that $\mathcal{S}(Q)_{\Sigma_0, \Sigma_0(1) \cap A} = \bigcup_{\sigma_I \in \Sigma_0} \mathcal{S}_I$ is a closed subset of $\mathcal{S}(Q)_{\Sigma, A}$.\\

 Geometrically, the passage from a cell $\mathcal{S}_I$ in an FLTZ skeleton to the corresponding cell in $\mathcal{S}(Q)_{\Sigma, A}$ replaces, for each $j\in J$, the half-cylinder $S^1_{\chi_j}\times\mathbb R_{\ge0}\theta_j^\vee$, formed by the ray direction and the circle $S^1_{\chi_j}\subset T_I$ of Corollary \ref{compatible splitting}, with the disk $D^2_j$ capping $S^1_{\chi_j}$. These disk factors should be thought of as corresponding to the $\mathbb{C}$ factors in the V/D decomposition defined in Section 2, see Definition \ref{v/d defn}.\\

Allow us to briefly unravel the singular topology of $\mathcal{S}(Q)_{\Sigma, A}$. In this passage, as in \cite{FLTZ} and \cite{gammage2022mirror}, for a set $I$ of rays we write $\sigma_I^\perp:=T_I\subset T$ for the subtorus orthogonal to them. When $\sigma_I$ is a cone this is its annihilator. One way to think about $\mathcal S(Q)_{\Sigma,A}$ is as a quotiented bordering of the FLTZ skeleton, as we indicate in the following lemma (whose proof is left to the reader).
 
 \begin{lem}
     Let $\overline{\mathbb L}_{\Sigma,A}$ be $-\mathbb L_\Sigma=\bigcup_{\sigma_I\in\Sigma} \sigma_I^\perp \times\sigma_I$ with each active cone direction compactified: i.e., with $T_I\times\sigma_J\times\sigma_K$ replaced by $T_I\times[0,\infty]^J\times\sigma_K$ and the same gluings. Then $\mathcal S(Q)_{\Sigma,A}$ is the quotient of $\overline{\mathbb L}_{\Sigma,A}$ by the relation collapsing the $S^1_{\chi_j}$-orbits, $j\in A$, in $\{t_j=\infty\}$ to a point.
 \end{lem}
 
\begin{figure}[t]
    \centering
\includegraphics[width=.8\linewidth]{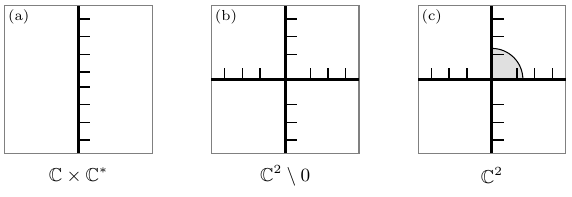}
    \caption{The FLTZ skeleta of the mirrors to $\mathbb C\times\mathbb C^*$ (a), $\mathbb C^2 \setminus 0$ (b), and $\mathbb C^2$ (c), drawn on the square torus.}
    \label{FLTZ skeleta}
\end{figure}

 The above elicits a graphical tool for drawing these (prospective) skeleta. Recall the standard way of drawing the FLTZ skeleton $-\mathbb L_\Sigma$ of $X^\vee_\Sigma$: namely, by drawing a square torus representation of $T$, stratified into the orthogonals $\sigma_I^\perp \subset T$ for each $\sigma_I \in \Sigma$, and with ``hair" in the shape of the conormal part $\sigma_I$ of the skeleton along each $\sigma_I^\perp$. For instance, Figure  \ref{FLTZ skeleta} shows the FLTZ skeleta of the mirrors to $\mathbb C\times\mathbb C^*$, $\mathbb C^2 \setminus 0$ and $\mathbb C^2$: each ray of the corresponding fan contributes a circle $\sigma_i^\perp$ to the skeleton, adorned with hair pointing in the direction of $\sigma_i$ itself.\\

The previous lemma tells us that $\mathcal S(Q)_{\Sigma,A}$ is obtained from the FLTZ picture by a single operation: put colloquially, that of choosing ``capped hairs’’. This suggests depicting $\mathcal S(Q)_{\Sigma,A}$ as a decorated FLTZ picture. Namely, we propose that in order to draw $\mathcal S(Q)_{\Sigma,A}$, we draw the FLTZ picture of $\Sigma$, color the span of the active rays in each conormal, and at each crossing $\sigma_{\{j,k\}}^\perp=\sigma_j^\perp\cap\sigma_k^\perp$ with $j\in A$ write the number $|\tilde E_{j,k}|$ of arrows between $v_j$ and $v_k$ (leaving it blank when they are not adjacent). A colored hair is read as a half-line whose $c_j$-circles are capped at its far end, and the stop of $\mathcal \X(Q)_{\Sigma,A}$ can be read off as the set of far ends of the uncolored hairs.\\

The labels record the only choice beyond the colored FLTZ picture, namely which circle in $\sigma_j^\perp$ is capped, and they record it through its position relative to the other $\sigma_k^\perp$. Indeed, at a crossing labelled $|\tilde E_{j,k}|$, each capped $c_j$-circle crosses $\sigma_k^\perp$ transversally in that many points. At a blank crossing, on the other hand, the capped $c_j$-circles lie in $\sigma_{\{j,k\}}^\perp$. In this case, if both rays are colored then the two families of caps there assemble into $S^3\times\overline T_{\{j,k\}}$ via the standard genus-one Heegaard splitting, and if, further, a fan is drawn there (i.e. $\sigma_{\{j,k\}}\in\Sigma$), then this $S^3$ is filled to $B^4\times\overline T_{\{j,k\}}$.\\

\begin{figure}[t]
    \centering
\includegraphics[width=.9\linewidth]{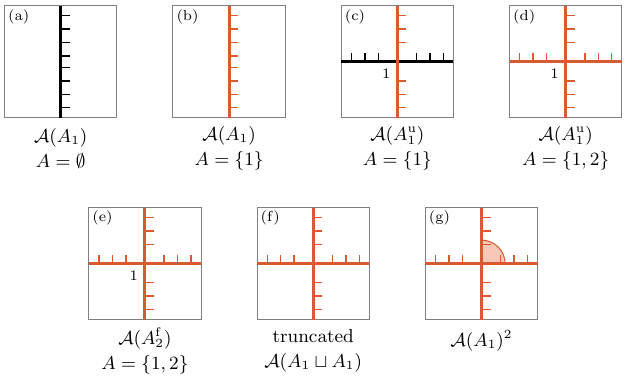}
    \caption{The decorated FLTZ skeleta of some $\mathcal S(Q)_{\Sigma,A}$, where orange colored hair is active. Top row: rank $2$. Bottom row: rank $4$, with one frozen vertex for each mutable vertex, projected onto the $2$-torus of the coordinates dual to the two rays.}
    \label{cluster mirror skeleta}
\end{figure}

Figure \ref{cluster mirror skeleta} illustrates these features in examples. In (a), which has no active rays, the cluster model $\mathcal A(A_1)_{\Sigma,\varnothing}$ is the toric variety $\mathbb C\times\mathbb C^*$ and the picture agrees with that of Figure \ref{FLTZ skeleta}. Coloring its hair corresponds to capping a meridian, which gives precisely the skeleton of Figure \ref{A1 skeleton} in (b). No stop, i.e. no uncolored hair, remains, since its mirror $\mathcal A(A_1)$ has no boundary.\\

Panels (c) and (d) treat $A_1^{\mathrm u}$, the quiver $v_1 \to v_2$ obtained from $A_1$ by unfreezing $v_2$. Doing this adds the ray $\rho_2$ to the corresponding fan. Since $v_1$ and $v_2$ are adjacent, the crossing is labeled $1$ and carries no fan. In (c), the capped circle $\sigma_1^\perp$ meets $\sigma_2^\perp$ once, so the stop (which we think of as mirror to the remaining boundary divisor $D_2$) meets the cap in a single boundary point. The resulting skeleton is $T^2$ with a disk glued along a meridian and a conormal cylinder glued along a latitude. In (d), the meridian and latitude are both capped and the two disks touch in one point. The skeleton is thus homotopy equivalent to $S^2$, and its prospective mirror is the full cluster variety $\mathcal{A}(A_1^{\mathrm{u}})$.\\

The bottom row shows why the labels are needed. Its examples have rank $4$, with one frozen vertex for each mutable vertex. To draw them we project $T$ onto the $2$-torus of the coordinates dual to the two rays, forgetting the frozen directions in which the capped circles also wind and which the labels then record. Notably, the three resulting pictures agree except in their labels and fans. In (e), $A_2^{\mathrm f}$ denotes the quiver with an arrow $v_1 \to v_2$ between mutable vertices, with frozen vertices $v_3$ and $v_4$, and with arrows $v_2\to v_3$ and $v_1\to v_4$, so that $c_1=-\theta_2-\theta_4$ and $c_2=\theta_1-\theta_3$. Each capped $c_1$-circle and each capped $c_2$-circle meets the $2$-torus $\sigma_{\{1,2\}}^\perp$ exactly once. Therefore the two families of disks touch pairwise, each $1$-disk meeting exactly one $2$-disk in a single boundary point, in a family of pairs parametrized by the $2$-torus $\sigma_{\{1,2\}}^\perp$, which is transverse to both families of circles rather than foliated by them. In (f), the vertices are not adjacent, and the two families of caps over the blank crossing \textit{do} assemble into an unfilled $S^3$, as described above.
In (g), the existence of a higher dimensional active cone fills this sphere to a $B^4$, and the skeleton is $(T^2\cup D^2)^2\simeq(S^2\vee S^1)^2$. In particular, without its label (e) would be indistinguishable from (f).\\

Now we define the actual stopped Weinstein manifold which is mirror to a cluster model $\mathcal{A}(Q)_{\Sigma, A}$. The sectors discussed in its construction will have certain collaredness properties which allow the construction to proceed at each step. To name them, given $A \subset \Sigma(1)$ we define a simplicial complex $\Delta_{\Sigma, A}$ by $$\Delta_{\Sigma, A} := \{K\:|\:\sigma_K \in \Sigma, K \cap A = \emptyset\}.$$

\noindent Note that $\Delta_{\Sigma, A}$ can be identified with the set of $\sigma_K \in \Sigma$ with $K$ satisfying the same. We also have that $\Delta_{\Sigma, \emptyset} = \Delta_\Sigma$, and that the simplices of $\Delta_{\Sigma, A}$ index the boundary strata $D_K$ of $\mathcal{A}(Q)_{\Sigma, A}$. Note that $link_{\Delta_{\Sigma, A}}(K) = \Delta_{\Sigma_{-K}, A_{-K}}$.

 \begin{con}[Mirrors to Cluster Models]\label{cluster mirrors} Fix an ordering $i_1<\dots<i_r$ of $A$, and put $A_s:=\{i_1,\dots,i_s\}$. For $\sigma_K\in\Sigma$ and $i\in\Sigma(1)$, put $$D_K^\vee:=\mathcal X(Q^{-K})_{\Sigma_{-K},A_{-K}} \quad \text{and} \quad H_i^\vee:=\mathcal X(Q^{/i})_{\Sigma_{/i},A_{/i}} .$$

These data have fewer vertices than $Q$, so these sectors are defined by induction on $m+n$. In the following construction the active set changes from stage to stage, and any use of the notation $D_K^\vee$ and $H_i^\vee$ in it is always formed with respect to the active set at the current step. For $i=i_s$ this causes no ambiguity in $D_i^\vee$ or $H_i^\vee$: since $\rho_i\notin\Sigma_{-i}(1)$, we have $(A_{s-1})_{-i}=(A_s)_{-i}$. For $\sigma_{K\cup i}\in\Sigma$, write $H^\vee_{i,K}$ for the stratum of $H_i^\vee$ at $K$ and $H_{i,K}$ for the corresponding stratum of $H_i$. Finally, we fix the notation $h_i := H^\vee_i \times T^* D^2$ throughout.\\

We mirror Construction \ref{GHK+} step by step, beginning with the mirror to a toric variety and performing a Weinstein handle attachment every time that construction performed an affine blowup:
$$X_0 := X_\Sigma^\vee \:\: \text{and} \:\:X_s :=X_{s-1}\cup_{i_s} h_i.$$

Here $\cup_{i_s}$ is the collared handle attachment of Definition \ref{collared handles}. Then, by definition we set $\mathcal \X(Q)_{\Sigma,A}:=X_r$ and $\mathcal \X(Q):=\mathcal \X(Q)_{\Sigma_Q,\Sigma_Q(1)}$.
\end{con}

\begin{prop}\label{A-side product splitting}
    Each handle attachment in Construction \ref{cluster mirrors} satisfies the hypotheses of Definition \ref{collared handles}, and so $\mathcal X(Q)_{\Sigma,A}$ is well defined. Moreover, for each $s$ we have 1) that $X_s$ is $\Delta_{\Sigma,A_s}$-collared, and its symplectic boundary stratum at $K\in\Delta^+_{\Sigma,A_s}$ is $F_K\cong D^\vee_K$, and 2) that if the datum is isomorphic to the product of $(\mathbb Zc,\{0\},\emptyset)$ with a datum $\mathcal D$, then $X_s\cong T^*S^1\times X_s(\mathcal D)$, compatibly with collars.
\end{prop}

\begin{proof}
   To prove the above claims we induct on $m + n$ and iterate on $s$. By the induction on $m+n$, the claims hold for every $D_K^\vee$ and $H_i^\vee$. Then at step 0, where $X_0 = X_\Sigma^\vee$, we see that claim $1)$ follows from Proposition \ref{lem:toric-collared}. Now assume we are at step $s$ and set $i := i_s$.  Then $1)$ for $X_{s-1}$ gives $F_i\cong D_i^\vee$. By Lemma \ref{product data decomposition}, the datum of $D_i^\vee$ is isomorphic to the product of $(\mathbb Zc_i,\{0\},\emptyset)$ with that of $H_i^\vee$. So $2)$ for $Q^{-i}$, which has fewer vertices, gives
$$F_i\cong T^*S^1_{\chi_i}\times H_i^\vee$$
compatibly with collars.\\

By Lemma \ref{handle ends} we learn, since $(\Delta_{\Sigma,A_{s-1}})_{\setminus i}=\Delta_{\Sigma,A_s}$, that $X_s$ is $\Delta_{\Sigma,A_s}$-collared. If $\sigma_{K\cup i}\notin\Sigma$, the stratum at $K$ is unchanged: that is, since $\rho_i\notin\Sigma_{-K}(1)$ we have $(A_s)_{-K}=(A_{s-1})_{-K}$, so it is still $D_K^\vee$.
On the other hand, if $\sigma_{K\cup i}\in\Sigma$ then the stratum at $K$ is $D_K^\vee$ for $A_{s-1}$, with $H^\vee_{i,K}\times T^*D^2$ attached at $i$. By $1)$ for $H_i^\vee$, and since $(Q^{/i})_{-K}=(Q^{-K})_{/i}$, while both fans are $\Sigma_{-(K\cup i)}$ and both active sets are $A_{-(K\cup i)}$, $H^\vee_{i,K}$ is the $H_i^\vee$ of the datum of $D_K$. So this attachment is the last step in the construction of $D_K^\vee$ for $A_s$, using the ordering induced on $(A_s)_{-K}$, in which $i$ comes last. This shows $1)$ for $X_s$.\\

Finally, for a product datum, the fan lies in the factor $\mathcal D$, so $X_0=T^*S^1\times X_0(\mathcal D)$. Inductively, the face at $i$ of $X_{s-1}$ is $T^*S^1$ times that of $X_{s-1}(\mathcal D)$. Contracting $i$ commutes with taking the product, so by $2)$ in smaller rank, $H_i^\vee$ is $T^*S^1$ times the $H_i^\vee$ of $\mathcal D$. Hence step $s$ attaches $T^*S^1$ times the handle of step $s$ for $\mathcal D$. Attaching a product handle commutes with taking the product with $T^*S^1$ so $X_s\cong T^*S^1\times X_s(\mathcal D)$, compatibly with collars. This implies $2)$. 
\end{proof}

It is not hard to see that $\mathcal{X}(Q)_{\Sigma, A}$ is a stably polarizable, and hence graded, stopped Weinstein manifold. To see this, give $X_0=T^*T$ and each $T^*D^2$ their cotangent-fiber polarizations, and give $H_i^\vee$ its polarization by induction. By claim $2)$ of Proposition \ref{A-side product splitting} and by Corollary \ref{compatible splitting}, the face $F_i\cong T^*S^1_{\chi_i}\times H_i^\vee$ inherits from $X_{s-1}$ the product of the fiber polarization on $T^*S^1_{\chi_i}$ with that of $H_i^\vee$, since the splitting of $T_i$ is linear. This is the restriction of the polarization of $h_i$, so the two glue.

\begin{rem}\label{naturality}
    We claim that $\mathcal \X(Q)_{\Sigma,A}$ does not depend on the ordering of $A$ or the dual system (up to Weinstein homotopy), i.e. that it depends only on the input cluster datum -- i.e., on $M$, $\Sigma$, $A$, and the columns $c_i$ for $i\in A$. Further, the construction is natural, via cotangent lifts, for isomorphisms of data, and it takes products of data to products. The claim follows quickly from unraveling the construction except for invariance under changing the ordering on $A$, and to see this it suffices to show invariance under swapping two consecutive steps $i$ and $k.$\\

    If $\sigma_{\{i,k\}}\notin\Sigma$, then $\mathcal N_i\cap\mathcal N_k=\emptyset$, and neither handle changes, since $\rho_i\notin\Sigma_{-k}(1)$. So both orders glue the same handles along disjoint regions. On the other hand, if $\sigma_{\{i,k\}}\in\Sigma$, then by the preceding remark and Lemma \ref{handle ends} for $H_k^\vee$, attaching at $i$ first enlarges $h_k$ by $H^\vee_{\{i,k\}}\times T^*D^2_i\times T^*D^2_k$, and symmetrically. These corner pieces agree, because $(Q^{/k})_{/i}=(Q^{/i})_{/k}$: both are $Q^{-\{i,k\}}$ with frozen lattice quotiented by $\mathbb Zc_i\oplus\mathbb Zc_k$, which is saturated since $Q$ is RFR. In either case one receives the same sector-with-corners before or after the swap.
\end{rem}

Constructions \ref{Weinstein from Fan} and \ref{cluster mirrors} are compatible.

\begin{prop}\label{issa skeleton}
    The skeleton of $\mathcal{X}(Q)_{\Sigma, A}$ is $\mathcal{S}(Q)_{\Sigma, A}$.
\end{prop}

\begin{proof}[Proof Sketch]
    Induct on $s$. For $s=0$, the relative skeleton of $X^\vee_\Sigma$ is $-\mathbb{L}_\Sigma=\mathcal S(Q)_{\Sigma,\emptyset}$. At step $s$, the skeleton changes only inside $N_i$, where it is $core(F_i)\times(\text{half-line})=S^1_{\chi_i}\times core(H_i^\vee)\times(\text{half-line})$. The handle replaces this by $core(H_i^\vee)\times D^2$, as for the single $A_1$ handle considered in Section 4 of\cite{GL}. On a cell with $\rho_i\in\sigma_I$, this replaces the circle $S^1_{\chi_i}$ split off by Corollary \ref{compatible splitting}, times the ray $\mathbb R_{\ge0}\theta_i^\vee$, by a disc. This is exactly the passage from $\mathcal S_I$ for $A_{s-1}$ to $\mathcal S_I$ for $A_s$.
\end{proof}

In the sequel we write $\partial^\infty\mathcal S(Q)_{\Sigma,A}$ for the Legendrian at infinity of $\mathcal S(Q)_{\Sigma,A}$, or equivalently the stop of $\mathcal \X(Q)_{\Sigma,A}$. By Definition \ref{collared} and the previous construction, $\partial^\infty\mathcal S(Q)_{\Sigma,A}$ is the core of $F_D$. Further, for $A=\emptyset$, $\partial^\infty\mathcal S(Q)_{\Sigma,A}=\Lambda_\Sigma$. When the input datum is clear, we will often suppress notation further and simply write $\partial^\infty \mathcal{S}$ for the same Legendrian.\\

We end this section by describing an $A$-sided analogue of the local form of cluster models established in Proposition \ref{chart lemma} and giving some examples.
 
\begin{cor} \label{cor:charts}
Assume given $\sigma_I\in\Sigma$, and put $J:=I\cap A$, $K:=I\setminus A$. Then $\mathcal S(Q)_{\langle\sigma_I\rangle,J}=\bigcup_{I'\subseteq I}\mathcal S_{I'}$ and $$\mathcal \X(Q)_{\langle\sigma_I\rangle,J}\cong\mathcal X(A_1)^J\times\bigl((\mathbb C\times\mathbb C^*)^\vee\bigr)^K\times T^*\overline{T}_{I,I},$$ which is the mirror of Proposition 2.28 factor by factor. Here $\mathcal X(A_1)$ is the mirror of $\mathcal A(A_1)$ from Section 4 of \cite{GL}, and $(\mathbb C\times\mathbb C^*)^\vee=(T^*S^1,\infty)\times T^*S^1$ is its toric model.
\end{cor}
 
\begin{proof}
The fan $\langle\sigma_I\rangle$ has rays $\rho_i$, $i\in I$, and active set $J$, so only the handles $h_j$, $j\in J$, are attached. By Lemma \ref{dual system lemma},

$$M=\bigoplus_{i\in I}(\mathbb Z\theta_i\oplus\mathbb Zc_i)\oplus B_{I,I},$$

and each ray $\theta_i^\vee$ pairs to zero with every $c_{i'}$ (since $I$ is independent) and with $B_{I,I}$. So the datum $(Q,\langle\sigma_I\rangle,J)$ splits as a product, with one factor for each $i\in I$ and one for $B_{I,I}$: for $j\in J$, the factor on $\mathbb Z\theta_j\oplus\mathbb Zc_j$, with the ray $\theta_j^\vee$ active, is the datum of $\mathcal A(A_1)$ with frozen variable $x^{c_j}$;
for $k\in K$, the same factor with the ray inactive is the datum of its toric model $\mathbb C\times\mathbb C^*$;
$B_{I,I}$ carries the trivial datum.\\

Replacing $c_j^\vee$ by $c_j^\vee-\sum_{i\in I}(\theta_i,c_j^\vee)\theta_i^\vee$ gives a dual system that respects this splitting. With it, Construction \ref{cluster mirrors} is the product of the constructions for the factors, by Remark \ref{naturality}. The statement about skeleta follows from the discussion after Construction \ref{Weinstein from Fan}.
\end{proof}

\begin{rem}
    The cells $\mathcal S_I$ are in bijection with the preimages in $\mathcal A(Q)_{\Sigma,A}$ of the torus orbits of $X_\Sigma$, with each disk factor $D^2_j$ corresponding to the affine factor $\mathbb C_{x'_j}$. For $\Sigma=\Sigma_Q$ and $A=\Sigma_Q(1)$ these are the V/D strata of Theorem 2.15. In particular, $\#\mathcal A(Q)_{\Sigma,A}(\mathbb F_p)=\sum_{\sigma_I\in\Sigma}p^{|J|}(p-1)^{m+n-|I|-|J|}$.
\end{rem}

Assume given $S\subset Q_{\mathrm{mut}}$. We note for future use that removing the pieces $l_k$ with $k\in S\setminus A$ from the stop of $\mathcal X(Q)_{\Sigma,A}$, as well as removing the cores of the handles $h_k$ with $k\in S\cap A$, gives $\mathcal X(Q^{\setminus S})_{\Sigma_{\setminus S},A_{\setminus S}}$, compatibly with collars. In particular, for $k\in A$, $\mathcal X(Q^{\setminus k})_{\Sigma_{\setminus k},A_{\setminus k}}$ is $\mathcal X(Q)_{\Sigma,A\setminus\{k\}}$ with $l_k$ removed. This is clear from the construction.\\

For $K\in\Delta_{\Sigma,A}$, write
$$\mathcal S^K_{\Sigma,A}:=\bigcup_{\sigma_I\in\Sigma,\ I\supseteq K}\mathcal S_I\ \subset\ \mathcal S(Q)_{\Sigma,A}$$

\noindent for the union of the cells of the cones containing $\sigma_K$. For $K\subset K'$ there are inclusions
$\mathcal S^{K'}_{\Sigma,A}\hookrightarrow\mathcal S^K_{\Sigma,A}.$\\

On the other hand, one may also consider the skeleton $\mathcal S(Q^{-K})_{\Sigma_{-K},A_{-K}}$ whose strata are indexed by the cones $\sigma_{I\setminus K}$ of $\Sigma_{-K}$; this set of cones is clearly in bijection with the indexing set of the union defining $\mathcal{S}^K_{\Sigma, A}$. Since $\Sigma$ is simplicial with coordinate rays, $\sigma_{I\setminus A}=\sigma_K\times\sigma_{(I\setminus K)\setminus A}$. Combining these with Corollary \ref{compatible splitting} gives inclusions
$$\mathcal S(Q^{-K'})_{\Sigma_{-K'},A_{-K'}}\times\sigma_{K'}=\bigl(\mathcal S(Q^{-K'})_{\Sigma_{-K'},A_{-K'}}\times\sigma_{K'\setminus K}\bigr)\times\sigma_K\hookrightarrow\mathcal S(Q^{-K})_{\Sigma_{-K},A_{-K}}\times\sigma_K .$$

Taking $K=\emptyset$, the image of such an inclusion agrees with the pieces defined above.

\begin{lem}
For $K\in\Delta_{\Sigma,A}$, the image of $\mathcal S(Q^{-K})_{\Sigma_{-K},A_{-K}}\times\sigma_K\hookrightarrow\mathcal S(Q)_{\Sigma,A}$ is $\mathcal S^K_{\Sigma,A}$. Under this identification the two inclusions just defined agree. Moreover, near $\mathcal S^K_{\Sigma,A}\cap\{\theta_k^\vee>0,\ k\in K\}$, the sector $\mathcal \X(Q)_{\Sigma,A}$ is the end
$$\mathcal N_K\cong D_K^\vee\times(T^*[0,1])^K,\qquad D_K^\vee=\mathcal X(Q^{-K})_{\Sigma_{-K},A_{-K}},$$

and these ends are nested accordingly.
\end{lem}

\begin{proof}
Fix $\sigma_I\in\Sigma$ with $I\supseteq K$, and put $I':=I\setminus K$. The cell of $I'$ for the deleted data is $D^J\times\overline T^{(-K)}_{I'}\times\sigma_{I'\setminus A}$, where $\overline T^{(-K)}_{I'}$ is built from $M_K\cap M_{I'}\cap\bigcap_{j\in J}\ker c_j^\vee=B_{I,J}$. So $\overline T^{(-K)}_{I'}=\overline T_I$, and the cell of $I'$ times $\sigma_K$ is $\mathcal S_I$. The faces match by Corollary \ref{compatible splitting}. The statement about ends is claim 1) of Construction \ref{cluster mirrors}, together with the nesting of ends after Definition \ref{collared}.
\end{proof}

We can rephrase this as a statement about a cover of the stop $\partial^\infty\mathcal S(Q)_{\Sigma,A}$.

\begin{cor}\label{stop cover}
   For $K\in\Delta^+_{\Sigma,A}$, let $\Omega_K$ be the open part of $\partial^\infty\mathcal S(Q)_{\Sigma,A}$ where $\theta_k^\vee>0$ for all $k\in K$. Then the $\Omega_k$ cover the stop, $\Omega_K=\bigcap_{k\in K}\Omega_k$, and
$$\Omega_K\cong\mathcal S(Q^{-K})_{\Sigma_{-K},A_{-K}}\times\operatorname{int}\partial^\infty\sigma_K,$$

with the inclusions $\Omega_{K'}\subset\Omega_K$ induced by those described in the previous lemma.
\end{cor}

Correspondingly, the stop is the core of the ribbon $F_D$ glued from the faces $F_k\cong D_k^\vee$ along the $F_K\cong D_K^\vee$. This is the mirror of the fact that the boundary $D$ of $\mathcal \A(Q)_{\Sigma,A}$ is the union of the strata $D_K\cong\mathcal A(Q^{-K})_{\Sigma_{-K},A_{-K}}$. By construction active ray directions do not carry any stops in $\mathcal{X}(Q)_{\Sigma, A}$, and so one sees (for example) that the stop of $\mathcal{X}(Q)$ is empty.
\section{HMS for Cluster Models}

%%We end this subsection by giving some examples of outputs of Construction \ref{Weinstein from Fan}. These are largely logically independent from the rest of the paper but illustrate its use.

%%\begin{exmp}
    %%Let $Q$ be any quiver, and write $\Sigma_Q(1)$ for the fan consisting only of rays $\rho_i$ as defined in Definition \ref{Q-fan}. Assume $A_{\Sigma_Q(1)} = \Sigma_Q(1)$. Then the content of \cite{GL} is that Construction \ref{Weinstein from Fan} applied to $(\Sigma_Q(1), A_{\Sigma_Q(1)})$ is exactly the mirror to truncated cluster variety corresponding to $Q$.
%%\end{exmp}

%%\begin{exmp}
    %%If $\Sigma$ is any smooth fan and $A = \emptyset$, then Construction \ref{Weinstein from Fan} applied to $(\Sigma, A)$ outputs the FLTZ mirror to the toric variety $X_\Sigma$. On the other hand, for $\Sigma_Q$ arising from a really full rank acyclic cluster variety $Q$, $(\Sigma_Q, A = (\Sigma_Q)_1$ is what we will prove to be the mirror to $\mathcal{A}(Q)$.
%%\end{exmp}

%%The cases where $A$ is not all or no rays will be important in the proofs of our main theorems. We also expect them to be of independent interest.
Comparing our two main constructions (Constructions \ref{GHK+} and \ref{cluster mirrors}), one already sees a panoply of steps and structures which behave quite similarly on the two sides of mirror symmetry. We tabulate a few important ones now for reference. Each row of the following table pairs an object on the $A$-side with its mirror on the $B$-side.\\
\begin{table}[ht]
\centering
\renewcommand{\arraystretch}{1.3}
\begin{tabular}{ll}
\hline
$A$-side & $B$-side \\
\hline
toric mirror $X_\Sigma^\vee=(T^*T^n,\partial^\infty(-\mathbb{L}_\Sigma))$ & toric variety $X_\Sigma$ \\
handle attachment of $h_i=H_i^\vee\times T^*D^2$ & affine blowup along $H_i$ \\
stop/handle removal $\mathfrak{sr}_S$ & freezing at $S$ (open restriction) \\
end $N_K\cong D_K^\vee\times(T^*[0,1])^K$ & boundary stratum $D_K\cong\mathcal{A}(Q^{-K})_{\Sigma_{-K},A_{-K}}$ \\
base $H_i^\vee=\mathcal{X}(Q^{/i})_{\Sigma_{/i},A_{/i}}$ & center $H_i\cong\mathcal{A}(Q^{/i})_{\Sigma_{/i},A_{/i}}$ \\
$\mathcal{X}(A_1)^J\times(T^*S^1,\infty)^K\times T^*T^r$ & chart $\mathcal{A}(A_1)^J\times\mathbb{C}^K\times(\mathbb{C}^*)^r$ \\
\hline
\end{tabular}
\end{table}

Given a cluster model $\A$ with mirror $\X$, the proof of HMS proceeds by plugging the iterative construction of the former by affine blowups and of the latter by handle attachments into the technology developed in \cite{GL}. In particular, we will repeatedly match affine blowups on the $B$-side to handle attachments on the $A$-side (both of which, by \cite{GL} and \cite{GPSI}, induce pushouts on categories) by identifying the corresponding spans as diagrams. More precisely still, given $i\in A$, write $A^-:=A\setminus\{i\}$, and write $\mathcal A^-:=\mathcal A(Q)_{\Sigma,A^-}$ and $\mathcal X^-:=\mathcal X(Q)_{\Sigma,A^-}$ for the cluster model and mirror with the same $Q$ and $\Sigma$ but active set $A^-$. We exhibit an equivalence of diagrams between
$$
\Fuk(\mathcal{X}^-)\xleftarrow{\,\iota\,}\Fuk(F_i)\xrightarrow{\,r\,}\Fuk(h_i) \quad \text{and}\quad \Coh(\mathcal{A}^-)\xleftarrow{\,(k^*)^L\,}\Coh(D_i)\xrightarrow{\,i^*_{H_i}\,}\Coh(H_i),
$$
which induces an equivalence between their pushouts. By construction these are $\Fuk(\X)$ and $\Coh(\A)$, completing the proof.\\

Recall that the left legs of these spans are pushforward along the end corresponding to a boundary divisor (on the $A$-side) and the left adjoint of restriction to that divisor (on the $B$-side), which are the two sides of the recursion along the boundary considered in Section 4.3 of \cite{gammage2022mirror} extended to cluster models. For this reason the proof's induction hypothesis must say, as in Section 7.2 of \cite{gammage2022mirror}, that restriction to boundary divisors is mirror to microlocalization to their ends. Passing to left adjoints, as in Theorem 7.4.1 of \cite{gammage2022mirror}, then matches the left legs. For this to survive to the next step, however, it must hold for all iterated strata at once and coherently, as the next spans are built from those strata. We therefore state the theorem as an equivalence of functors on the poset $\Delta_{\Sigma,A}^{op}$, not just of individual spans.

\subsection{The Main HMS Theorem} 

Given a cluster datum $(Q, \Sigma, A)$, consider again the poset $\Delta_{\Sigma, A} = \{\sigma_K \in \Sigma \:|\: K \cap A = \emptyset\}$, which we think of as the poset of inactive cones in $\Sigma$. For $K \in \Delta_{\Sigma, A}$, we also consider again the two spaces $D_K := \mathcal{A}(Q^{-K})_{\Sigma_{-K}, A_{-K}}$ and $F_K \cong D_K^\vee := \mathcal{X}(Q^{-K})_{\Sigma_{-K}, A_{-K}}$ studied in Sections 2 and 3, respectively, as well as the end $N_K \cong F_K \times (T^*[0,1])^K$. Then, one receives two functors $$\Coh(D_\bullet), \Fuk(D_\bullet^\vee): \Delta_{\Sigma, A}^{op} \rightarrow Cat_\infty$$ \noindent defined by sending $K \subset K'$ to $(k^*_{K\subset K'})^L: \Coh(D_{K'}) \rightarrow \Coh(D_K)$, the left adjoint of restriction with respect to the natural closed immersion $D_{K'} \hookrightarrow D_K$ (a shifted pushforward, by Remark \ref{GV duality remark}), and to $\iota_{K \subset K'}: \Fuk(D_{K'}^\vee) \rightarrow \Fuk(D_K^\vee)$, the corresponding sectorial pushforward.\\

To establish mirror symmetry we first show that the functors we are considering out of $\Delta_{\Sigma, A}^{op}$ are, in the first place, coherently compatible with the relevant pushout formulae. In the sequel, given $i\in A$ we write $A^-:=A\setminus\{i\}$ and $D^-_K$, $(D_K^-)^
\vee$ for the $K$-th stratum of $\mathcal{A}(Q)_{\Sigma,A^-}$ and  $\mathcal{X}(Q)_{\Sigma,A^-}$, respectively. Write $link(i)$ for the link of $i$ in $\Delta_{\Sigma,A^-}$, i.e. for the simplicial complex of those $K\in\Delta_{\Sigma,A}$ with $\sigma_{K\cup i}\in\Sigma$. We regard $D^-_\bullet$ as a functor on $\Delta_{\Sigma,A}^{op}$ by restriction from $\Delta_{\Sigma,A^-}^{op}$. Similarly, we regard $D^-_{\bullet\cup i}$, $H_{i,\bullet}$, $H^\vee_{i,\bullet}$ and $h_i\cap N_\bullet$ as functors on $link(i)^{op}$, and in the sequel we will also extend them by $0$ to functors on $\Delta_{\Sigma,A}^{op}$. Here $H_{i,K}$ and $H^\vee_{i,K}$ are the $K$-strata of $H_i$ and $H_i^\vee$, as in Construction 3.20.

\begin{lem}\label{lem:Bglue}
The functor $\Coh(D_\bullet)$ is the pushout, in the category $\mathrm{Fun}(\Delta^{op}_{\Sigma,A},\mathrm{Cat}_\infty)$, of the span
 $$
\Coh(D^-_\bullet)\xleftarrow{(k^*)^L}\Coh(D^-_{\bullet\cup i})\xrightarrow{i^*_{H_i}}\Coh(H_{i,\bullet}),$$
where all transition functors are the left adjoints of restrictions. The cocone realizing this pushout is given componentwise by $\pi^*$ and $\iota_{E*}q^*[-1]$ (see after Theorem \ref{affine blowup}).
\end{lem}
 
\begin{proof}
By Proposition \ref{quiver geometry} applied to the deleted data we see both that $D^-_K\cap D_i=D^-_{K\cup i}$ transversally when $\sigma_{K\cup i}\in\Sigma$ and is empty otherwise, and that $D_K=\mathrm{Bl}^{\mathrm{aff}}_{D^-_K\cap H_i}D^-_K$. So pointwise the statement is Theorem \ref{affine blowup} (in its form where the functor from $\Coh(D)$ to $\Coh(X)$ is taken to be $(k^*)^L\simeq k_*[-1]$). The legs are natural in $\bullet$ by base change, since all squares involved are cartesian and Tor-independent (being products in the charts of Proposition 2.28), so the span is a diagram in $\mathrm{Fun}(\Delta^{op}_{\Sigma,A}, \mathrm{Cat}_\infty)$. Since colimits there are computed pointwise, the claim follows from the pointwise statement.
\end{proof}
 
\begin{lem}\label{lem:Aglue}The functor 
$\Fuk(D^\vee_\bullet)$ is the pushout, in the category $\mathrm{Fun}(\Delta^{op}_{\Sigma,A},\mathrm{Cat}_\infty)$, of the span
$$
\Fuk((D^-_\bullet)^\vee)\xleftarrow{\ \iota\ }\Fuk((D^-_{\bullet\cup i})^\vee)\xrightarrow{\ r\ }\Fuk(h_i\cap  N_\bullet).
$$
Here $h_i\cap  N_K\cong H^\vee_{i,K}\times T^*D^2$ by Lemma \ref{handle ends}. The cocone realizing this pushout is given componentwise by $\jmath_*$ and $\kappa_*$.
\end{lem}
 
\begin{proof}
We need only show the above for $K\in\Delta_{\Sigma,A}$ with $\sigma_{K\cup i}\in\Sigma$, since otherwise $h_i$ misses $ N_K$ and there is nothing to prove. By Construction \ref{cluster mirrors} and Lemma \ref{handle ends} the $K$-stratum of $\mathcal{X}(Q)_{\Sigma,A}$ is the collared handle attachment $(D^-_K)^\vee\cup_i(H^\vee_{i,K}\times T^*D^2)$, all of whose pieces are Weinstein. So pointwise the claim is Theorem \ref{thm:handle-pushout}, using Corollary 4.9 of \cite{GL}. All functors are pushforwards along inclusions of subsectors, which commute with the end inclusions by Lemma \ref{handle ends}, and $\Fuk$ is ($\infty$-categorically) functorial for these by the main theorem of \cite{GPSI}. The claim then follows from the pointwise statement as in the previous lemma.
\end{proof}

The proof of HMS will be able to match the left legs of the above diagrams iteratively by induction, and so we are left as a final input to show that the right legs match as well.

\begin{lem}\label{lem:Hleg}
Let $i\in A$, and let $\Phi''\colon \Fuk(H^\vee_{i,\bullet})\simeq \Coh(H_{i,\bullet})$ be an equivalence of functors on $link(i)^{\mathrm{op}}$. Then, under the splittings
$$
(D^-_{\bullet\cup i})^\vee\cong T^*S^1_{\chi_i}\times H^\vee_{i,\bullet}\quad \text{and} \quad D^-_{\bullet\cup i}\cong\C^*_{y_i}\times H_{i,\bullet},$$
the right legs $r$ and $i^*_{H_i}$ of Lemmas~\ref{lem:Aglue} and~\ref{lem:Bglue} are intertwined. That is,
$$
i^*_{H_i}\circ(\Phi_{\C^*}\otimes\Phi'')\ \simeq\ (\Phi_{\mathrm{pt}}\otimes\Phi'')\circ r,
$$
naturally in $\bullet$, where $\Phi_{\C^*}$ denotes HMS for $\C^*$, and $\Phi_{\mathrm{pt}}$ the equivalence from $\Fuk(T^*D^2)\simeq \Coh(\mathrm{pt})$ utilized in 
\cite{GL}.
\end{lem}
 
\begin{proof}
The idea is that both right legs are the $A_1$ leg studied in \cite{GL} times the identity on the contracted data. Indeed, on the $B$-side we see that Proposition \ref{quiver geometry} applied to the deleted data gives $D^-_{K\cup i} \cong \mathbb{C}^*_{y_i} \times H_{i,K}$ with $D^-_K\cap H_i=\{y_i=-1\}\times H_{i,K}$. On the $A$-side, the face at $i$ is $F_i\cong T^*S^1_{\chi_i}\times H_i^\vee$, compatibly with collars. Therefore, restricting to strata gives $(D^-_{K\cup i})^\vee\cong T^*S^1_{\chi_i}\times H^\vee_{i,K}$, and so by Lemma \ref{handle ends} we have that $h_i\cap N_K\cong T^*D^2\times H^\vee_{i,K}$. Under these splittings, we see that
$$
i^*_{H_i}=\rho\boxtimes\mathrm{id},\qquad r=r_{A_1}\times\mathrm{id},
$$
\noindent where $\rho\colon \Coh(\C^*)\to \Coh(\mathrm{pt})$ is restriction to $y_i=-1$ and $r_{A_1}\colon \Fuk(T^*S^1)\to \Fuk(T^*D^2)$ is the collar inclusion of the $A_1$ handle. K\"unneth then turns these into tensor products of functors.\\
 
We are left to show naturality in $
\bullet$. For $K\subset K'$ in $link(i)$, in both cases the relevant transition functors act on the second factor only: on the $B$-side this follows by Proposition 
\ref{chart lemma}, where base change holds because the squares are products and hence Tor-independent, and on the 
$A$-side this holds because the end inclusions respect the splitting, by compatibility with collars. So both sides are coherently of the form $(\text{fixed }A_1\text{ data})\otimes(\text{diagram over } link(i))$, where the former term is the $A_1$ square
$$
\rho\circ\Phi_{\C^*}\simeq\Phi_{\mathrm{pt}}\circ r_{A_1}.
$$
Tensoring a natural transformation with a fixed square preserves naturality, and so the claim follows.
\end{proof}

Now we may prove our main theorem.

\begin{thm}\label{HMS}
The spaces $\mathcal{X}(Q)_{\Sigma, A}$ and $\mathcal{A}(Q)_{\Sigma, A}$ are homological mirrors compatibly with their boundary strata. That is, for every cluster datum $(Q, \Sigma, A)$, there is an equivalence of functors on $\Delta_{\Sigma, A}^{op}$ $$\Phi: \Fuk(D_\bullet^\vee) \simeq \Coh(D_\bullet),$$ \noindent whose restriction to $\emptyset \in \Delta_{\Sigma, A}$ induces an equivalence $$\Phi: \Fuk(\mathcal{X}(Q)_{\Sigma, A}) \simeq \Coh(\mathcal{A}(Q)_{\Sigma, A}).$$

This equivalence satisfies the following additional properties.

\begin{enumerate}
    \item Given $K \in \Delta_{\Sigma, A}$, the restriction of $\Phi$ to the upper set $\{K \subset K'\} \cong \Delta_{\Sigma_{-K}, A_{-K}}$ agrees with the equivalence for the cluster datum $(Q^{-K}, \Sigma_{-K}, A_{-K})$. Therefore sectorial pushforward along the end $N_K 
    \hookrightarrow \mathcal{X}(Q)_{\Sigma, A}$ is mirror to the left adjoint $(k_K^*)^L\simeq k_{K*}[-|K|]$ of restriction to $D_K$, and its right adjoint is mirror to $k_K^*$.

    \item Blowups are mirror to handle attachments, naturally in $\Delta_{\Sigma, A}$. That is, writing $A^-:=A\setminus\{i\}$ given $i \in A$, and using the notation of Theorem \ref{affine blowup} for the affine blowup $\mathcal \A(Q)_{\Sigma,A}=\mathrm{Bl}^{\mathrm{aff}}_{H_i}\mathcal \A(Q)_{\Sigma,A^-}$ and that of Theorem \ref{thm:handle-pushout} for the handle attachment $\mathcal \X(Q)_{\Sigma,A}=\mathcal \X(Q)_{\Sigma,A^-}\cup_ih_i$, there are natural transformations of functors on $\Delta_{\Sigma,A}^{\mathrm{op}}$ 
$$\begin{aligned}&\pi^*\colon \Coh(D^-_\bullet)\Rightarrow \Coh(D_\bullet), &&\iota_{E*}q^*[-1]\colon \Coh(H_{i,\bullet})\Rightarrow \Coh(D_\bullet),\\ &\jmath_*\colon \Fuk((D^-_\bullet)^\vee)\Rightarrow \Fuk(D^\vee_\bullet), &&\kappa_*\colon \Fuk(h_i\cap\mathcal N_\bullet)\Rightarrow \Fuk(D^\vee_\bullet),\end{aligned}$$

\noindent which $\Phi$ intertwines:

$$\Phi\circ\jmath_*\simeq\pi^*\circ\Phi^-,\qquad\Phi\circ\kappa_*\simeq\iota_{E*}q^*[-1]\circ\Phi_h.$$

\noindent Here $\Phi^-:=\Phi_{Q,\Sigma,A^-}$ and $\Phi_h:=\Phi_{\mathrm{pt}}\otimes\Phi''$, where $\Phi''$ is the equivalence for $(Q^{/i},\Sigma_{/i},A_{/i})$ supplied by induction via $H_i\cong\mathcal A(Q^{/i})_{\Sigma_{/i},A_{/i}}$ and $H_i^\vee=\mathcal X(Q^{/i})_{\Sigma_{/i},A_{/i}}$. Diagrams for $A^-$ are restricted to $\Delta_{\Sigma,A}$, and those for $H_i$ and $h_i$ are extended by $0$ off $link(i)$.

\item $\Phi$ is compatible with lattice isomorphisms of the data, and for a product datum we have $\Phi_{\mathbb C^*\times\mathcal D}\simeq\Phi_{\mathbb C^*}\otimes\Phi_{\mathcal D}$ via K\"unneth.
\end{enumerate}
\end{thm}

\begin{proof}
    By induction on $m + n$, and then on $|A|$. When $A = \emptyset$ then every item of the above follows from standard results on toric HMS. In particular $(1)$ follows from the main theorem of \cite{gammage2022mirror}, $(2)$ is vacuous, and $(3)$ is the last sentence of Theorem \ref{toric HMS natural}.\\

   In general we set $
   \Phi$ to be the result of comparing the spans in Lemmas \ref{lem:Bglue} and \ref{lem:Aglue} termwise via components which are comparison functors that have been supplied by induction. That is, we set $\Phi$ to be the map of colimits of spans induced by the diagram $$
 \begin{tikzcd}
{\Coh(D_\bullet^-)} &  & \Coh(D_{\bullet \cup i}^-) \arrow[rr, "i^*_{H_i}"] \arrow[ll, "(k^*)^L"'] &  & \Coh(H_{i,\bullet}) \\
\Fuk((D_\bullet^-)^\vee) \arrow[u, "\Phi^-"]                          &  & \Fuk((D_{\bullet \cup i}^-)^\vee) \arrow[ll, "\iota"'] \arrow[rr, "r"] \arrow[u, "\Phi^-_{\bullet\cup i}"]                                                      &  & \Fuk(h_i \cap N_\bullet) \arrow[u, "\Phi_h"]
\end{tikzcd}.$$

Here $\Phi^-_{\bullet\cup i}$ denotes the restriction of $\Phi^-$ along $K\mapsto K\cup i$. Therefore, by the induction hypothesis and item (1) we see that the left square in the definition of $\Phi$ commutes, since its horizontal maps are the transition functors of the diagrams for $A^-$ along $K\cup i\to K$. We also see that said square has vertical legs which are equivalences. On the other hand, the right vertical arrow is actually an equivalence because of K\"unneth and induction. Therefore Lemma \ref{lem:Hleg} implies that the right square also commutes, and so the theorem itself holds for $
\Phi$.\\

$(1)$ holds because the spans, when restricted to $\{K'\supseteq K\}$, are the spans for the deleted data. Meanwhile, $(3)$ holds because every ingredient is compatible with products and tensor products commute with colimits. Therefore we are left to show $(2)$, which requires considering the cases where $\{i\} = A \setminus A^-$ and otherwise, respectively. For $i$ itself, the four maps in (2) are the structure maps into the colimits of Lemmas \ref{lem:Bglue} and \ref{lem:Aglue}. Since $\Phi$ is induced by a map of spans it commutes with these. For another $i'\in A$ we use that blowups and handle attachments at $i$ and $i'$ commute, so that the structure maps for $i'$ are induced, corner by corner, by the same maps one stage earlier: for $A^-$ on the corners $D^-_\bullet$ and $D^-_{\bullet\cup i}$, and for the contracted data on the corner $H_i$. On the $B$-side this uses base change, which holds since the squares are Tor-independent. By induction on $|A|$ and on $m+n$, $\Phi$ intertwines these corner by corner, hence also on colimits.
\end{proof}

We enumerate consequences. First, as a corollary we learn homological mirror symmetry for RFR acyclic cluster varieties themselves.

\begin{cor}
    For any RFR acyclic quiver $Q$, $$\Fuk(\mathcal{X}(Q)) \simeq \Coh(\mathcal{A}(Q)).$$\end{cor} 
    
    Since the construction of $\mathcal{X}(Q)$ converts all conic $\sigma_I^\vee$ loci into disks, we see that $\mathcal{X}(Q)$ is simply a Weinstein manifold (i.e., with no stop). As noted in §1.5, Hu and Shende announce that $\mathcal X(Q)$ is Weinstein deformation equivalent to $\mathcal A(Q)$ \cite{HS}; this is ansatz (*). Granting it, we learn, put into a slogan, that RFR acyclic cluster varieties are \textit{self-mirror}. The remainder of this section establishes that this duality enjoys all of the categorical and combinatorial naturality one could hope for.\\

    To start in this direction we first couch the main theorem microlocally. Since $\mathcal{X}(Q)_{\Sigma,A}$ is stably polarized with core $\mathcal{S} := \mathcal{S}(Q)_{\Sigma,A}$ (Proposition \ref{issa skeleton}), the category $\mu sh_{-\mathcal{S}}(-\mathcal{S})$ of microlocal sheaves on $\mathcal{S}$ is defined, and $\mathrm{Fuk}(\mathcal{X}) \simeq \mu sh_{-\mathcal{S}}(-\mathcal{S})^c$ by Theorem 1.4 of \cite{GPSIII}. Then, given $K \in \Delta_{\Sigma, A}$, we also write $\mathcal{S}_K$ to denote the core of the end $N_K$, so that the obvious inclusions $\mathcal
    S_K \subset S$ render into a functor $\mu sh: \Delta_{\Sigma, A} \rightarrow Cat_\infty$ whose transition maps are microlocal restrictions. The following is then implicit in the constructions, once one applies the main theorems of \cite{GPSIII}, by applying the descent properties of both sides.

\begin{thm}\label{microlocal}
        There is an equivalence of functors $$\Fuk(D_\bullet^\vee) \simeq \mu sh_{-\mathcal{S}_\bullet}(-\mathcal{S}_\bullet)^c$$ \noindent out of $\Delta_{\Sigma, A}$, where the transition functors on the left are right adjoint to sectorial pushforwards and those on the right are microlocalization.
    \end{thm}

    Combining with Theorem \ref{HMS}, this implies that closed restriction is mirror to microlocalization for all cluster models vis-a-vis their boundary divisors. We largely content ourselves to state things about Fukaya categories in the remainder, though the reader is encouraged to reinterpret in terms of microlocal sheaves as they see fit.

    \subsection{Combinatorial Naturalities of Cluster HMS} In this subsection we put the pieces together, proving the appropriate analogue of Theorem \ref{toric HMS natural} for all cluster models. In fact we receive a few extra pieces, to be commented on later, from the fact that cluster HMS is natural under contraction of quivers and from the fact that all cluster models are quasi-affine.\\

    This subsection is essentially one long proof, albeit rife with explanation, which is written to mirror the pieces of Theorem \ref{toric HMS natural}. Therefore, we assume given a cluster model $\A:= \mathcal{A}(Q)_{\Sigma, A}$ and its mirror $\X := \mathcal{X}(Q)_{\Sigma, A}$ throughout, with $S := \mathcal{S}(Q)_{\Sigma, A}$ the skeleton and $\partial S := \partial^\infty \mathcal{S}(Q)_{\Sigma, A}$ the corresponding stop. Also let $\Phi$ be the equivalence of Theorem \ref{HMS}.\\

    \noindent \textit{Freezing is mirror to stop/handle removal.} Assume given $S \subset Q^{mut}$ a set of mutable vertices in $Q$ with corresponding frozen cluster datum $(Q^{\setminus S},\Sigma_{\setminus S},A_{\setminus S})$. Set
$$
\A_S: =\A(Q)_{\Sigma,A}\setminus V(\prod_{e\in S}x_e)= \A(Q^{\backslash S})_{\Sigma_{\backslash S},A_{\backslash S}}\quad\text{and}\quad \X_S:=\X(Q^{\backslash S})_{\Sigma_{\backslash S},A_{\backslash S}}.$$
Write $j_S: \A_S \hookrightarrow \A$ for the induced open immersion of varieties, which gives an open restriction functor $j^*_S: \Coh(\A) \rightarrow \Coh(\A_S)$. We define a prospective mirror functor $\mathfrak sr_S$ as follows. When $S = \{k\}$ contains a single element we let $\mathfrak sr_{\{k\}}\colon \Fuk(\X)\to \Fuk(\X_{\{k\}})$ be either stop removal at $l_k$, if $k \notin A$, or the Verdier quotient along the image of the sectorial pushforward $\Fuk(h_k) \rightarrow \Fuk(\X)$, if $k \in A$. By commuting the latter quotient past the colimit defining $\Fuk(\X)$ from the span $\Fuk(\mathcal{X}(Q)_{\Sigma, A \setminus \{k\}}) \leftarrow \Fuk(F_k) \rightarrow \Fuk(h_k)$, and using the fact that if a stop has a ribbon then stop-removal is precisely the Verdier localization at the image of the sectorial pushforward from a collar of that ribbon, we see in this case that $\mathfrak sr_k$ is actually the stop removal of $l_k$ from $\Fuk(\mathcal{X}(Q)_{\Sigma, A \setminus \{k\}})$.\\

For general $S$, then, we define $\mathfrak sr_S$ to be the composite of the $\mathfrak sr_{\{e\}}$, $e\in S$, in any order. Finally, for $S\subset S'$ we define $\mathfrak sr_{S\subset S'}: \Fuk(\X_S)\to \Fuk(\X_{S'})$ and $j^*_{S\subset S'}: \Coh(\A_S)\rightarrow \Coh(\A_{S'})$ in the same way for the frozen data at $S$. As a last piece of notation we write $\mathcal{P}_Q$ for the poset of sets of mutable vertices in $Q$ ordered by inclusion.
 
\begin{thm}[cf. Theorem \ref{toric HMS natural} (1)]\label{cor:freezing}
The assignments $S\mapsto \Fuk(\X_S)$ and $S \mapsto \Coh(\A_S)$, with transition functors $\mathfrak sr_{S\subset S'}$ and $j^*_{S\subset S'}$, assemble into functors $\mathcal{P}_Q \to Cat_\infty$. Further, there is an equivalence of functors on $\mathcal{P}_Q$
$$
\Phi_\bullet\colon \Fuk(\X_\bullet)\xrightarrow{\ \simeq\ } \Coh(\A_\bullet)$$
whose components are each the appropriate HMS equivalence. In particular $\Phi_S\circ \mathfrak sr_S\simeq j_S^*\circ\Phi$ for every $S$.
\end{thm}

Therefore we learn that restriction onto an open subset of $\A$ corresponding to a frozen subquiver of $Q$ is mirror to stop/handle removal away from the corresponding cells of $\X$.

\begin{proof}
    Both sides are Verdier quotients by certain objects, and so it suffices to match these objects under the equivalence $\Phi$ of Theorem \ref{HMS}. Since Verdier quotients are determined by their kernels, we can reduce to the case of a single vertex removal.\\

    When $k \notin A$, the $B$-side kernel is the image of $\Coh(D_k)$ under the closed pushforward into $\A$ and, since $l_k$ admits a ribbon, the $A$-side kernel is the sectorial image of $\Fuk(N_k)$. By the naturality established in Theorem \ref{HMS} this implies that the two localizations agree. When $k \in A$ the $B$-side kernel is the image of $\Coh(E_k)$ and the $A$-side kernel is the image of $\Fuk(h_k)$, and in this case the same naturality completes the proof.
\end{proof}

On the $B$-side, since categories of coherent sheaves satisfy Zariski descent we see from Proposition \ref{chart lemma} that $$\Coh(\A) \simeq \lim_{\Sigma^{op}} \Coh(\mathcal{A}(Q)_{\Sigma, A}|_{U_{\sigma_I}}).$$

By composing with $\Phi$ and utilizing the previous theorem this implies a corresponding local to global principle for Fukaya categories of mirrors to cluster models. To state it we assume the notation of Corollary \ref{cor:charts}.

\begin{cor}\label{descent}
For any cluster model mirror $\X$, there is an equivalence $$\Fuk(\X) \simeq \lim_{\Sigma^{op}} \Fuk(\mathcal{X}(Q)_{\langle \sigma_I\rangle, I \cap A})$$ \noindent with structure functors the $\mathfrak sr_{S\subset S'}$ of the previous theorem.
\end{cor}

\begin{exmp}
    If $v_i$ and $v_j$ are adjacent, then $\sigma_{\{i,j\}}\notin\Sigma$, so $V(x_i)\cap V(x_j)=\emptyset$ (Corollary 2.29). Hence $\mathcal A=\mathcal A_{\{i\}}\cup\mathcal A_{\{j\}}$, and by Theorem \ref{cor:freezing}

$$\Fuk(\X)\simeq \Fuk(\X_{\{i\}})\times_{\Fuk(\X_{\{i,j\}})} \Fuk(\X_{\{j\}})$$

\noindent along the $\mathfrak{sr}$’s.
\end{exmp}

\noindent \textit{Deletion is mirror to microlocalization.} This piece of our package follows somewhat more directly from the formalism of the previous subsection than the last one. To state it, we return to our usual poset $\Delta_{\Sigma, A}$ of inactive cones in $\Sigma$, with elements $K$, and the notation $D_K$ and $D_K^\vee$ of Theorem \ref{HMS}.

\begin{thm}[cf. Theorem \ref{toric HMS natural} (2)]\label{cluster HMS 2}
There is an equivalence of functors $$\mu sh_{-\mathcal{S}_\bullet}(-\mathcal{S}_\bullet)^c \simeq \Coh(D_\bullet)$$ on $\Delta_{\Sigma, A}$, where the former has structure functors microlocalization and the latter has structure functors restriction.
\end{thm}

In particular we learn that restriction to a closed boundary stratum in $\A$ is mirror to microlocalization onto the corresponding stratum of $\mathcal{S}$.

\begin{proof}
    Combine Theorems \ref{HMS} and \ref{microlocal}.
\end{proof}

\noindent \textit{A generic cotangent fiber maps to the structure sheaf, and accoutrements.} Assume given $x \in T \setminus \pi(\Lambda_\Sigma)$, which is a connected subset of $\mathcal{S}$ since the rays of $\Sigma$ are coordinate rays. The cotangent fibers $T^*_xT$ at any such $x$ are all Hamiltonian isotopic, and so equivalent as objects in $\Fuk(\X)$. We refer to the corresponding such object as the \textit{generic cotangent fiber} of $\X$.
 
\begin{prop}[cf. Theorem \ref{toric HMS natural} (3)]\label{cotangent generation}
For any $x$ as above,
$\Phi(T^*_xT) \simeq \mathcal O_\A$. Hence $$HW^*(T^*_xT,T^*_xT) \simeq H^*(\A, \mathcal{O}_\A).$$ This is natural under stop-removal and sectorial pushforward, in the sense that $1)$ for every set $S$ of mutable vertices, $\mathfrak sr_S(T^*_xT)$ is the generic cotangent fiber of $\X_S$, and
$$\Phi_S\bigl(\mathfrak sr_S(T^*_xT)\bigr)\simeq j_S^*\mathcal O\simeq\mathcal O_{\A_S};$$ $2)$ for $K\in\Delta^+_{\Sigma,A}$ with $T^*_{x_K}T_K$ the generic cotangent fiber of $D_K^\vee$,
$$\Phi\bigl(\iota_K(T^*_{x_K}T_K)\bigr)\simeq(k_K^*)^L\mathcal O_{D_K};$$ and $3)$ for $i\in A$, writing $L_i:=T^*_{x_i}T_i\times C$ for the product of the generic cotangent fiber of $H_i^\vee$ with the cocore $C$ of $T^*D^2$, we have $\Phi(\kappa_*L_i)\simeq\iota_{E*}\mathcal O_E[-1]$.
\end{prop}
 
The object $\iota_K(T^*_{x_K}T_K)=T^*_{x_K}T_K\times(\text{linking arc})^K$ in $2)$ is the linking disk of the stop at a smooth point over $T_K\times\operatorname{int}(\sigma_K)$. In particular, the first clause of this proposition implies that, when $\A$ is just a cluster variety, the wrapped Floer cohomology of the generic cotangent fiber of $\X$ is isomorphic to the ring of global sections on $\A$. 
 
\begin{proof}
By induction on $|A|$, with base case the toric theorem itself. Inducting, let $i \in A$ be the last element of $A$ in the order of handle attachment in Construction \ref{cluster mirrors} (any order will do, by Remark \ref{naturality}) and set $A^-:=A\setminus\{i\}$. Since $T^*_xT$ misses $\mathcal N_i$ it lies in the subsector $\mathcal \X(Q)_{\Sigma,A^-}$, and so $\jmath_*T^*_xT=T^*_xT$. Therefore, by Theorem \ref{HMS} and induction,
$$\Phi(T^*_xT)\simeq\Phi(\jmath_*T^*_xT)\simeq\pi^*\Phi^-(T^*_xT)\simeq\pi^*\mathcal O=\mathcal O.$$

Claims $1)$ and $2)$ then follow from the compatibilities already established, paired with the fact that stop/handle removal clearly preserve generic cotangent fibers. Finally, for $3)$, $(2)$ of Theorem \ref{HMS} gives $\Phi\circ\kappa_*\simeq\iota_{E*}q^*[-1]\circ\Phi_h$, where $\Phi_h=\Phi_{pt} \otimes \Phi''$. The contracted datum $(Q^{/i},\Sigma_{/i},A_{/i})$ has fewer vertices, so by induction on $m+n$ the main statement of this proposition holds for it: $\Phi''(T^*_{x_i}T_i)\simeq\mathcal O_{H_i}$. Also $\Phi_{pt}(C)\simeq\mathcal O_{pt}$ \cite{GL}. Hence $\Phi(\kappa_*L_i)\simeq\iota_{E*}q^*\mathcal O_{H_i}[-1]=\iota_{E*}\mathcal O_E[-1]$.
\end{proof}

We established previously that all cluster models are quasi-affine, which implies that $\mathrm{Perf}(\A) \simeq \Coh(\A)$ is generated by the single object $\mathcal{O}_\A$. Therefore we learn the same of $\Fuk(\X)$: that is, denoting $CW^*$ for the wrapped Lagrangian Floer complex between a pair of Lagrangians, we have

\begin{cor}\label{single object gen}
    $\Fuk(\X) \simeq \mathrm{Perf}(CW^*(T^*_x T, T^*_x T))$.
\end{cor}

\noindent \textit{HMS for the Cluster Boundary.} Because our arguments have trafficked in an understanding of the boundary pieces of both the $A$- and $B$-sides throughout, it requires only a bit more work to also establish a mirror symmetry statement for boundaries of cluster models. In the sequel we write $\partial \A := \bigcup_{i \in \Sigma(1) \setminus A} D_i$, previously denoted $D$, to denote the union of $\A$'s cluster boundary strata and $\partial \X := F_D$ to denote the smoothed symplectic boundary of the convexification of $\X$.\\

We claim that, on either side, these geometric boundaries are expressible as colimits of pieces we understand. On the $B$-side this is particularly easy to see.

\begin{lem}\label{colim B-side}
$\partial \A = \colim_{K \in (\Delta_{\Sigma, A}^+)^{op}} D_K$ in the category of schemes, where the diagram has transition morphisms the closed immersions. Therefore
    $$\Coh(\partial \A) \simeq \colim_{K \in (\Delta_{\Sigma, A}^+)^{op}} \Coh(D_K)$$ \noindent with transition functors the $k_{K \subset K' *}$.
\end{lem}

\begin{proof}
    Proposition \ref{chart lemma} implies that $\partial \A$ is a simple normal crossings divisor, and such divisors are always colimits of the Cech nerve of their components. The result then follows as in Theorem 7.4.1 of \cite{gammage2022mirror}, by utilizing descent properties of IndCoh established in \cite{GR} and then passing to compact objects.
\end{proof}

\begin{rem}
  By the natural transformation whose components are the shifts $[\lvert K\rvert]$, the diagram of the previous lemma is equivalent to the same diagram with transition functors $k_{K\subset K'*}[-| K'\setminus K|]$. But $k_{K \subset K' *}[-|K' \setminus K|] \simeq (k^*_{K \subset K'})^L$. Therefore we learn that $\Coh(\partial \A)$ is also the colimit of the same diagram with transition functors the $(k^*_{K \subset K'})^L$.
\end{rem}

We establish a similar descent statement for $\partial \X$.

\begin{prop}\label{prop:ribbon}
$$\Fuk(\partial \X)\simeq \colim_{K \in (\Delta^+_{\Sigma,A})^{op}}\Fuk(D_K^\vee),$$ with transition functors the $\iota_{K\subset K'}$.
Moreover, the stop of $\X$ is the core of $\partial \X$, i.e. $\partial \X$ is a ribbon of the whole stop, and under the equivalence above the sectorial pushforward functor $\iota_\partial: \Fuk(\partial \X)\to \Fuk(\X(Q)_{\Sigma,A})$ is induced by the $\iota_K$.
\end{prop}
 
\begin{proof}
By Proposition \ref{A-side product splitting} and the paragraph after Definition \ref{collared}, $\partial\mathcal X=F_D$ is covered by slight enlargements of the symplectic faces $F_k\cong D^\vee_k$ (see Example 12.24 of \cite{GPSI}), with $|K|$-fold intersections $D^\vee_K\times T^*\mathbb R^{|K|-1}$ (in standard coordinates). This cover is Weinstein, so the equivalence is Theorem \ref{sectorial} together with K\"{u}nneth.
\end{proof}

Now we have explicitly identified sufficient descent properties on both sides to also prove HMS for boundaries of (mirrors to) cluster models. To state the theorem, we fix one final piece of notation. First, recall the notation $\A^{red}_{\Sigma, A} := \mathcal{A}(Q^{\backslash (\Sigma(1) \setminus A)})_{\Sigma_{\backslash (\Sigma(1) \setminus A)}, A}$ for $\mathcal{A}(Q)_{\Sigma, A} \setminus \partial \mathcal{A}(Q)_{\Sigma, A}$ introduced in Proposition \ref{CY prop}. Mirroring this, write $\X_{\Sigma, A}^{red}$ to denote the Weinstein manifold that arises from removing all stops from $\mathcal{X}(Q)_{\Sigma, A}$. Then we have the following
 
\begin{thm}\label{boundary HMS} The restriction of Theorem \ref{HMS}'s equivalence $\Phi$ to $(\Delta^+_{\Sigma,A})^{op}$ induces an equivalence $$\Phi_\partial\colon \Fuk(\partial \X) \simeq \Coh(\partial \A),$$ satisfying $\Phi\circ\iota_\partial\simeq\partial k_*\circ\Phi_\partial$, where $\partial k:\partial\mathcal A\hookrightarrow\mathcal A$ denotes the boundary inclusion.\\
    
Under this equivalence, the Verdier quotient $\Fuk(\partial \X) \xrightarrow{\iota_\partial} \Fuk(\mathcal{X}(Q)_{\Sigma,A})\xrightarrow{\mathfrak sr_S} \Fuk(\X^{red}_{\Sigma, A})$, which removes the entirety of the stop from $\Fuk(\mathcal{X}(Q)_{\Sigma,A})$ is sent to the Verdier quotient $\Coh(\partial \A) \xrightarrow{\partial k_*} \Coh(\mathcal{A}(Q)_{\Sigma,A})\xrightarrow{j^*} \Coh(\mathcal{A}(Q)_{\Sigma,A}\setminus \partial \A) \simeq \Coh(\A^{red}_{\Sigma, A})$ arising from deleting the boundary from the log Calabi-Yau pair $(\A, \partial \A)$ (see Proposition \ref{CY prop}).
\end{thm}
 
\begin{proof}
The equation and statement about intertwining pushforwards follows from Theorem \ref{HMS} by Lemma \ref{colim B-side} (in the left adjointed form of the remark the following it) and Proposition \ref{prop:ribbon}. Matching the Verdier quotients then follows from the fact that each kernel is generated by the images of the intertwined functors.
\end{proof}

We rephrase this microlocally, leaving the interested reader to deduce the appropriate naturalities.

\begin{cor}[cf. Theorem \ref{toric HMS natural} (4)]
$\mu sh_{-\partial^\infty\mathcal S}(-\partial^\infty\mathcal S)^c \simeq \Coh(\partial \A)$. \end{cor}

\begin{proof}
    By Theorem \ref{cluster HMS 2} and the preceding three results.
\end{proof}

The point of much of the above discussion is that $\partial \A$ behaves like a toric boundary. For toric varieties, the mirror of the boundary is a fiber of the Hori-Vafa superpotential. In fact we expect something similar to be true here, with the torus replaced by the mirror of the interior $\A^{red}_{\Sigma, A} = \A(Q)_{\Sigma,A}\setminus \partial \A$ and the HV monomials by theta functions. Exhibiting this is beyond the scope of the present work.

\subsection{Novelties of Cluster HMS, and Computations.}

Because cluster models are closed under more operations than toric varieties (in particular, affine blowups along certain centers), their mirror symmetry admits further naturalities involving handle attachments. We note the following two here, as a sort of appendix to the previous subsection.

\begin{enumerate}
    \item \textit{Blowups are mirror to handle attachments.} Affine blowup along the non-toric center $H_i=D_i\cap\{x^{c_i}=-1\}$ is mirror to attaching the handle $h_i$, naturally in the boundary strata ($(2)$ of Theorem \ref{HMS}). On the face at $i$, restriction to the slice $\{y_i=-1\}\subset D_i\cong\mathbb C^*_{y_i}\times H_i$ is mirror to capping $S^1_{\chi_i}$ (Lemma \ref{lem:Hleg}).

    \item \textit{Contraction is mirror to passing to the handle.} The center $H_i\cong\mathcal A(Q^{/i})_{\Sigma_{/i},A_{/i}}$ is mirror to the base $H_i^\vee=\mathcal X(Q^{/i})_{\Sigma_{/i},A_{/i}}$ of $h_i=H_i^\vee\times T^*D^2$, and $\Fuk(h_i) \simeq \Fuk(H_i^\vee)$. So $(2)$ of Theorem \ref{HMS} says that $\kappa_*$ is mirror to $\iota_{E*}q^*[-1]\simeq(\iota_E^*)^Lq^*$, i.e. to pulling back to the exceptional divisor and then taking the left adjoint of restriction. Further, the cocore family over the generic fiber of $H_i^\vee$ is mirror to $\iota_{E*}\mathcal O_E[-1]$ ($(3)$ of Proposition \ref{cotangent generation}).
\end{enumerate}

We expect these to be of broader interest. Another interesting feature of cluster HMS is that it is \textit{mutation invariant}. Recall that if $Q$ and $Q'$ are mutation equivalent RFR acyclic quivers, then $\mathcal{A}(Q)\cong\mathcal{A}(Q')$. By applying Theorem \ref{HMS} we learn the following 
\begin{cor}\label{mut invar}
    Fukaya categories of mirrors to cluster varieties are mutation invariant: $$\Fuk(\mathcal{X}(Q)) \simeq \Fuk(\mathcal{X}(Q')),$$ matching generic fibers, for any mutation equivalent RFR acyclic quivers $Q$ and $Q'$.
\end{cor}

We do not know a Weinstein homotopy from $\mathcal{X}(Q)$ to $\mathcal{X}(Q')$.\\

It is worth noting an unsatisfying feature of our HMS formalism, which is that by Proposition \ref{chow vanishing} there are no interesting line bundles on cluster models (and therefore, for example, no correspondence between small pushoffs of cotangent fibers at smooth points of strata and the corresponding Bondal collection line bundles). We expect, however, that a torus equivariant version of the statements in this work will hold, which will restore interesting coherent-constructible geometry involving elements of the torus equivariant Picard group.\\

Similarly, there is no obvious analogue of convolution once one leaves behind $T^* T$. Therefore our statements are not monoidal.\\

We are, however, able to do some computations. In particular, we can compute both the symplectic (co)homology of certain mirrors to cluster models and the $0$-th $K$-theory of the mirror to \textit{any} cluster model.

\begin{cor}
    Assume given a cluster model interior $\A^{red}_{\Sigma, A} := \mathcal \A(Q)_{\Sigma,A}\setminus\partial\mathcal A$ with stopless mirror $\X^{red}_{\Sigma, A}$ (per Theorem \ref{boundary HMS}). Then
    
    $$SH^k(\mathcal \X^{red}_{\Sigma, A})\cong HH^k(\A^{red}_{\Sigma, A})\cong\bigoplus_{p+q=k}H^p\bigl(\A^{red}_{\Sigma, A},\textstyle\bigwedge^qT_{\A^{red}_{\Sigma, A}}\bigr).$$

In degree $0$, the closed–open map $SH^0(\mathcal \X^{red}_{\Sigma, A}) \rightarrow HH^0(\A^{red}_{\Sigma, A})\simeq \Gamma(\A^{red}_{\Sigma, A},\mathcal O)$ is a ring isomorphism.
\end{cor}

\begin{proof}
    $\mathcal \X(Q)^{red}_{\Sigma,A}$ is Weinstein, so the closed–open map is an isomorphism by the main theorem of \cite{GanatraThesis}.
\end{proof}

    When $\A^{red}_{\Sigma, A}$ is just a cluster variety then all higher cohomology on the right hand side vanish and the symplectic cohomology of $\X$ simplifies dramatically.\\

    We end by transporting the computation of the K-theory of cluster models across $\Phi$.

    \begin{cor}
        For any cluster model mirror $\X$, $$K_0(\Fuk(\X)) = \mathbb{Z}.$$
    \end{cor}

    \begin{proof}
        By Proposition \ref{chow vanishing}, transported across the equivalence $\Phi$.
    \end{proof}

\section{Applications, Examples, and Weak $k$-Booleanness}\label{apps}

Throughout this section we assume the ansatz (*) (see \cite{HS}), which asserts that for an acyclic really full rank cluster variety $\mathcal{A}(Q)$ there is a Weinstein homotopy between $\X(Q)$ and the Stein manifold $\mathcal{A}(Q)$. Everything will be true independently of (*) if one prefers to read the theorems as establishing equivalences $\Fuk(\mathcal{X}(Q)) \simeq \Coh(\Pi_{u,w})$ (equivalently, $\mu sh_{-\mathcal{S}}(-\mathcal{S})^c \simeq \Coh(\Pi_{u,w})$) for acyclic $Q$ such that $\mathcal{A}(Q) \cong \Pi_{u,w}$.

\subsection{Acyclic Open Positroid Cells}

Our first application of the main theorem is to dense open positroid cells, i.e. those of the form $\Pi_{k,n}$. Recall that a Grassmannian $Gr(k,n)$ is called \textit{finite-type} if $k = 2, n \in \mathbb{N}$ or $k = 3, n \in \{6, 7, 8\}$. This is equivalent to $\Pi_{k,n}$ being a finite-type cluster variety, meaning it admits only finitely many seeds.

\begin{thm}
    Assume that $\Pi_{k,n}$ is open in a finite-type Grassmannian. Then $$\Fuk(\Pi_{k,n}) \simeq \Coh(\Pi_{k,n}).$$
\end{thm}

An application of work by Castronovo establishes one side of \textit{closed} HMS for $Gr(2, p)$ with $p$ prime and $Gr(3,7)$, namely that $\Fuk(Gr(k,n)) \simeq MF(\Pi_{k,n}, W)$ in such cases (see \cite{castronovo2020fukaya}). Our statement is an open analogue that also applies to the remaining finite-type cases.

\begin{proof}
    By \cite{SSBW}, the cluster structures on $\Pi_{2,n}$, $\Pi_{3,6}$, $\Pi_{3,7}$, and $\Pi_{3,8}$ arise from quivers whose non-frozen subgraphs are of type $A_{n-3}$, $D_4$, $E_6$, and $E_8$, resp. Therefore they are all acyclic, and Theorem \ref{Acyclic HMS} applies.
\end{proof}

On the other hand, most $\Pi_{k,n}$'s are provably not acyclic cluster varieties. For example, no $\Pi_{k,2k+1}$ with $k \geq 5$ is acyclic because of cohomological considerations, as explained in \cite{LSII}.

\subsection{Weak $k$-Boolean Permutations and Braids}

A number of interesting examples to which Theorem \ref{Acyclic HMS} applies become available to us through Constructions \ref{con:DES Construction} and \ref{Pasha Construction}. In particular, the following conditions on permutations and braid words will ensure that the outputs of these constructions produce acyclic cluster varieties.

\begin{defn}
    A permutation $w \in S_n$ or positive braid word $\beta \in Br_n^+$ is called \textit{weak $k$-boolean} if there exists a reduced word for $w$ (resp. $\beta)$ which only uses each elementary transposition $s_1, ..., s_{n-1}$ (resp. $\sigma_1, ..., \sigma_{n-1})$ at most $k$ times. A weak  $1$-boolean permutation is called \textit{boolean}.
\end{defn}

Given a pair $u \leq w$ in $S_n$ with $w$ Grassmannian, we expect that the minimal weak $k$-booleanness of $w$ is a measure of the ``cluster complexity" of $\Pi_{u,w}$. To furnish this idea, in the sequel we show that if $w$ is weak $2$-boolean then $\Pi_{u,w}$ is an acyclic cluster variety.

\begin{rem}
    Gao and H\"anni considered a similar property of permutations in \cite{gao2020boolean}. They define $w \in S_n$ to be \textit{strong $k$-boolean} if \textit{every} reduced word for $w$ only uses each elementary transposition at most $k$ times. Gao and H\"anni call this condition ``$k$-booleanness," but we append the adjective \textit{strong} to clearly distinguish their condition from our (literally) weaker one. A permutation is strong $1$-boolean if and only if it is weak $1$-boolean, but in general strong and weak $k$-booleanness differ.\\
    
    Gao and H\"anni establish that strong $2$-booleanness is determined by pattern avoidance criteria, generalizing the same fact for boolean permutations. However, they also show that this cannot be the case for strong $(\geq 3)$-booleanness: for example, the permutation $436512$ is a subpattern of $4357612$, and the former is not strong $3$-boolean while the latter is.\end{rem}

    We introduce higher weak booleanness for two reasons. One is that our theorems apply to positive braid words with certain weak, rather than strong, booleanness properties (we will see that the difference does not matter for any of the \textit{permutations} we consider).
    The other is that we think weak $k$-booleanness is a genuinely interesting combinatorial condition in its own right -- those who are not so inspired may skip ahead to Lemma 5.6. For example, it turns out that weak $2$-booleanness is a pattern avoidance condition:
    
    \begin{lem}
        Let $w$ be weak $k$-boolean, $k \geq 2$, and let $w'$ be a subpattern of $w$ defined by deleting a single letter $m$ from $w$. Then $w'$ is weak $(2k - 2)$-boolean.
    \end{lem}

    \begin{proof}
       Consider the wiring diagram of a reduced word for $w$ which uses each elementary transposition at most $k$ times, thought of as a braid diagram on $n$ strands (i.e., a weave of $n$ strands connected to the top and bottom of a square, labeled from 1 to $n$ in order of where they meet the top boundary) minus the information of over- and under-crossings. We implicitly label a crossing between a strand occupying column $i$ and one in column $i+1$ by the transposition $s_i$. Then, by assumption, the wiring diagram of $w'$ comes from that of $w$ by removing a single strand, say the one whose top is labeled $m$ (and henceforth simply referred to as $m$). After the removal, a crossing which had $m$ in a column to its left now occupies columns one lower than it did before, while a crossing which had $m$ to its right keeps its columns. Since any two strands cross at most once, the result is again a reduced word.\\

Let $A_i$ denote the number of times $s_i$ appears and $m$ is in a column smaller than $i$, $B_i$ the number of times $s_i$ appears and $m$ participates in the crossing (i.e., $m$ is in column $i$ or $i+1$), and $C_i$ the number of times $s_i$ appears and $m$ is in a column larger than $i+1$. By assumption, $A_i + B_i + C_i \le k$. Further, one sees that the number of times $s_i$ appears in $w'$ is $C_i + A_{i+1}$.\\

If $C_i$ and $A_{i+1}$ are both non-zero, then both $B_i$ and $B_{i+1}$ must be non-zero, because at some point $m$ will have to pass between column $i$ and column $i+2$ (in some order) to contribute to both counts. Therefore $0 \le C_i \le k-1$ and $0 \le A_{i+1} \le k-1$, so $C_i + A_{i+1} \le 2k-2$. If instead one of $C_i$ or $A_{i+1}$ is zero then $C_i + A_{i+1} \le k \le 2k-2$, since $k \ge 2$, so the claim follows.
    \end{proof}
    
    Computational evidence communicated to us by Yibo Gao in the case $k = 3$ indicates that weak $3$-booleanness is likely also a pattern avoidance condition, albeit one with infinitely many forbidden subpatterns. Indeed, a refinement of the argument above (which we omit due to space and relevance to the work writ large) shows that if $w$ admits a reduced expression which uses every elementary transposition but one at most $3$ times, and which uses the remaining transposition only \textit{twice}, then any single letter deletion of $w$ is weak $3$-boolean. Since the $k = 3$ case is the critical one for the failure of strong $k$-boolean pattern avoidance, altogether we are led to make the following
    
    \begin{conj}
        Weak $k$-booleanness is a pattern avoidance condition.
    \end{conj}
    
    Note, for example, that both $436512$ and $4357612$ are \textit{weak} $3$-boolean, so Gao and H\"anni's example does not rule this conjecture out.\\

We now list some necessary lemmata about higher booleanness, and then return to geometry in the next subsection. Our first lemma is clear from the subword characterization of Bruhat order.

\begin{lem}
    Assume that $u \leq w$. If $w$ is a weak $k$-boolean permutation, then so is $u$.
\end{lem}

Let $c$ denote the alternating orientation Coxeter element of $S_{2h}$, i.e. a permutation using each elementary transposition exactly once and satisfying $c^h = w_0$. The existence of such an element implies, in tandem with the previous lemma, the following

\begin{cor}
    Every element of $S_n$ is weak $\ceil{\frac{n}{2}}$-boolean.
\end{cor}

Often we can do much better. For example, the longest $k$-Grassmannian permutation $(n - k + 1, ..., n, 1, ..., n - k)$ has $(s_k \dots s_1)\cdots(s_{n-1}\dots s_{n-k})$ as a reduced expression, and so the previous lemma implies that every $k$-Grassmannian permutation is at most weak $k$-boolean.\footnote{We note that by definition, fully commutative permutations are weak $k$-boolean if and only if they're strong $k$-boolean. In particular this means that weak and strong booleanness are the same for $k$-Grassmannian permutations.} We will utilize the following substantial refinement of this fact repeatedly in the sequel.

\begin{prop}
    Let $w \in S_n$ be a $k$-Grassmannian permutation. Then $w$ is weak (equivalently strong, per the footnote) $i$-boolean if and only if the Young diagram $\lambda$ corresponding to $w$ does not contain an $(i + 1) \times (i + 1)$ square.
\end{prop}

\begin{proof}
    Recall from Construction \ref{**} that one can build a reduced expression for $w$ by filling the boxes of $\lambda$ with elementary transpositions and reading them off in a certain order (as depicted in Figure \ref{Le depiction}). In particular, as illustrated in Figure \ref{Le depiction}, two boxes of $\lambda$ contain the same transposition $s_i$ if and only if they lie on the same antidiagonal of the ambient lattice $\mathbb{Z}^2$. The proposition then follows because $\lambda$ contains an antidiagonal with $i$ boxes if and only if it contains an $i\times i$ square. \end{proof}

\begin{cor}
    Any $k$-Grassmannian permutation $w \in S_n$ is $\min(k, n-k)$-boolean.
\end{cor}

\begin{defn}
    A Schubert, Richardson, positroid, or double Bott--Samelson variety $X_w$, $R_{u,w}$, $\Pi_{u,w}$, or $X(\Delta \beta)$ resp. is said to be \textit{weak $k$-boolean} if $w$ or $\beta$ is weak $k$-boolean. A weak $1$-boolean variety is simply called \textit{boolean}.
\end{defn}

We will largely be interested in weak $2$-boolean positroid and double Bott--Samelson varieties, but let us first say a word on the geometry of booleanness.\\

Boolean permutations have been well studied, and are known to enjoy many interesting Coxeter theoretic properties. For example, a permutation $w \in S_n$ is boolean if and only if the Bruhat interval $[e, w]$ is a boolean lattice, if and only if it is $321$ and $3412$ avoiding. Booleanness is also equivalent to $X_w$ being toric. Note that the pattern avoidance description of boolean Schubert varieties shows that any toric Schubert variety is smooth, since a $321$ is a subpattern of $4231$.\\

More generally, a Richardson variety $R_{u,w}$ is toric with respect to the standard torus action of $(\mathbb{C}^*)^{n-1}$ on $Fl(n)$ if and only if the Bruhat interval $[u,w]$ is a (not necessarily boolean) lattice and $\ell(w) - \ell(u) \le n-1$. Recent work of Gorsky, Kim and Sherman-Bennett (see \cite{GKSB}) removes the dimension bound: $R_{u,w}$ is toric with respect to some torus action if and only if $[u,w]$ is a lattice. Interestingly, they also show that this is equivalent to the open Richardson variety $R^\circ_{u,w}$, equivalently (when $w$ is $k$-Grassmannian) the positroid cell $\Pi_{u,w}$, being a torus.\\

Upon tracing through Construction \ref{Pasha Construction}, one sees readily that a positroid pair $(u, w)$ with associated Le-diagram $D_{u,w}$ having at most one dot in each antidiagonal of the relevant underlying skew-shape has a corresponding quiver $Q$ which has only frozen vertices. Conversely, if any antidiagonal has two or more dots, then $Q$ will have a mutable vertex, and so the associated positroid variety is not a torus. Therefore we see the following

\begin{cor}
    A closed Richardson variety $R_{u,w}$ with $w$ $k$-Grassmannian is toric if and only if the Le diagram $D_{u,w}$ associated to the positroid pair $(u, w)$ has at most one dot in every antidiagonal.
\end{cor}

One sees similarly, albeit tracing now through Construction \ref{con:DES Construction}, that a double Bott-Samelson variety $X(\Delta \beta)$ is a torus if and only if $\beta$ is boolean.

%%Let $X(\Delta \beta)$ be a boolean double Bott--Samelson variety. In this case, Construction $\ref{con:DES Construction}$ shows that the quiver for $X(\Delta \beta)$ is a collection of frozen vertices, so that $X(\Delta \beta) \cong (\mathbb{C}^*)^n$ for some $n$.\\

%%In fact, a cluster variety $\mathcal{A}(Q)$ is a torus if and only if it has a single seed, which is the same as saying that $Q$ has no edges. Tracing through the algorithm in Construction \ref{con:DES Construction}, the assumption that $X(\Delta \beta)$ is a torus therefore implies that $\beta$ is boolean. All in all this implies the following

%%\begin{prop}
    %%A braid word $\beta$ is boolean if and only if $X(\Delta \beta)$ is a torus.
%%\end{prop}

%%\begin{rem}
    %%In particular, the previous proposition shows that a permutation $w$ is boolean if and only if $\Pi_{e,w}$ is a torus. Since the anticanonical divisor of a closed positroid variety $X$ is exactly the complement of the corresponding open positroid cell in $X$, we see in this case that in fact $\Pi_{e,w}$ is the dense open torus in the toric variety $X_w$. We are unaware of a reference for this equivalent characterization of booleanness for $w$.
%%\end{rem}

%%To end this example we note that a Schubert variety is toric if and only if it is Boolean, and that tracing through the appropriate identifications shows that the maximal torus orbit in a Boolean $X_w$ is exactly the locus $\Pi_{e,w} \subset X_w$.\\

\subsection{On HMS for Positroid and Double Bott--Samelson Varieties}

As was alluded to in the previous subsection, we expect that the (minimal) weak $k$-booleanness of a permutation $w$ or braid word $\beta$ is a reasonable measure of the complexity of the positroid or double Bott--Samelson variety $\Pi_{u,w}$ or $X(\Delta \beta)$. For example, while weak 2-boolean double Bott-Samelson varieties are not tori, they are only slightly more complex.

\begin{prop}
    Let $\beta$ be a weak $2$-boolean braid word. Then $X(\Delta\beta) \cong \mathcal{A}(Q_\beta)$, where the underlying undirected graph on the mutable portion of $Q_\beta$ is a disjoint union of paths. In particular, $Q_\beta$ is acyclic.
\end{prop}

\begin{proof}
    We tease this out of the algorithm in Construction \ref{con:DES Construction}. If $\beta$ is weak $2$-boolean, then that algorithm produces a quiver $Q_\beta$ with at most one mutable vertex colored by a given integer $i$. The only time this introduces edges between vertices of different colors is in Step 3), and such edges can only go between vertices with adjacent (i.e., differing by exactly $1$) colors. Therefore, coloring embeds the underlying graph of the mutable portion of $Q_\beta$ into a path, and so it is a disjoint union of paths.
\end{proof}

The preceding proposition does not admit a clear converse. That is, for arbitrarily large $k$ there are braid words $\beta$ which are \textit{not} weak $k$-boolean but such that $X(\Delta \beta)$ has a cluster structure defined by a seed whose quiver is a disjoint union of copies of $A_n$. For example, assume that $\beta = \sigma_{i_1} \cdots \sigma_{i_\ell}$ has no \textit{adjacent transpositions}, i.e. that if $\sigma_i$ appears (possibly with repetition) amongst the $\sigma_{i_j}$ then $\sigma_{i-1}$ and $\sigma_{i+1}$ do not. Then one can see from Construction \ref{con:DES Construction} that $X(\Delta \beta)$ is a product of $A_n$ cluster varieties.\\

It would be interesting to characterize which $X(\Delta \beta)$ are $\mathcal{A}(Q)$'s for $Q$ with mutable portion a disjoint union of paths; as per our expectation at the beginning of this subsection, we expect this class to be not much larger than the set of $X(\Delta \beta)$ where $\beta$ is weak $2$-boolean.

\begin{cor}
    Weak $2$-boolean double Bott--Samelson varieties are self-mirror.
\end{cor}

Let us now consider positroid varieties. We first prove Corollary \ref{Squareless HMS}.

\begin{prop}[Corollary \ref{Squareless HMS}]
    Let $u \leq w$ be a Bruhat comparable pair of permutations in $S_n$, and assume that $w$ is $k$-Grassmannian. Assume further that the Le-diagram corresponding to $u \leq w$ contains no antidiagonals with $3$ or more dotted boxes. Then $\Pi_{u,w}$ is an acyclic cluster variety, and so $$\Fuk(\Pi_{u,w}) \simeq \Coh(\Pi_{u,w}).$$
\end{prop}

\begin{proof}
    Consider the quiver $Q_{u,w}$ for the cluster structure on $\Pi_{u,w}$ resulting from Construction \ref{Pasha Construction}. Color each vertex of $Q_{u,w}$ by the antidiagonal of the underlying Young diagram of the Le-diagram for $u \leq w$ in which it lives, starting from the bottom left with antidiagonal $1$ (e.g., in Figure \ref{Le quiver}, $F_2$ is colored 1 and $F_6$ is colored 5). It follows from the construction that the bottom right vertex of a given color will be frozen in $Q_{u,w}$.\\

    As such, the assumptions ensure that $Q_{u,w}$ has only a single mutable vertex of a given color. Further, the allowed local moves in the construction only produce edges between vertices with adjacent colors. Therefore the same argument as the previous proposition shows that $Q_{u,w}$ is acyclic.
\end{proof}

We record in particular the following slightly weaker but pithy version of the preceding proposition, which says (interpreted diagrammatically) that any positroid cell contained in a Grassmannian Schubert variety $X_\lambda$ satisfying $\lambda \leq (m, m, 2, ..., 2)$ for a thickened hook with arbitrarily long first two rows and arbitrarily many rows of length at most $2$ is self-mirror.

\begin{figure}
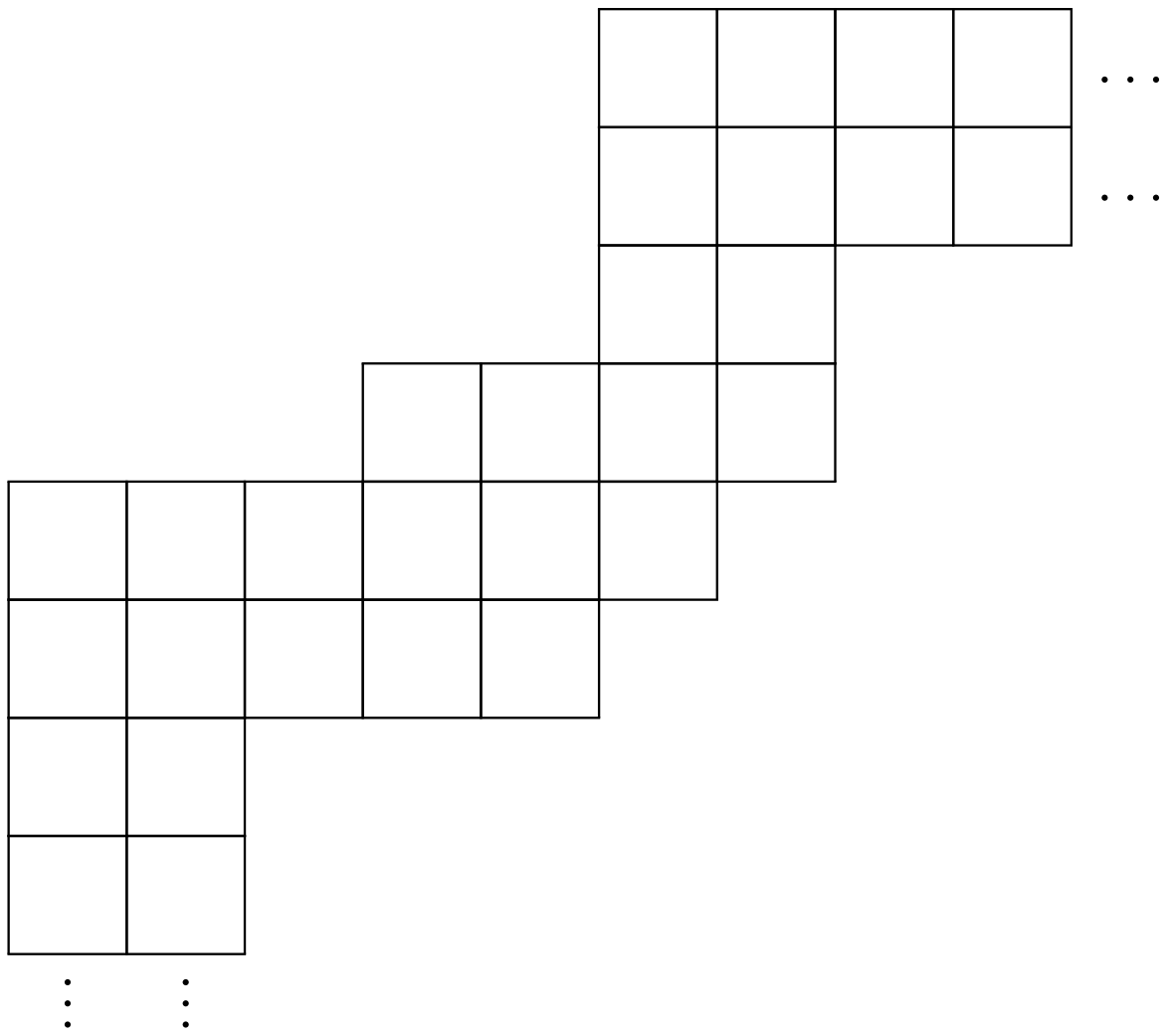

    \centering
    \ytableausetup
{mathmode, boxframe=normal, boxsize=2em}
\begin{ytableau}
\none & \none & \none & \none & \none & & & & & \none[\hdots]\\
\none & \none & \none & \none & \none & & & & & \none[\hdots]\\
\none & \none & \none & \none & \none & &\\
\none & \none & \none & & & &\\
 & & & & & \\
 & & & & \\
 & \\
 & \\
 \none[\vdots] & \none[\vdots]
\end{ytableau}
    
    \caption{\label{BSD} A \textit{border strip diagram} (BSD) is a skew-shape that contains no $2 \times 2$ squares. A \textit{thick BSD} is a skew-shape that contains no $3 \times 3$ squares. Corollary \ref{Skew Squareless} says that if a skew shape $\lambda \setminus \mu$ fits inside of a thick BSD, a generic example of which is depicted in this figure, then $\Pi_{\lambda \setminus \mu}$ is self-mirror. Corollary \ref{2-boolean Schuberts} is the stronger but less general statement that if the (non-skew) Young diagram $\lambda_w$ associated to a permutation $w$ is a thick BSD, then \textit{any} $\Pi_{u,w}$ (whether it is a skew-shaped positroid variety or not) is self-mirror.}
\end{figure}

\begin{cor}\label{2-boolean Schuberts}
     Let $w \in S_n$ be a weak $2$-boolean $k$-Grassmannian permutation, and let $u \leq w$ in Bruhat order. Then $\Pi_{u,w}$ is self-mirror.
\end{cor}

Another convenient specialization of Corollary \ref{Squareless HMS} involves skew-shapes. Let $\mu \leq \lambda$ be a pair of comparable Young diagrams, and let $w_\mu \leq w_\lambda$ be the corresponding positroid pair (notably where $w_\mu$ and $w_\lambda$ are $k$-Grassmannian for the same $k$).  The \textit{skew-shaped positroid variety} $\Pi_{\lambda \setminus \mu}$ associated to $\mu \leq \lambda$ is then defined to be $\Pi_{w_\mu, w_\lambda}$. We note that skew-shaped positroid varieties are exactly the dense open positroid varieties in Grassmannian Richardson varieties. Then, with this notation fixed we have the following

\begin{cor}\label{Skew Squareless}
    Let $\lambda \setminus \mu$ denote a skew-shape, and let $\Pi_{\lambda \setminus \mu}$ denote the associated skew-shaped positroid variety. If $\lambda \setminus \mu$ does not contain a $3 \times 3$ square, then $\Pi_{\lambda \setminus \mu}$ is self-mirror.
\end{cor}

\begin{proof}
    The Le-diagram corresponding to a skew-shaped positroid pair is the skew-shape $\lambda \setminus \mu$ decorated with dots in \textit{every box}. We therefore see that Corollary \ref{Squareless HMS} applies whenever $\lambda \setminus \mu$ has no length $3$ antidiagonals, which is equivalent to the assumptions.
\end{proof}

We now begin an extended computation showing that every positroid cell in a finite-type Grassmannian is self-mirror.

\begin{cor}
    Let $\Pi_{u,w}$ be a positroid cell in $Gr(2, n)$. Then $$\Fuk(\Pi_{u,w}) \simeq \Coh(\Pi_{u,w}).$$ 
\end{cor}

\begin{proof}
    Theorem \ref{projected} implies that $w$ can be chosen to be $2$-Grassmannian, and so Corollary \ref{Squareless HMS} applies.
\end{proof}

Corollary \ref{Squareless HMS} establishes self-HMS for any positroid cell in $Gr(3,6)$ of the form $\Pi_{u,w}$ with $w \leq 356124$ (equivalently, any positroid variety which is contained in the Schubert divisor of $Gr(3,6)$). A computation strengthens this statement into the following

\begin{prop}\label{3,6 acyclic}
    Let $\Pi_{u,w}$ be a positroid cell in $Gr(3,6)$. $\Pi_{u,w}$ is an acyclic cluster variety, and so $$\Fuk(\Pi_{u,w}) \simeq \Coh(\Pi_{u,w}).$$
\end{prop}
    
\begin{proof}
    Since every Young diagram inside of a $3 \times 3$ rectangle save the whole rectangle itself doesn't contain any antidiagonals with $3$ boxes, Corollary \ref{Squareless HMS} implies that it suffices to show that any positroid cell of the form $\Pi_{u, 456123}$ is acyclic.\\
    
    $\Pi_{e, 456123}$ was analyzed in Theorem 5.1. One computes from Construction \ref{Pasha Construction} that the quiver $Q_{s_1, 456123}$ is the top one in Figure \ref{3,6 quivers} (where $F_1, F_2, F_3, F_5, F_8$ are frozen). Mutating at $F_6$ gives an acyclic quiver. Similarly, the quiver $Q_{s_2, 456123}$ is the bottom one in Figure \ref{3,6 quivers} (where $F_1, F_2, F_3, F_5, F_8$ are frozen), which is already acyclic. By symmetry, this handles $Q_{s_4, 456123}$ and $Q_{s_5, 456123}$ as well, and the Le-diagram for $s_3 \leq 456123$ does not contain any antidiagonals with more than $2$ dots and is therefore handled by Corollary \ref{Squareless HMS}.\\

    For any $u$ with length greater than $1$, one can see that the number of non-frozen vertices of $Q_{u, 456123}$ is less than or equal to $2$. But all cluster varieties with at most two mutable vertices are acyclic, and so we're done.

    \begin{figure}
\def\scl{1.3}
\def\Fscl{0.9}
\def\Fsep{4pt}
\def\arrlen{0.5}
\begin{tikzpicture}[scale=\scl]
\def\nodescl{0.7}

\node[circle,scale=0.8,inner sep=0pt,fill=white] (A1) at (0.5,0.5) {$\strut$};
\node[circle,scale=0.8,inner sep=0pt,fill=white] (A2) at (1.5,0.5) {$\strut$};
\node[circle,scale=0.8,inner sep=0pt,fill=white] (A3) at (2.5,0.5) {$\strut$};
\node[circle,scale=0.8,inner sep=0pt,fill=white] (A4) at (1.5,1.5) {$\strut$};
\node[circle,scale=0.8,inner sep=0pt,fill=white] (A5) at (2.5,1.5) {$\strut$};
\node[circle,scale=0.8,inner sep=0pt,fill=white] (A6) at (0.5,2.5) {$\strut$};
\node[circle,scale=0.8,inner sep=0pt,fill=white] (A7) at (1.5,2.5) {$\strut$};
\node[circle,scale=0.8,inner sep=0pt,fill=white] (A8) at (2.5,2.5) {$\strut$};

\node[scale=\Fscl,blue,anchor=north west,inner sep=\Fsep] (F1) at (A1.south east) {$F_1$};
\node[scale=\Fscl,blue,anchor=north west,inner sep=\Fsep] (F2) at (A2.south east) {$F_2$};
\node[scale=\Fscl,blue,anchor=north west,inner sep=\Fsep] (F3) at (A3.south east) {$F_3$};
\node[scale=\Fscl,blue,anchor=north west,inner sep=\Fsep] (F4) at (A4.south east) {$F_4$};
\node[scale=\Fscl,blue,anchor=north west,inner sep=\Fsep] (F5) at (A5.south east) {$F_5$};
\node[scale=\Fscl,blue,anchor=north west,inner sep=\Fsep] (F6) at (A6.south east) {$F_6$};
\node[scale=\Fscl,blue,anchor=north west,inner sep=\Fsep] (F7) at (A7.south east) {$F_7$};
\node[scale=\Fscl,blue,anchor=north west,inner sep=\Fsep] (F8) at (A8.south east) {$F_8$};
\def\qlw{1pt}
\def\qop{1}
\draw[->,blue,opacity=\qop,line width=\qlw] (F4)--(F1);
\draw[->,blue,opacity=\qop,line width=\qlw] (F2)--(F4);
\draw[->,blue,opacity=\qop,line width=\qlw] (F4)--(F3);
\draw[->,blue,opacity=\qop,line width=\qlw,dashed] (F3)--(F2);
\draw[->,blue,opacity=\qop,line width=\qlw] (F5)--(F4);
\draw[->,blue,opacity=\qop,line width=\qlw,dashed] (F3)--(F5);
\draw[->,blue,opacity=\qop,line width=\qlw] (F6)--(F4);
\draw[->,blue,opacity=\qop,line width=\qlw] (F8)--(F7);
\draw[->,blue,opacity=\qop,line 
width=\qlw,dashed] (F5)--(F8);
\draw[->,blue,opacity=\qop,line 
width=\qlw] (F7)--(F5);
\draw[->,blue,opacity=\qop,line 
width=\qlw] (F4)--(F7);
\draw[->,blue,opacity=\qop,line 
width=\qlw] (F7)--(F6);
\end{tikzpicture}

\begin{tikzpicture}[scale=\scl]
\def\nodescl{0.7}

\node[circle,scale=0.8,inner sep=0pt,fill=white] (A1) at (0.5,0.5) {$\strut$};
\node[circle,scale=0.8,inner sep=0pt,fill=white] (A2) at (1.5,0.5) {$\strut$};
\node[circle,scale=0.8,inner sep=0pt,fill=white] (A3) at (2.5,0.5) {$\strut$};
\node[circle,scale=0.8,inner sep=0pt,fill=white] (A4) at (1.5,1.5) {$\strut$};
\node[circle,scale=0.8,inner sep=0pt,fill=white] (A5) at (2.5,1.5) {$\strut$};
\node[circle,scale=0.8,inner sep=0pt,fill=white] (A6) at (0.5,2.5) {$\strut$};
\node[circle,scale=0.8,inner sep=0pt,fill=white] (A7) at (1.5,2.5) {$\strut$};
\node[circle,scale=0.8,inner sep=0pt,fill=white] (A8) at (2.5,2.5) {$\strut$};

\node[scale=\Fscl,blue,anchor=north west,inner sep=\Fsep] (F1) at (A1.south east) {$F_1$};
\node[scale=\Fscl,blue,anchor=north west,inner sep=\Fsep] (F2) at (A2.south east) {$F_2$};
\node[scale=\Fscl,blue,anchor=north west,inner sep=\Fsep] (F3) at (A3.south east) {$F_3$};
\node[scale=\Fscl,blue,anchor=north west,inner sep=\Fsep] (F4) at (A4.south east) {$F_4$};
\node[scale=\Fscl,blue,anchor=north west,inner sep=\Fsep] (F5) at (A5.south east) {$F_5$};
\node[scale=\Fscl,blue,anchor=north west,inner sep=\Fsep] (F6) at (A6.south east) {$F_6$};
\node[scale=\Fscl,blue,anchor=north west,inner sep=\Fsep] (F7) at (A7.south east) {$F_7$};
\node[scale=\Fscl,blue,anchor=north west,inner sep=\Fsep] (F8) at (A8.south east) {$F_8$};
\def\qlw{1pt}
\def\qop{1}
\draw[->,blue,opacity=\qop,line width=\qlw,dashed] (F2)--(F1);
\draw[->,blue,opacity=\qop,line width=\qlw] (F2)--(F4);
\draw[->,blue,opacity=\qop,line width=\qlw] (F4)--(F3);
\draw[->,blue,opacity=\qop,line width=\qlw,dashed] (F3)--(F2);
\draw[->,blue,opacity=\qop,line width=\qlw] (F5)--(F4);
\draw[->,blue,opacity=\qop,line width=\qlw,dashed] (F3)--(F5);
\draw[->,blue,opacity=\qop,line width=\qlw] (F1)--(F6);
\draw[->,blue,opacity=\qop,line width=\qlw] (F8)--(F7);
\draw[->,blue,opacity=\qop,line 
width=\qlw,dashed] (F5)--(F8);
\draw[->,blue,opacity=\qop,line 
width=\qlw] (F7)--(F5);
\draw[->,blue,opacity=\qop,line 
width=\qlw] (F4)--(F7);
\draw[->,blue,opacity=\qop,line 
width=\qlw] (F7)--(F6);
\draw[->,blue,opacity=\qop,line 
width=\qlw] (F6)--(F2);
\end{tikzpicture}
\caption{\label{3,6 quivers} Quivers for certain positroid cells in $Gr(3,6)$, whose acyclicity essentially constitutes a proof of Proposition \ref{3,6 acyclic}.}
\end{figure}
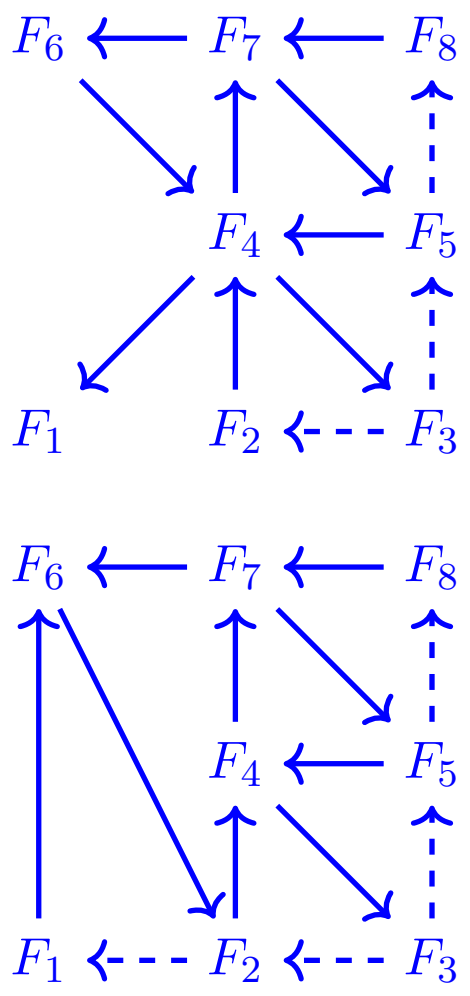
\end{proof}

One could hope that similar computations can be performed for any positroid cell in a finite-type Grassmannian. This is true, but even for $Gr(3,7)$ they are sufficiently cumbersome as to take up untoward amounts of space. We are grateful to Pavel Galashin for access to code computing the quivers $Q_{u,w}$ attached to Le-diagrams. Running it on all positroid cells in $\mathrm{Gr}(3,7)$ and $\mathrm{Gr}(3,8)$ verifies their acyclicity, and so we obtain the following result.

\begin{thm}
    Any positroid cell in a finite-type Grassmannian is self-mirror.
\end{thm}

It is not clear to us, given a positroid variety $\Pi_{u,w}$ which is an acyclic cluster variety, whether any positroid subvariety of it is also an acyclic cluster variety. For this reason we are unable to give a conceptual explanation of this final theorem. We would be quite interested to see one, however.

\bibliographystyle{amsplain}
\bibliography{refs}

@article{STWZ,
  author  = {Shende, Vivek and Treumann, David and Williams, Harold and Zaslow, Eric},
  title   = {Cluster varieties from {L}egendrian knots},
  journal = {Duke Mathematical Journal},
  volume  = {168},
  number  = {15},
  year    = {2019},
  pages   = {2801--2871},
  doi     = {10.1215/00127094-2019-0027},
  eprint  = {1512.08942},
  archivePrefix = {arXiv},
  primaryClass  = {math.SG}
}

@article{FLTZ,
  author  = {Fang, Bohan and Liu, Chiu-Chu Melissa and Treumann, David and Zaslow, Eric},
  title   = {{$T$}-duality and homological mirror symmetry for toric varieties},
  journal = {Advances in Mathematics},
  volume  = {229},
  number  = {3},
  year    = {2012},
  pages   = {1873--1911},
  doi     = {10.1016/j.aim.2011.10.022},
  eprint  = {0811.1228},
  archivePrefix = {arXiv},
  primaryClass  = {math.SG}
}

@article{GPSI,
  author  = {Ganatra, Sheel and Pardon, John and Shende, Vivek},
  title   = {Sectorial descent for wrapped {F}ukaya categories},
  journal = {Journal of the American Mathematical Society},
  volume  = {37},
  number  = {2},
  year    = {2024},
  pages   = {499--635},
  doi     = {10.1090/jams/1012},
  eprint  = {1809.03427},
  archivePrefix = {arXiv},
  primaryClass  = {math.SG}
}

@article{GPSIII,
  author  = {Ganatra, Sheel and Pardon, John and Shende, Vivek},
  title   = {Microlocal {M}orse theory of wrapped {F}ukaya categories},
  journal = {Annals of Mathematics},
  volume  = {199},
  number  = {3},
  year    = {2024},
  pages   = {969--1119},
  doi     = {10.4007/annals.2024.199.3.2},
  eprint  = {1809.08807},
  archivePrefix = {arXiv},
  primaryClass  = {math.SG}
}

@article{GPS1,
  author  = {Ganatra, Sheel and Pardon, John and Shende, Vivek},
  title   = {Covariantly functorial wrapped {F}loer theory on {L}iouville sectors},
  journal = {Publications Math{\'e}matiques de l'IH{\'E}S},
  volume  = {131},
  year    = {2020},
  pages   = {73--200},
  eprint  = {1706.03152},
  archivePrefix = {arXiv},
  primaryClass  = {math.SG}
}

@article{Mul,
  author  = {Muller, Greg},
  title   = {Locally acyclic cluster algebras},
  journal = {Advances in Mathematics},
  volume  = {233},
  year    = {2013},
  pages   = {207--247},
  doi     = {10.1016/j.aim.2012.10.002},
  eprint  = {1109.1881},
  archivePrefix = {arXiv},
  primaryClass  = {math.RA}
}

@article{LSII,
  author  = {Lam, Thomas and Speyer, David E.},
  title   = {Cohomology of cluster varieties. {II}. {A}cyclic case},
  journal = {Journal of the London Mathematical Society},
  volume  = {108},
  number  = {6},
  year    = {2023},
  pages   = {2377--2414},
  eprint  = {2111.12759},
  archivePrefix = {arXiv},
  primaryClass  = {math.AG}
}

@article{Kuwagaki,
  author  = {Kuwagaki, Tatsuki},
  title   = {The nonequivariant coherent-constructible correspondence for toric stacks},
  journal = {Duke Mathematical Journal},
  volume  = {169},
  number  = {11},
  year    = {2020},
  pages   = {2125--2197},
  doi     = {10.1215/00127094-2020-0011},
  eprint  = {1610.03214},
  archivePrefix = {arXiv},
  primaryClass  = {math.AG}
}

@article{GaLa,
  author  = {Galashin, Pavel and Lam, Thomas},
  title   = {Positroid varieties and cluster algebras},
  journal = {Annales Scientifiques de l'{\'E}cole Normale Sup{\'e}rieure (4)},
  volume  = {56},
  number  = {3},
  year    = {2023},
  pages   = {859--884},
  eprint  = {1906.03501},
  archivePrefix = {arXiv},
  primaryClass  = {math.CO}
}

@article{CGGLSS,
  author  = {Casals, Roger and Gorsky, Eugene and Gorsky, Mikhail and Le, Ian and Shen, Linhui and Simental, Jos{\'e}},
  title   = {Cluster structures on braid varieties},
  journal = {Journal of the American Mathematical Society},
  volume  = {38},
  number  = {2},
  year    = {2025},
  pages   = {369--479},
  doi     = {10.1090/jams/1048},
  eprint  = {2207.11607},
  archivePrefix = {arXiv},
  primaryClass  = {math.RT}
}

@article{KLS,
  author  = {Knutson, Allen and Lam, Thomas and Speyer, David E.},
  title   = {Positroid varieties: juggling and geometry},
  journal = {Compositio Mathematica},
  volume  = {149},
  number  = {10},
  year    = {2013},
  pages   = {1710--1752},
  eprint  = {1111.3660},
  archivePrefix = {arXiv},
  primaryClass  = {math.AG}
}

@misc{GLSBSI,
  author  = {Galashin, Pavel and Lam, Thomas and Sherman-Bennett, Melissa and Speyer, David},
  title   = {Braid variety cluster structures, {I}: {3D} plabic graphs},
  note    = {arXiv:2210.04778},
  year    = {2022},
  eprint  = {2210.04778},
  archivePrefix = {arXiv},
  primaryClass  = {math.CO}
}

@article{MR,
  author  = {Marsh, Robert J. and Rietsch, Konstanze},
  title   = {The {B}-model connection and mirror symmetry for {G}rassmannians},
  journal = {Advances in Mathematics},
  volume  = {366},
  year    = {2020},
  pages   = {107027},
  doi     = {10.1016/j.aim.2020.107027},
  eprint  = {1307.1085},
  archivePrefix = {arXiv},
  primaryClass  = {math.AG}
}

@misc{RW,
  author  = {Rietsch, Konstanze and Williams, Lauren},
  title   = {A superpotential for {G}rassmannian {S}chubert varieties},
  note    = {arXiv:2409.00734},
  year    = {2024},
  eprint  = {2409.00734},
  archivePrefix = {arXiv},
  primaryClass  = {math.AG}
}

@article{GHK,
  author  = {Gross, Mark and Hacking, Paul and Keel, Sean},
  title   = {Birational geometry of cluster algebras},
  journal = {Algebraic Geometry},
  volume  = {2},
  number  = {2},
  year    = {2015},
  pages   = {137--175},
  eprint  = {1309.2573},
  archivePrefix = {arXiv},
  primaryClass  = {math.AG}
}

@article{GL,
  author  = {Gammage, Benjamin and Le, Ian},
  title   = {Mirror symmetry for truncated cluster varieties},
  journal = {SIGMA. Symmetry, Integrability and Geometry: Methods and Applications},
  volume  = {18},
  year    = {2022},
  pages   = {Paper 055, 30 pp.},
  doi     = {10.3842/SIGMA.2022.055},
  eprint  = {2103.12232},
  archivePrefix = {arXiv},
  primaryClass  = {math.SG}
}

@article{SSBW,
  author  = {Serhiyenko, Khrystyna and Sherman-Bennett, Melissa and Williams, Lauren},
  title   = {Cluster structures in {S}chubert varieties in the {G}rassmannian},
  journal = {Proceedings of the London Mathematical Society},
  volume  = {119},
  number  = {6},
  year    = {2019},
  pages   = {1694--1744}
}

@article{gao2020boolean,
  author  = {Gao, Yibo and H{\"a}nni, Kaarel},
  title   = {Boolean elements in the {B}ruhat order},
  journal = {Combinatorial Theory},
  volume  = {5},
  number  = {3},
  year    = {2025},
  pages   = {Paper No. 13},
  eprint  = {2007.08490},
  archivePrefix = {arXiv},
  primaryClass  = {math.CO}
}

@article{castronovo2020fukaya,
  author  = {Castronovo, Marco},
  title   = {Fukaya category of {G}rassmannians: rectangles},
  journal = {Advances in Mathematics},
  volume  = {372},
  year    = {2020},
  pages   = {107287}
}

@article{gammage2022mirror,
  author  = {Gammage, Benjamin and Shende, Vivek},
  title   = {Mirror symmetry for very affine hypersurfaces},
  journal = {Acta Mathematica},
  volume  = {229},
  number  = {2},
  year    = {2022},
  pages   = {287--346},
  eprint  = {1707.02959},
  archivePrefix = {arXiv},
  primaryClass  = {math.SG}
}

@article{zhou2019twisted,
  author  = {Zhou, Peng},
  title   = {Twisted polytope sheaves and coherent--constructible correspondence for toric varieties},
  journal = {Selecta Mathematica (N.S.)},
  volume  = {25},
  number  = {1},
  year    = {2019},
  pages   = {Paper No. 1}
}

@misc{castronovo2024cluster,
  author  = {Castronovo, Marco and Gorsky, Mikhail and Simental, Jos{\'e} and Speyer, David E.},
  title   = {Cluster deep loci and mirror symmetry},
  note    = {arXiv:2402.16970},
  year    = {2024},
  eprint  = {2402.16970},
  archivePrefix = {arXiv},
  primaryClass  = {math.AG}
}

@misc{hlavinka2026chow,
  author  = {Hlavinka, Josephine},
  title   = {Chow vanishing and motives of cluster varieties},
  note    = {arXiv:2609.19744},
  year    = {2026},
  eprint  = {2609.19744},
  archivePrefix = {arXiv},
  primaryClass  = {math.AG}
}

@misc{SSWZ,
  author = {Shen, Linhui and Sun, Zhe and Weng, Daping and Zhou, Peng},
  title  = {Homological mirror symmetry for semi-truncated cluster varieties},
  year   = {2026},
  note   = {Preprint}
}

@article{HK,
  author  = {Hacking, Paul and Keating, Ailsa},
  title   = {Homological mirror symmetry for log {C}alabi--{Y}au surfaces},
  journal = {Geometry \& Topology},
  volume  = {26},
  number  = {8},
  year    = {2022},
  pages   = {3747--3833},
  eprint  = {2005.05010},
  archivePrefix = {arXiv},
  primaryClass  = {math.SG}
}

@article{Asplund,
  author  = {Asplund, Johan},
  title   = {Simplicial descent for {C}hekanov--{E}liashberg dg-algebras},
  journal = {Journal of Topology},
  volume  = {16},
  number  = {2},
  year    = {2023},
  pages   = {489--541},
  eprint  = {2112.01915},
  archivePrefix = {arXiv},
  primaryClass  = {math.SG}
}

@phdthesis{GammageThesis,
  author  = {Gammage, Benjamin},
  title   = {Microlocal sheaves and mirror symmetry},
  school  = {University of California, Berkeley},
  year    = {2019}
}

@phdthesis{GanatraThesis,
  author  = {Ganatra, Sheel},
  title   = {Symplectic cohomology and duality for the wrapped {F}ukaya category},
  school  = {Massachusetts Institute of Technology},
  year    = {2012},
  eprint  = {1304.7312},
  archivePrefix = {arXiv},
  primaryClass  = {math.SG}
}

@book{GR,
  author    = {Gaitsgory, Dennis and Rozenblyum, Nick},
  title     = {A study in derived algebraic geometry. {V}ol. {I}. {C}orrespondences and duality},
  series    = {Mathematical Surveys and Monographs},
  volume    = {221},
  publisher = {American Mathematical Society, Providence, RI},
  year      = {2017}
}

@misc{GKSB,
  author  = {Gorsky, Eugene and Kim, Soyeon and Sherman-Bennett, Melissa},
  title   = {Unexpected toric {R}ichardson varieties},
  note    = {arXiv:2603.29260},
  year    = {2026},
  eprint  = {2603.29260},
  archivePrefix = {arXiv},
  primaryClass  = {math.AG}
}

@misc{HS,
  author = {Hu, Mingyuan and Shende, Vivek},
  title  = {On sinks and mirrors},
  year   = {2026},
  note   = {Preprint}
}

\end{document}